\documentclass[letterpaper,9pt,twocolumn,oneside,final,journal]{IEEEtran}
\usepackage{url}
\usepackage{showkeys}
\usepackage{pifont}
\usepackage{amssymb,amsthm,epsfig}
\usepackage[noadjust]{cite}
\usepackage[cmex10]{amsmath}
\usepackage{graphicx}
\usepackage{float}
\usepackage{color}
\usepackage{amsfonts}
\usepackage{multirow}
\usepackage{bm}
\usepackage{caption}
\usepackage{algorithm} 
\usepackage{algorithmic} 
\usepackage{makecell}
\usepackage{xr}
\usepackage{stfloats}

\ifCLASSOPTIONcompsoc
\usepackage[caption=false,font=normalsize,labelfon
t=sf,textfont=sf]{subfig}
\else
\usepackage[caption=false,font=footnotesize]{subfig}
\fi

\theoremstyle{plain}
\newtheorem{theorem}{\hspace{1em}\emph{Theorem}}
   \newtheorem{Lemma}{Lemma}[section]

\theoremstyle{remark}

\theoremstyle{definition}
\newtheorem{Example}{Example}[section]

\DeclareMathOperator{\rank}{rank} \DeclareMathOperator{\diag}{diag}

\DeclareMathOperator{\tr}{tr}

\usepackage{soul}

\begin{document}

\title{Refined Schur Method for Robust Pole Assignment with Repeated Poles}

\author{Zhen-Chen Guo, Jiang Qian, Yun-feng Cai and  Shu-fang Xu
\thanks {Z.C.~Guo, Y.F.~Cai and S.F.~Xu are with LMAM \& School of Mathematical Sciences,
Peking University, Beijing, China, 100871 (e-mail: guozhch06@pku.edu.cn, yfcai@math.pku.edu.cn, xsf@pku.edu.cn).}
\thanks {J.~Qian is with School of Sciences, Beijing University of Posts and Telecommunications, Beijing, China, 100876 (e-mail: jqian104@126.com).}

\thanks{This research was supported in part by NSFC under grant 11301013.}
}
\date{}
\maketitle

\begin{abstract}
Schur-type methods in \cite{Chu2} and \cite{GCQX} solve the robust pole assignment problem
by employing the departure from normality of  the closed-loop system matrix as the measure of robustness.
They work well generally when all poles to be assigned are simple.
However, when some poles are close or even repeated, the eigenvalues of the computed closed-loop system matrix might be inaccurate.
In this paper, we present a refined Schur method, which is able to deal with the case when some or all of the poles to be assigned are repeated.
More importantly, the refined Schur method can still be applied when \verb|place| \cite{KNV} and \verb|robpole|  \cite{Tits} fail to output a solution
when the  multiplicity of some repeated poles is greater than the input freedom.
\end{abstract}


\IEEEpeerreviewmaketitle


\section{Introduction}\label{section1}
\IEEEPARstart{T}{he}  behavior of the state feedback control system in engineering is essentially
determined by the eigen-structure  of the closed-loop system matrix.
Such observation ultimately evokes the arising of the pole assignment problem,
which  can be mathematically  stated as follows.
Denote the dynamic state equation of the time invariant linear system by
\begin{displaymath}
\dot{x}(t)=Ax(t)+Bu(t),
\end{displaymath}
where $A\in \mathbb{R}^{n\times n}$ is the open-loop system matrix and
$B\in \mathbb{R}^{n\times m}$ is  the input matrix.
In control theory, the {\bf State-Feedback Pole Assignment Problem (SFPA)}
is to find a state feedback matrix $F\in\mathbb{R}^{m\times n}$ such that
the eigenvalues of the closed-loop system matrix $A_c=A+BF$,
associated with the closed-loop system
\begin{align*}
\dot{x}(t)=Ax(t)+Bu(t)=(A+BF)x(t)=A_c x(t),
\end{align*}
are the given poles in $\mathfrak{L}=\{\lambda_1, \lambda_2, \ldots, \lambda_n\}$,
which is closed under complex conjugate.
Many valuable contributions have been made to the {\bf SFPA}.
We refer readers to \cite{BD, FO, Ka, MP2, MP3, PM, PCK1, RM, Var, XU} for details.
It is well known that the {\bf SFPA} is solvable for any $\mathfrak{L}$ if and only if $(A,B)$ is controllable \cite{Wonh, XU}.
Through the rest of this paper, {\em we will always assume that $(A,B)$ is controllable}.

When $m>1$, the solution to the {\bf SFPA} is generally not unique.
It then leads to the problem on how to explore the freedom of $F$
such that the closed-loop system achieves some desirable properties.
An important engineering application is to find an appropriate solution $F\in\mathbb{R}^{m\times n}$
to the \textbf{SFPA} such that the eigenvalues of the closed-loop system matrix $A_c=A+BF$ are
as insensitive to perturbations on $A_c$ as possible, which is known as the
{\bf State-Feedback Robust Pole Assignment Problem (SFRPA)}.
%

To solve the {\bf SFRPA}, it is imperative to choose an
appropriate  measure of robustness to characterize the ``insensitivity" quantitatively.
Based on different measures, various methods
\cite{BN, CB, Chu2, Dic, FR, GR, GCQX, KN, KNV, LV, LW, RG, RFBT, SNNP, SAY, SEPB, Tits, Var1, XU}
are put forward. The most attractive methods
might be those given by Kautsky, Nichols, and Van Dooren  \cite{KNV}, where the adopted measures are closely related to the condition number of the eigenvectors matrix of $A_c$.
Method $1$ in \cite{KNV} is implemented as the function \verb|place| in the MATLAB control system toolbox.
Method $0$ in \cite{KNV} may not converge, and then
Tits and Yang \cite{Tits} posed a new approach upon it,
which tends to maximize the absolute value of the determinant of the eigenvectors matrix of $A_c$ and is implemented as the function \verb|robpole| (from SLICOT).
Based on recurrent neural networks, a method  recently is put forward in \cite{LW},
where  many parameters need to be adjusted in order to achieve fast convergence.
Notice that these methods can deal with both simple and repeated poles.
However, they are iterative methods and hence can be expensive.
Moreover, in these methods, the  multiplicity of any repeated pole $\lambda\in\mathfrak{L}$ must not exceed the input freedom $m$. Otherwise, they will fail to give a solution. There exist feasible methods (\cite{RFBT, SNNP}) when the multiplicity of some repeated pole exceeds the input freedom $m$. They also tend to minimize the condition number of the eigenvectors matrix of $A_c$. In both methods,
the real Jordan canonical form of the closed-loop system matrix is employed, and
the size of each Jordan block of the repeated poles is assumed to be known in prior, which is, however, generally hard to obtain. Additionally, both methods could be numerical unstable since the computation of
the Jordan canonical form of a matrix is usually suspected.

Another type of methods uses the departure from normality of $A_c$ as the measure  of robustness. It is firstly proposed as the \verb|SCHUR| method in \cite{Chu2}. Some variations can also be found there. Recently, the authors \cite{GCQX} made some improvements to the methods proposed in \cite{Chu2}, especially for placing complex  conjugate  poles, which is referred to as the \verb|Schur-rob| method. All these Schur-type methods are designed for the case when all poles to be assigned are simple. If some poles are close or even repeated, these methods can still output a solution $F$, but the relative errors of the eigenvalues of the computed closed-loop system matrix $A_c=A+BF$,
compared with the entries in $\mathfrak{L}$, might be fairly large.

%
%

In this paper, we intend to propose a refined version of the
\verb|Schur-rob| method  \cite{GCQX} specifically for repeated poles.
It is well known that a defective eigenvalue,
whose geometric multiplicity is less than its algebraic multiplicity,
is generally more sensitive to perturbations than a semi-simple one,
whose geometric and algebraic multiplicities are identical.
So in the present refined Schur method,
we manage to keep the geometric multiplicities of the repeated poles as large as
possible by constructing the real Schur form  of $A_c$ in more special form,
and then attempt to minimize the departure from normality of $A_c$.
The present refined Schur method can achieve higher relative
accuracy of the placed poles than those Schur-type methods in \cite{GCQX, Chu2} for repeated poles. Moreover,
it still works well when methods in \cite{KNV, Tits} fail
in the case where the  multiplicity of some poles is greater than $m$.
Numerical examples illustrate the superiorities of our approach.


The rest of this paper is organized as follows.
Section \ref{section2} displays some useful preliminaries for solving the {\bf SFRPA}.
Our refined Schur method to assign repeated poles is developed in Section \ref{section3}.
Several illustrative examples are presented in Section \ref{section4} to illustrate the performance of our method. Some concluding remarks are finally drawn in Section \ref{section5}.

\section{Preliminaries and Notations}\label{section2}

We first briefly review the parametric solutions to the {\bf SFPA} \cite{Chu2, GCQX}
using the real Schur decomposition of the closed-loop system matrix $A_c=A+BF$.
Let
\begin{equation}\label{eqrealschur}
A + BF = XTX^{\top}
\end{equation}
be the real Schur decomposition of $A_c$,
where $X\in \mathbb{R}^{n\times n}$ is orthogonal and
$T\in \mathbb{R}^{n\times n}$ is upper quasi-triangular.
Without loss of generality, assume
that {\em $B$ is of  full column rank} and let
$B=Q\begin{bmatrix}R^{\top}& 0\end{bmatrix}^{\top}=\begin{bmatrix}Q_1&Q_2\end{bmatrix}\begin{bmatrix}R^{\top}& 0\end{bmatrix}^{\top}=Q_1R$
be the QR decomposition of $B$, where $Q\in\mathbb{R}^{n\times n}$ is orthogonal, $R\in\mathbb{R}^{m\times m}$ is nonsingular
upper triangular, and $Q_1\in \mathbb{R}^{n\times m}$.
Then  with $X$ and $T$ satisfying
\begin{align}\label{eqXT}
Q_2^{\top}(AX -XT)=0,
\end{align}
the parametric solutions to the {\bf SFPA} can be given by
\[
F=R^{-1}Q_1^{\top}(XTX^{\top}-A).
\]
Consequently, once the orthogonal $X$ and the upper quasi-triangular
$T$ satisfying \eqref{eqXT} are obtained, $F$ will be acquired immediately.

When solving the {\bf SFRPA}, we employ the departure from normality of $A_c$ as the
measure of robustness, which can be specified as (\cite{SSun})
\[
\Delta_F(A_c)=\sqrt{\|A_c\|_F^2-\sum_{j=1}^{n}|\lambda_j|^2},
\]
where $\lambda_j$, $j=1, \dots, n$, are the poles to be placed.
As in \cite{GCQX}, we write $T=D+N$, where $D$ and $N$ are the block diagonal part and the strictly upper quasi-triangular part of $T$, respectively.
Let the $2\times2$ diagonal blocks in $D$ be of the form
$\begin{bmatrix}\begin{smallmatrix}\mbox{Re}(\lambda)&\delta\mbox{Im}(\lambda)\\
-\frac{1}{\delta}\mbox{Im}(\lambda)&\mbox{Re}(\lambda)\end{smallmatrix}\end{bmatrix}$
with $\mbox{Im}(\lambda)\neq 0$, $0\neq \delta\in\mathbb{R}$.
Then $\Delta_F(A_c)$ can be reformulated as
\begin{equation}\label{departure}
\Delta_F(A_c)=\sqrt{\|N\|_F^2+\sum_{\mbox{Im}(\lambda)\neq 0}(\delta-\frac{1}{\delta})^2{\mbox{Im}(\lambda)}^2},
\end{equation}
where the summation is over all $2\times2$ diagonal blocks in $D$. Hence, if some poles to be assigned are non-real, it is not only the corresponding part in $N$ that contributes to $\Delta_F(A_c)$, but also that in $D$.
Our method displayed in the next section is designed to solve the {\bf SFRPA} by finding
some appropriate $X$ and $T$, which  satisfy \eqref{eqXT}, such that the departure from normality of $A_c$,
specified in  \eqref{departure}, is minimized.
Acquiring an optimal solution to $\min\Delta_F(A_c)$ is  rather difficult.
So instead of obtaining a global optimal solution, we prefer to get a suboptimal one with lower computational costs. The matrices $X$ and $T$ satisfying \eqref{eqXT} are computed column by column via solving a series of optimization problems.
Specifically, corresponding to a real pole $\lambda_j$ (the $j$-th diagonal element in $D$), the objective function to be minimized, associated with
$\Delta^2_F(A_c)$, is $\|v_j\|_2^2$, where $\breve{v}_{j}=\begin{bmatrix}v_{j}^\top&0\end{bmatrix}^\top$ with
$v_{j}\in\mathbb{R}^{j-1}$ is the $j$-th column of $N$;
while corresponding to a pair of complex conjugate poles $\lambda_j, \lambda_{j+1}=\bar{\lambda}_j$, it  is
\begin{align}\label{dep.}
\|v_{j}\|_2^2 + \|v_{j+1}\|_2^2 + {\mbox{Im}(\lambda_j)}^2 (\delta -\frac{1}{\delta})^2,
\end{align}
where $\breve{v}_{j+k}=\begin{bmatrix}v_{j+k}^\top&0\end{bmatrix}^\top$ with
$v_{j+k}\in\mathbb{R}^q$, $q\leq j$, are the $(j+k)$-th columns
of $N$ for $k=0,1$, and
$\begin{bmatrix}\begin{smallmatrix}\mbox{Re}(\lambda_{j})&\delta\mbox{Im}(\lambda_{j})\\
-\frac{1}{\delta}\mbox{Im}(\lambda_{j})&\mbox{Re}(\lambda_{j})\end{smallmatrix}\end{bmatrix}$ is
the corresponding $2\times 2$ diagonal block in $D$.

%

The following two lemmas are needed when assigning complex conjugate poles.

\begin{Lemma}\label{Lemma2.1}
Let $A, B\in \mathbb{R}^{n\times n}$ be symmetric, then there exist a diagonal matrix
$\Theta=\diag(\theta_1, \theta_2, \ldots,  \theta_n)$ with $\theta_j\geq 0\,\,(j=1,2,\ldots,n)$
and an orthogonal matrix $U\in\mathbb{R}^{2n\times 2n}$,
whose $j$-th column $u_{j}$  and  $(n+j)$-th column $u_{n+j}$
satisfy $u_{n+j}=\begin{bmatrix}&-I_n\\I_n&\end{bmatrix}u_{j}$, such that
\begin{align}\label{Ham}
\begin{bmatrix}A&B\\B&-A\end{bmatrix}=U\diag(\Theta, -\Theta) U^\top.
\end{align}
Furthermore,  it holds that
$\begin{bmatrix}B&-A\\-A&-B\end{bmatrix}=U\begin{bmatrix}0&-\Theta\\-\Theta&0\end{bmatrix}U^\top$.
\end{Lemma}

Lemma \ref{Lemma2.1} can be verified directly by utilizing properties of Hamiltonian matrices, and we skip the proof here.

\begin{Lemma}(Jacobi Orthogonal Transformation \cite{GCQX})\label{Lemma2.2}
Assume that $x, y \in\mathbb{R}^n$ are linearly independent,
then there exists an orthogonal matrix $Q\in\mathbb{R}^{2\times 2}$,
such that $\tilde{x}^\top \tilde{y}=0$ with
$\begin{bmatrix}\tilde{x}& \tilde{y}\end{bmatrix}=\begin{bmatrix}x& y\end{bmatrix}Q$.
\end{Lemma}

Actually, the $2\times 2$ orthogonal matrix $Q$ in Lemma~\ref{Lemma2.2} can be obtained as follows.
Let $\varrho_1=\|x\|_2^2,\ \varrho_2=\|y\|_2^2, \ \gamma=x^{\top}y$, $\tau=\frac{\varrho_2-\varrho_1}{2\gamma}$ and define $t$ as \begin{displaymath}
t=\left\{ \begin{array}{ll}
1/(\tau+\sqrt{1+\tau^2}), & \text{ if}\quad \tau\geq0,\\
-1/(-\tau+\sqrt{1+\tau^2}), & \text{ if}\quad \tau<0.\\
\end{array} \right.
\end{displaymath}
Then the required $Q$ is $Q=\begin{bmatrix}c & s\\ -s&c\end{bmatrix}$, where $c=1/\sqrt{1+t^2}$ and $s=tc$.

Throughout this paper, we denote the space spanned by the columns of a
matrix $M$ by $\mathcal{R}(M)$,  the null space by $\mathcal{N}(M)$,
and the set of eigenvalues of $M$ by $\lambda(M)$.
The MATLAB expression, which specifies the submatrix with
the colon notation, will be used when necessary,
that is, $M(k:l, s:t)$ refers to the submatrix of $M$ formed
by rows $k$ to $l$ and columns $s$ to $t$. We denote $X=\begin{bmatrix}x_1&x_2&\cdots&x_{n}\end{bmatrix}$ and $X_j=\begin{bmatrix}x_1&\cdots&x_{j}\end{bmatrix}$.
Write the strictly upper quasi-triangular part $N$ of $T$ as
$N=\begin{bmatrix}\breve{v}_1&\breve{v}_2&\cdots&\breve{v}_n\end{bmatrix}$.
For simplicity, we also denote $T(1:j,1:j)$ by $T_j$.


\section{Refined Schur method for repeated poles}\label{section3}
The method in \cite{GCQX} can dispose both simple and repeated poles. However, the repeated eigenvalues of the computed $A_c$,
compared with the entries in $\mathfrak{L}$, might be inaccurate. So this paper is specifically dedicated to repeated poles, both real and non-real. As pointed out in the Introduction part,
a semi-simple eigenvalue is less sensitive to
perturbations than a defective one.
Thus when solving the {\bf SFRPA}, we would keep the geometric multiplicities of repeated poles,
as eigenvalues of $A_c$, as large as possible,
which is actualized by setting  special structure in the upper quasi-triangular matrix $T$ in \eqref{eqrealschur}.

Analogously to  \cite{Chu2, GCQX}, we compute $X$ and $T$ satisfying \eqref{eqXT}
column by column, minimizing corresponding functions associated with $\Delta_F^2(A_c)$ for real poles or complex conjugate poles.
We start with the first pole $\lambda_1$, which is  assumed to be repeated  with  multiplicity $a_1(>1)$, that is, it appears exactly $a_1$ times in $\mathfrak{L}$.

\subsection{Assigning repeated poles $\lambda_1$}\label{subsection3.1}
The strategies vary depending on whether $\lambda_1$ is real or non-real.

\subsubsection{$\lambda_1$ is real}\label{subsection3.1.1}

As an eigenvalue of $A_c=A+BF$, denote its geometric multiplicity by $g_1$.
It then follows that $g_1\le m$ (\cite{KNV}).
If $a_1\le m$, the methods in \cite{KNV,Tits} can be applied, assigning $\lambda_1$ as a semi-simple eigenvalue.
Otherwise, that is $a_1>m$, those methods will fail. In our refined Schur method, if $a_1\le m$, $\lambda_1$ can also
be placed as a semi-simple eigenvalue of $A_c$ with $g_1=a_1$;
if $a_1>m$, $\lambda_1$ can still be assigned with $g_1=m$. Notice that
geometric multiplicity issues are not involved in those Schur-type methods in \cite{Chu2, GCQX}.

Comparing the first $a_1$ columns of \eqref{eqXT} brings
\begin{align}\label{eq2}
Q_2^{\top}AX_{a_1}=Q_2^{\top}X_{a_1}T_{a_1},
\end{align}
where $X_{a_1}=X(:,1:a_1)$ satisfying $X_{a_1}^\top X_{a_1}=I_{a_1}$ and
$T_{a_1}=T(1:a_1, 1:a_1)$ with $\lambda(T_{a_1})=\{ \underbrace{ \lambda_1, \ldots, \lambda_1}_{a_1}\}$
are to be determined. More specifically, to maximize the geometric multiplicity $g_1$, we take $T_{a_1}$ in the special form of
\begin{align}\label{T1form}
\begin{array}{lll}
& \begin{array}{llll}
\quad n_1&\qquad \quad n_2  & \quad \cdots&\quad \ n_l
\end{array}&\\
T_{a_1}=&\begin{bmatrix}D_{11}(\lambda_1)&*&\cdots&*\\ & D_{22}(\lambda_1)&\cdots&*\\&&\ddots&\vdots\\&&&D_{ll}(\lambda_1)\end{bmatrix}&
\begin{array}{l}
n_1\\n_2\\\vdots\\n_l
\end{array}
\end{array}
\end{align}
with $D_{kk}(\lambda_1)=\lambda_1I_{n_k}$, $k=1, \ldots, l$,  $n_1+\dots+n_l=a_1$.
The integers $n_k$,  $k=1, \ldots, l,$ are also to be specified.
Once such $X_{a_1}$ and $T_{a_1}$ satisfying \eqref{eq2} are found,
the geometric multiplicity of $\lambda_1$ will be  no less than $\max\{n_k: \ k=1,  \dots,   l \}$.
So we shall make these $n_k$ as large as possible.
In the following, we show how to set these $n_1, \dots, n_l$
and obtain  the corresponding columns of $X_{a_1}$ and $T_{a_1}$ meanwhile.

Since $D_{11}(\lambda_1)=\lambda_1 I_{n_1}$, by equalling the first $n_1$ columns in both sides of the  equation
in \eqref{eq2} and noticing the orthonormal requirements on columns of $X$, it shows
 that the first $n_1$ columns of $X$ should satisfy
\begin{align}\label{initial}
\begin{array}{l}
M_1\begin{bmatrix}x_1&\cdots&x_{n_1}\end{bmatrix}=0, \\
\begin{bmatrix}x_1&\cdots&x_{n_1}\end{bmatrix}^\top\begin{bmatrix}x_1&\cdots&x_{n_1}\end{bmatrix}=I_{n_1},
\end{array}
\end{align}
where
\begin{align}\label{M_1}
M_1=Q_2^\top(A-\lambda_{1} I_n).
\end{align}
Here, $M_1$ is of full row rank by the controllability of the matrix pencil $(A, B)$,
which implies that $\dim(\mathcal{N}(M_1))=m$.
Let the columns of $S\in\mathbb{R}^{n\times m}$ be an orthonormal basis of $\mathcal{N}(M_1)$. We then display
how to determine $n_1$ and find corresponding $X_{n_1}=\begin{bmatrix}x_1&\cdots&x_{n_1}\end{bmatrix}$ by
distinguishing  two different situations.

\paragraph{Situation \uppercase\expandafter{\romannumeral1} --- $a_1\leq m$ }
In this situation, we set $n_1=a_1$. Then by selecting $x_1, x_2, \ldots, x_{a_1}\in \mathcal{R}(S)$
with $\begin{bmatrix}x_1& x_2& \cdots& x_{a_1}\end{bmatrix}^{\top}\begin{bmatrix}x_1& x_2& \cdots& x_{a_1}\end{bmatrix}=I_{a_1}$,
we have already assigned all $\lambda_1$ and then proceed to the next pole as described in
the next subsection --- Subsection \ref{subsection2.2}.
It is worthwhile to point out that with such choice, the geometric multiplicity $g_1$ of $\lambda_1$ is just $a_1$, that is, $\lambda_1$ is a semi-simple eigenvalue of $A_c$.

\paragraph{Situation \uppercase\expandafter{\romannumeral2} --- $a_1>m$ }
In this situation, we can at most choose $m$ orthonormal vectors from $\mathcal{N}(M_1)$.
So we set $n_1=m$, and then choose $X_{n_1}=SZ$ with $Z\in\mathbb{R}^{m\times m}$ being some orthogonal matrix.

Now assume that we have already obtained $X_q=\begin{bmatrix}x_1&\cdots&x_q\end{bmatrix}$
and $T_q=T(1:q, 1:q)$  with
\begin{align*}
\begin{array}{lll}
& \begin{array}{llll}
\quad n_1&\qquad \quad n_2  & \quad \cdots&\quad \quad n_{k-1}
\end{array}&\\
T_q=&\begin{bmatrix}D_{11}(\lambda_1)&*&\cdots&*\\ & D_{22}(\lambda_1)&\cdots&*\\&&\ddots&\vdots\\&&&D_{k-1,k-1}(\lambda_1)\end{bmatrix}&
\begin{array}{l}
n_1\\n_2\\\vdots\\n_{k-1}
\end{array}
\end{array},
\end{align*}
where $k>1$, $\sum_{j=1}^{k-1}n_j=q$, $n_1=m$ and $D_{jj}(\lambda_1)=\lambda_1I_{n_j}$, $j=1, \ldots, k-1$.
We will show how to determine $n_k$,  the corresponding columns of $X$ and the
corresponding strictly block upper triangular part $T(1:q, q+1:q+n_k)$ in $T$.

From \eqref{eq2} and \eqref{T1form}, the $(q+1)$-th, $\ldots$, $(q+n_k)$-th columns of $X$ and $N$ must satisfy
\begin{align}\label{xv}
\begin{bmatrix}x_{q+j}^\top&v_{q+j}^\top\end{bmatrix}^\top\in\mathcal{N}(M_{q,q}),
\end{align}
where $\breve{v}_{q+j}$, the $(q+j)$-th column of $N$, is $\breve{v}_{q+j}=\begin{bmatrix}v_{q+j}^\top&0\end{bmatrix}^\top$ with $v_{q+j}\in\mathbb{R}^q$ for $j=1,\ldots, n_k$, and
\begin{align}\label{M_m}
M_{q,q}=\begin{bmatrix}Q_2^\top(A-\lambda_{1}I_n)& -Q_2^\top X_q\\X_q^\top&0\end{bmatrix}.
\end{align}
Suppose that the columns of
\begin{align}\label{S1q}
S_{q,q}=\begin{bmatrix}S_{q,q}^{(1)}\\S_{q,q}^{(2)}\end{bmatrix} \quad \text{with }
S_{q,q}^{(1)}\in\mathbb{R}^{n\times m}, \ S_{q,q}^{(2)}\in\mathbb{R}^{q\times m},
\end{align}
form an orthonormal basis of $\mathcal{N}(M_{q,q})$, where $\dim(\mathcal{R}(S_{q,q}))=m$ is guaranteed by Theorem \ref{theorem} in  Subsection \ref{subsection2.3}.   Let $S_{q,q}^{(1)}=U_{q,q}\Sigma_{q,q} V_{q,q}^{\top}=
U_{q,q}\begin{bmatrix}\Sigma_{q,q}^{1}&0\\0&0\end{bmatrix}V_{q,q}^{\top}$
be the Singular Value Decomposition (SVD) of $S_{q,q}^{(1)}$ with $\rank(S_{q,q}^{(1)})=r_q$
and $\Sigma_{q,q}^{1}=\diag(\sigma_{1,q}, \cdots, \sigma_{r_q,q})$, $\sigma_{1,q}\ge\cdots\ge\sigma_{r_q,q}>0$.
Keep in mind that $a_1-q$ is the number of the pole $\lambda_1$ to be assigned,
and $r_q$ is the rank of $S_{q,q}^{(1)}$, which is the maximum number of orthonormal vectors $x_{q+j}$ satisfying \eqref{xv}.
We then need to distinguish whether $a_1-q\leq r_q$ or
not these two cases to discuss how to determine $n_k$ and get
those $x_{q+j}, v_{q+j}, j=1, \ldots, n_k$.
Note that if $r_q=0$, there does not exist nonzero vector $x_{q+j}$ satisfying \eqref{xv},
and hence the method will terminate. Fortunately,  Theorem \ref{theorem} in Subsection \ref{subsection2.3} can assure that $r_q$ is always nonzero.


\begin{itemize}
\item {\bf Case \romannumeral1: $(a_1-q) \leq r_q$.}\quad
In this case, we can set $n_k=a_1-q$, that is, we can assign the remaining $\lambda_1$ together.
From \eqref{xv}, to minimize the departure from normality in \eqref{departure}, it is natural to solve the following optimization problem
\begin{subequations}\label{opt1}
\begin{align}
&\min_{}\|\begin{bmatrix}v_{q+1}& v_{q+2}&\cdots&v_{a_1}\end{bmatrix}\|_F^2 \label{eqreal-opt-equal} \\
\mbox{s.t.}&\left\{ \begin{array}{l}
M_{q,q}\begin{bmatrix}x_{q+1}& x_{q+2}&\cdots&x_{a_1}\\v_{q+1}& v_{q+2}&\cdots&v_{a_1}\end{bmatrix}=0,\\
\begin{bmatrix}x_{q+1}&\cdots&x_{a_1}\end{bmatrix}^{\top}
\begin{bmatrix}x_{q+1}&\cdots&x_{a_1}\end{bmatrix}=I_{a_1-q}, \label{eqreal-opt-constrain}
\end{array}
\right.
\end{align}
\end{subequations}
for $x_{q+1}, \ldots, x_{a_1}$ and $v_{q+1}, \ldots, v_{a_1}$.
By the definition of $S_{q,q}$ we know that
there exists $Z\in\mathbb{R}^{m\times (a_1-q)}$ being of full column rank, such that
\begin{equation}\label{eqrealxv}
\begin{split}
&\begin{bmatrix}x_{q+1}& x_{q+2}&\cdots&x_{a_1}\end{bmatrix}=S_{q,q}^{(1)}Z, \\
&\begin{bmatrix}v_{q+1}& v_{q+2}&\cdots&v_{a_1}\end{bmatrix}=S_{q,q}^{(2)}Z.
\end{split}
\end{equation}
Hence, the optimization problem \eqref{opt1} is equivalent to
\begin{equation}\label{eqreal-opt-equal-2}
\min_{Z^\top S_{q,q}^{(1)\top} S_{q,q}^{(1)}Z=I_{a_1-q}}\tr(Z^\top S_{q,q}^{(2)\top} S_{q,q}^{(2)}Z).
\end{equation}
Let $\hat{Z}=V_{q,q}^{\top}Z$ with
$\hat{Z}=\begin{bmatrix}\hat{Z}_1^\top& \hat{Z}_2^\top\end{bmatrix}^\top$,
$\hat{Z}_1\in \mathbb{R}^{r_q\times {(a_1-q)}}$.
Using $S_{q,q}^{(1)\top}S_{q,q}^{(1)}+S_{q,q}^{(2)\top}S_{q,q}^{(2)}=I_m$,
then the problem \eqref{eqreal-opt-equal-2} is equivalent to
\begin{equation}\label{eqreal-opt-equal-3}
\min_{\hat{Z}_1^\top \Sigma_{q,q}^{1^2}\hat{Z}_1=I_{a_1-q}}\tr(\hat{Z}^{\top}\hat{Z}).
\end{equation}
Write $\tilde{Z}_1=\Sigma_{q,q}^{1}\hat{Z}_1$,
then \eqref{eqreal-opt-equal-3} equals to
\begin{equation}\label{eqreal-opt-equal-4}
\min_{\tilde{Z}_1^{\top}\tilde{Z}_1=I_{a_{1}-q}}\tr(\tilde{Z}_1^{\top}{(\Sigma_{q,q}^{1})}^{-2}\tilde{Z}_1),
\end{equation}
with $\hat{Z}_2=0$. Obviously, the minimum value $\sum_{j=1}^{a_1-q}\frac{1}{\sigma_{j,q}^{2}}$ of \eqref{eqreal-opt-equal-4}
is obtained when $\tilde{Z}_1=\begin{bmatrix}e_1&\cdots&e_{a_1-q}\end{bmatrix}$,
suggesting that \eqref{eqreal-opt-equal-2} achieves its minimum when
\[
Z=V_{q,q}\begin{bmatrix}e_1&\cdots&e_{a_1-q}\end{bmatrix}\diag( \frac{1}{\sigma_{1,q}}, \ldots, \frac{1}{\sigma_{a_1-q,q}}).
\]
Once such $Z$ is obtained, $x_{q+1}, \ldots, x_{a_1}$ and $v_{q+1}, \ldots, v_{a_1}$  can be computed by \eqref{eqrealxv}.
We may then update $X_q$ and $T_q$ as
\begin{equation}\label{updatereal_a_11}
\begin{split}
X_{a_1}&=\begin{bmatrix}X_q &x_{q+1}&x_{q+2}&\cdots&x_{a_1}\end{bmatrix}\in\mathbb{R}^{n\times a_1}, \\
T_{a_1}&=\left[\begin{array}{c|c}
T_q \mathbf & \begin{array}{cccc}
v_{q+1}&v_{q+2}&\cdots&v_{a_1}\\
\end{array}\\
& \\[-2mm]
\hline
& \\[-2mm]
\mathbf &\lambda_{1}I_{a_1-q}\\
\end{array}\right]
\in\mathbb{R}^{a_1\times a_1},
\end{split}
\end{equation}
and proceed with the next pole $\lambda_2$.


\smallskip

\item{\bf Case \romannumeral2: $(a_1-q)> r_q$.} \quad
In this case, we can choose at most $r_q$ orthonormal $x_{q+j}$, $j\geq 1$.
So we set $n_k=r_q$ and let
\begin{align*}
&\begin{bmatrix}x_{q+1}&\cdots&x_{q+r_q}\end{bmatrix}=U_{q,q}(\ :\ ,1:r_q),\\ &\begin{bmatrix}v_{q+1}&\cdots&v_{q+r_q}\end{bmatrix}=S_{q,q}^{(2)}V_{q,q}(\ :\ ,1:r_q){(\Sigma_{q,q}^{1})}^{-1}.
\end{align*}
It can be easily verified that such $x_{q+j},  v_{q+j}$, $j=1, \ldots, r_q$, satisfy \eqref{xv}.
It is worthwhile to point out that in this case we do not need to solve an optimization problem similar to \eqref{opt1} in {\bf Case \romannumeral1},
because the value of the objective function now is a constant when the constraints are satisfied.
We can then update $X_q$ and $T_q$ as 
\begin{equation}\label{updatereal_m_r}
\begin{split}
&X_{q+n_k}=X_{q+r_q}\\
=&\begin{bmatrix}X_q &x_{q+1}&x_{q+2}&\cdots&x_{q+r_q}\end{bmatrix}\in\mathbb{R}^{n\times {(q+r_q)}}, \\
\\
&T_{q+n_k}=T_{q+r_q}\\
=&\left[\begin{array}{c|c}
T_q \mathbf & \begin{array}{cccc}
v_{q+1}&\cdots&v_{q+r_q}\\
\end{array}\\
& \\[-2mm]
\hline
& \\[-2mm]
\mathbf &\lambda_{1}I_{r_q}\\
\end{array}\right]
\in\mathbb{R}^{{(q+r_q)}\times {(q+r_q)}}.
\end{split}
\end{equation}

In this case, some $\lambda_1$ are still unassigned.
We can then pursue a similar process either in {\bf Case \romannumeral1} or
{\bf Case \romannumeral2} until all $\lambda_1$ are placed.
\end{itemize}

Eventually, $T_{a_1}$ being of the form \eqref{T1form} would be acquired.
And this procedure is summarized in  Algorithm \ref{algorithm1}.

\begin{algorithm}
\caption{ \ Assigning real $\lambda_{1}$ }
\renewcommand{\algorithmicrequire}{\textbf{Input:}}
\renewcommand\algorithmicensure {\textbf{Output:} }
\begin{algorithmic}[1]\label{algorithm1}
\REQUIRE ~~\\
$A, Q_2 $, $\lambda_1\in\mathbb{R}$ and $a_1$ (the multiplicity of $\lambda_1$).
\ENSURE ~~\\
Orthogonal $X_{a_1}$ and upper  triangular $T_{a_1}$.
\STATE Find $S\in \mathbb{R}^{n\times m}$,  whose columns are an orthonormal basis of $\mathcal{N}(M_1)$ defined in \eqref{M_1}.
       \IF { $a_1\leq m$}
            \STATE Set $X_{a_1}=SZ$ with $Z\in\mathbb{R}^{m\times a_1}$
                   satisfying $Z^\top Z=I_{a_1}$ and $T_{a_1}=\lambda_1I_{a_1}$.
       \ELSE
            \STATE Set  $X_{a_1}(:,1:m)=S$, $T_{a_1}(1:m,1:m)=\lambda_1I_m$, $q=m$;
            \WHILE{$q<a_1$}
              \STATE Find $S=\begin{bmatrix}S_1\\S_2\end{bmatrix}$ with $S_1\in\mathbb{R}^{n\times m}, S_2\in\mathbb{R}^{q\times m}$,
                          whose columns are an orthonormal basis of $\mathcal{N}(M_{q,q})$ in \eqref{M_m};\label{alg:back_initial}
                   \IF {$(a_1-q)\leq \rank(S_1)$}
                       \STATE Solve the optimization problem \eqref{opt1};
                       \STATE Update $X_{a_1}(:,1:q)$ and $T_{a_1}(1:q,1:q)$ by \eqref{updatereal_a_11}, set $q=a_1$.
                   \ELSE
                       \STATE Update $X_{a_1}(:,1:q)$ and $T_{a_1}(1:q,1:q)$  by \eqref{updatereal_m_r}, set $q=q+\rank(S_1)$.
                  \ENDIF
            \ENDWHILE
       \ENDIF
\end{algorithmic}
\end{algorithm}

\bigskip
\subsubsection{$\lambda_1$ is non-real}\label{subsection3.1.2}
Let $\lambda_1=\alpha_1+i\beta_1$, where $\alpha_1, \beta_1\in\mathbb{R}$ and $\beta_1\neq 0$. As the eigenvalue of $A_c$, its algebraic multiplicity is denoted by $a_1$. Then $\bar{\lambda}_1=\alpha_1-i\beta_1$ is also an eigenvalue of $A_c$ with algebraic multiplicity $a_1$. We are to assign all $a_1$  complex conjugate pairs $\{\lambda_1, \bar{\lambda}_1\}$ in turn,
where the  complex conjugate  poles $\lambda_1$ and  $\bar{\lambda}_1$ are placed simultaneously.

Comparing the first  $2a_1$ columns of \eqref{eqXT} and recalling that $X$ is orthogonal, one can show that $T_{2a_1}$ and $X_{2a_1}$ must satisfy
\begin{align}\label{xt}
Q_2^\top AX_{2a_1}-Q_2^\top X_{2a_1}T_{2a_1}=0, \qquad X_{2a_1}^\top X_{2a_1}=I_{2a_1},
\end{align}
with $\lambda(T_{2a_1})=\{ \underbrace{ \lambda_1, \ldots, \lambda_1}_{a_1},
\underbrace{ \bar{\lambda}_1, \ldots, \bar{\lambda}_1}_{a_1}\}$.
There is  a classical strategy  in \cite{GCQX} to get $T_{2a_1}$ and $X_{2a_1}$ satisfying \eqref{xt}.
Here, the substantial refinement on the strategy in \cite{GCQX} is taking
the geometric multiplicities of $\lambda_1$  and $\bar{\lambda}_1$  into account.
That is,  we would choose $T_{2a_1}$ in a more special form:
\begin{align}\label{T_2a1}
\begin{array}{lll}
& \begin{array}{llll}
\ \ 2n_1& \ \qquad 2n_2&\quad  \cdots&\quad 2n_l
\end{array}&\\
T_{2a_1}=&\begin{bmatrix}D_{11}(\lambda_1)&*&\cdots&*\\ & D_{22}(\lambda_1)&\cdots&*\\&&\ddots&\vdots\\&&&D_{ll}(\lambda_1)\end{bmatrix}&
\begin{array}{l}
2n_1\\2n_2\\\ \vdots\\2n_l
\end{array}
\end{array},
\end{align}
where $D_{kk}(\lambda_1)=\diag(D(\delta_{1,k}(\lambda_1)), \ldots, D(\delta_{n_k,k}(\lambda_1)))$ with
\begin{equation}\label{D_delta1}
\begin{split}
&D(\delta_{p,k}(\lambda_1))\\
=&\begin{bmatrix}\mbox{Re}(\lambda_1)&\delta_{p,k}(\lambda_1)\mbox{Im}(\lambda_1)\\
-\frac{1}{\delta_{p,k}(\lambda_1)}\mbox{Im}(\lambda_1)&\mbox{Re}(\lambda_1)\end{bmatrix},
\quad 0 \neq \delta_{p,k}(\lambda_1)\in\mathbb{R}
\end{split}
\end{equation}
for $p=1, \ldots, n_k$, $k=1, \ldots, l$,  and $\sum_{k=1}^{l}n_k=a_1$.
With such special form of $T_{2a_1}$, the geometric multiplicity $g_1$ of $\lambda_1$ ( and $\bar{\lambda}_1$), as a repeated eigenvalue of $A_c$, is no less than $\max \{n_k: \  k=1, \ldots, l \}$.

Similarly to the case when $\lambda_1$ is real,
we then tend to choose $\max\{n_k: \  k=1, \ldots, l \}$ as large as possible while
computing $T_{2a_1}$ and $X_{2a_1}$ satisfying \eqref{xt}.
However, the placing procedure for the case when $\lambda_1$ is real can not be easily extended to this non-real case. The reason is that for the repeated and non-real poles, it is not only those columns in $N$ that contribute to $\Delta_F(A_c)$, but also those $\delta_{p,k}$ in the diagonal blocks $D(\delta_{p,k}(\lambda_1))$ in $D$, which may differ in each $2\times 2$ blocks of $D$. Let us take the first $2n_1$ columns of $X$ and $T$ as an illustration. Assume that $n_1$ is known
(Indeed, $n_1$ is also a parameter to be determined. We will discuss how to set $n_1$ later.), then to find the first $2n_1$ columns of $X$ and $T$ simultaneously, we need to solve the following optimization problem originated from minimizing $\Delta_F(A_c)$ defined in \eqref{departure}:
\begin{subequations}
\begin{align}
\min_{\delta_{1,1}(\lambda_1), \ldots, \delta_{n_1,1}(\lambda_1)}& \beta_1^2((\delta_{1,1}(\lambda_1)-\frac{1}{\delta_{1,1}(\lambda_1)})^2+\cdots \\
  & \qquad +(\delta_{n_1,1}(\lambda_1)-\frac{1}{\delta_{n_1,1}(\lambda_1)})^2)\\
\mbox{s.t.}\ \qquad &Q^\top_2(AX_{2n_1}-X_{2n_1}D_{11}(\lambda_1))=0,\\
&X_{2n_1}^\top X_{2n_1}=I_{2n_1}.
\end{align}
\end{subequations}
The above optimization problem is fairly difficult to solve.
The associate optimization problems corresponding to other $D_{kk}(\lambda_1)$, $k>1$ are even more ticklish to solve.
Be aware that in the case considered in the above part
when $\lambda_1$ is real, those $\delta_{p,1}(\lambda_1)$ vanish, and we only need to find the
columns of $X$ and $T$ satisfying the two constraints.
Hence, rather than acquiring the columns of $X$ and $T$ corresponding to each $D_{kk}(\lambda_1)$ straightway,
we shall compute those associated with $D(\delta_{p,k}(\lambda_1))$, $p=1, \ldots, n_k$, $k=1, \ldots, l$,  alternately.
That is, in each step, we only compute two more columns of $X$ and $T$ corresponding
to $D(\delta_{p,k}(\lambda_1))$.
Bear in mind that those $n_1, \dots, n_l$ are also to be determined in
the assigning process such that $\max \{n_k: \ k=1, \ldots, l\}$ is as large as possible.

We start with the first two columns of $X$ and $T$.
Comparing the first two columns of \eqref{xt}, we have
\begin{align}
&Q_2^{\top}A\begin{bmatrix}x_1&x_2\end{bmatrix}=Q_2^{\top}\begin{bmatrix}x_1&x_2\end{bmatrix}
\begin{bmatrix}\alpha_1&\delta_{1,1}(\lambda_1)\beta_1\\-\frac{1}{\delta_{1,1}(\lambda_1)}\beta_1&\alpha_1\end{bmatrix}, \label{firstxt}\\
&x_1^\top x_2=0, \ \|x_1\|_2=\|x_2\|_2=1. \label{xt2}
\end{align}
Note that  the corresponding strictly upper quasi-triangular part in $T$ vanishes here,
and the corresponding objective  function \eqref{dep.} now
becomes $\beta_1^2(\delta_{1,1}(\lambda_1) -\frac{1}{\delta_{1,1}(\lambda_1)})^2$. Apparently, it achieves its minimum value $0$ at $\delta_{1,1}(\lambda_1)=1$. We then show how to find $x_1$ and $x_2$ satisfying \eqref{firstxt} and \eqref{xt2} with $\delta_{1,1}(\lambda_1)=1$. Similarly as in \cite{GCQX}, it is equivalent to find $x_1$ and $x_2$ such that
\begin{align}\label{xt1}
Q_2^\top (A-\lambda_1I_n)( x_1+ix_2)=0
\end{align}
with \eqref{xt2} holding.
%

It holds that $\dim(\mathcal{N}(Q_2^{\top}(A-\lambda_1I_n)))=m$  since $(A, B)$ is controllable. Assume that the columns of $S\in\mathbb{C}^{n\times m}$ form an orthonormal basis of $\mathcal{N}(Q_2^{\top}(A-\lambda_1I_n))$.
Define $S_1=\mbox{Re}(S)$, $S_2=\mbox{Im}(S)$. Then \eqref{xt1} implies that $x_1+ix_2=(S_1+iS_2)(y_1+iy_2)$ for some $y_1,y_2\in\mathbb{R}^m$, or equivalently
\begin{align}\label{x1x2}
x_1=S_1y_1-S_2y_2,\quad x_2=S_1y_2+S_2y_1.
\end{align}
If we can choose $y_1$ and $y_2$ to satisfy $x_1^{\top}x_2+x_2^{\top}x_1=0$ and $x_1^{\top}x_1-x_2^{\top}x_2=0$,
then the normalized $x_1$ and $x_2$ will satisfy \eqref{xt2} and \eqref{xt1}. Direct calculations show that
\begin{equation}\label{initial_com}
\begin{split}
x_1^{\top}x_2+x_2^{\top}x_1=\begin{bmatrix}y_1^\top &y_2^\top\end{bmatrix}H_1\begin{bmatrix}y_1^\top&y_2^\top\end{bmatrix}^\top,\\
x_1^{\top}x_1-x_2^{\top}x_2=\begin{bmatrix}y_1^\top &y_2^\top\end{bmatrix}H_2\begin{bmatrix}y_1^\top&y_2^\top\end{bmatrix}^\top,
\end{split}
\end{equation}
with
\begin{align*}
&H_1=\begin{bmatrix}S_1^{\top}S_2+S_2^{\top}S_1&S_1^{\top}S_1-S_2^{\top}S_2\\
S_1^{\top}S_1-S_2^{\top}S_2&-(S_1^{\top}S_2+S_2^{\top}S_1) \end{bmatrix},\\
&H_2=\begin{bmatrix}S_1^{\top}S_1-S_2^{\top}S_2&-(S_1^{\top}S_2+S_2^{\top}S_1)\\
-(S_1^{\top}S_2+S_2^{\top}S_1)&S_2^{\top}S_2-S_1^{\top}S_1 \end{bmatrix}.
\end{align*}
Since $S^{*}S=I_m$, it can be easily verified that $S_1^\top S_2=S_2^\top S_1$ and $S_1^\top S_1+S_2^\top S_2=I_m$.
If $S_1^\top S_2=0$ and $S_1^\top S_1=\frac{1}{2}I_m$, then $x_1^\top x_2=0$ and $\|x_1\|_2=\|x_2\|_2$ for
any $y_1\in\mathbb{R}^m$ and $y_2\in\mathbb{R}^m$ due to \eqref{initial_com}. In this case,
we may arbitrarily choose $y_1$ and $y_2$ with $\|y_1\|_2=\|y_2\|_2=1$, then $x_1$
and $x_2$ computed by \eqref{x1x2} satisfy   
\eqref{xt2} and \eqref{xt1} as required. If $S_1^\top S_2\neq 0$ or $S_1^\top S_1\neq \frac{1}{2}I_m$, then
$\rank(H_1)\geq 1$.
Now by Lemma \ref{Lemma2.1}, assume that
\begin{align*}
H_1=U\diag(\Theta, -\Theta) U^\top, \qquad H_2=U\begin{bmatrix}0&-\Theta\\-\Theta&0\end{bmatrix}U^{\top},
\end{align*}
where $U$ is orthogonal whose $j$-th column $u_{j}$  and $(m+j)$-th column $u_{m+j}$ satisfy $u_{m+j}=\begin{bmatrix}&-I_m\\I_m&\end{bmatrix}u_{j}$, $j=1, \ldots, m$, and
$\Theta=\diag(\theta_1, \theta_2, \ldots, \theta_m)$ with $\theta_j\geq 0$, $j=1, \ldots, m$ and $\theta_1>0$.
Then with
\begin{align}\label{gammazeta}
\begin{bmatrix}y^\top_1& y^\top_2\end{bmatrix}^\top=U\begin{bmatrix}\mu&1&0&\cdots&0 &-\mu&1&0&\cdots&0 \end{bmatrix}^{\top},
\end{align}
where $\mu=\sqrt{\theta_2/\theta_1}$, one can show that $x_1$ and $x_2$ computed by \eqref{x1x2} satisfy $x_1^\top x_2=0$ and $\|x_1\|_2=\|x_2\|_2$. Thus the normalized $x_1$ and $x_2$, i.e. $x_1\triangleq x_1/\|x_1\|_2,\, x_2\triangleq x_2/\|x_2\|_2$, are
the vectors  desired. Overall,  we can obtain $X_2=\begin{bmatrix}x_1 &x_2\end{bmatrix}$ and $T_2=D(\delta_{1,1}(\lambda_1))=D_0(\lambda_1)\triangleq\begin{bmatrix}\alpha_1&\beta_1\\
-\beta_1&\alpha_1\end{bmatrix}$ in either case.

\smallskip

Now assume that the first $2q$ ($1\le q<a_1$) columns of $X$ and $T$ have already been obtained with
\begin{align}\label{QAX-2q}
Q_2^\top AX_{2q}=Q_2^{\top}X_{2q}T_{2q}&, \quad X_{2q}^\top X_{2q}=I_{2q},
\end{align}
we are to find the subsequent $(2q+1)$-th and $(2q+2)$-th columns of $X$ and $T$. Here $T_{2q}$ is of the form   similar as \eqref{T_2a1}:
\begin{equation}\label{T_2q}
\begin{split}
 &T_{2q}\\
=&\begin{array}{ll}
  \begin{array}{llll}\ \ \ 2n_1& \ \ \ \   \cdots &\,\,\, \ \ 2n_{k-1} &\ \ \ \qquad 2n_k\end{array}& \\
     \begin{bmatrix}D_{11}(\lambda_1)&\ \cdots&*&*\\ & \ddots&\vdots&\vdots\\&&D_{k-1,k-1}(\lambda_1)&*\\&&&D_{kk}(\lambda_1)\end{bmatrix}&
     \begin{array}{l}2n_1\\\,\,\, \vdots\\2n_{k-1}\\2n_k\end{array}
\end{array},
\end{split}
\end{equation}
where $D_{11}(\lambda_1), \dots, D_{kk}(\lambda_1)$ are block diagonal with  $2\times 2$ matrices being of the form  \eqref{D_delta1} as the diagonal blocks and $n_1+\dots+n_k=q$. Notice that $n_1,\dots,n_{k-1}$ have already been determined, while $n_k$ might still be updated when  computing the $(2q+1)$-th and $(2q+2)$-th columns of $X$ and $T$.
More specifically, denote
\[
T_p=\begin{bmatrix}D_{11}(\lambda_1) & \ \cdots &*\\&\ddots &\vdots\\&&D_{k-1,k-1}(\lambda_1)\end{bmatrix}
\]
with $p=2n_1+\dots+2n_{k-1}$ and let
$D_{kk}(\lambda_1)$$=$$\diag$$(D(\delta_{1,k}(\lambda_1))$, $\ldots$, $D(\delta_{j,k}(\lambda_1) ))$, then
the resulted $T_{2q+2}$ could be in the form of
\begin{equation}\label{updatereal_a_1}
\begin{split}
&T_{2q+2}=\left[\begin{array}{c|c|c}
T_p \mathbf & *   \mathbf &
\begin{array}{cc}
v_{2q+1}&v_{2q+2}\\
\end{array}  \\
& \\[-2mm]
\hline
& \\[-2mm]
\mathbf &  D_{kk}(\lambda_1) \mathbf & 0\\
& \\[-2mm]
\hline
& \\[-2mm]
\mathbf &  \mathbf & D(\delta_{j+1, k}(\lambda_1))\\
\end{array}\right],\\
&v_{2q+1}, \ v_{2q+2}\in\mathbb{R}^{p},
\end{split}
\end{equation}
or in the form of
\begin{equation}\label{updatereal_a_2}
T_{2q+2}=\left[\begin{array}{c|c}
T_{2q} \mathbf & \begin{array}{cc}
v_{2q+1}&v_{2q+2}\\
\end{array}\\
& \\[-2mm]
\hline
& \\[-2mm]
\mathbf & D(\delta_{1,k+1}(\lambda_1))\\
\end{array}\right],\quad v_{2q+1},v_{2q+2}\in\mathbb{R}^{2q}.
\end{equation}

If $T_{2q+2}$ is in the form of \eqref{updatereal_a_1}, $n_k$ will be increased by $1$, meaning that
$n_k$ would be updated as $n_k\triangleq n_k+1$; while if $T_{2q+2}$ is in the form of \eqref{updatereal_a_2}, $n_k$ is fixed and $n_{k+1}$ is initially set to be 1. Taking the geometric multiplicity $g_1$ of
$\lambda_1$ (and $\bar{\lambda}_1$)
into account, we incline to  make $n_k$ as large as possible,
suggesting that we would prefer $T_{2q+2}$ in the form of \eqref{updatereal_a_1} whenever possible.

We now turn to show how to determine whether \eqref{updatereal_a_1} is possible and how to find
the $(2q+1)$-th and  $(2q+2)$-th columns of $X$ and $T$ accordingly.
Provided that $T_{2q+2}$ is in the form of \eqref{updatereal_a_1},
then by comparing the $(2q+1)$-th and $(2q+2)$-th columns of \eqref{xt} and noting that $X$ is orthogonal, we have
\begin{align}\label{eqxq}
\left\{\begin{array}{l}
Q_2^\top (A \begin{bmatrix}x_{2q+1}& x_{2q+2}\end{bmatrix}
-X_{p} \begin{bmatrix}v_{2q+1} &v_{2q+2}\end{bmatrix} \\
 \qquad - \begin{bmatrix}x_{2q+1}& x_{2q+2}\end{bmatrix}D(\delta_{j+1,k}(\lambda_1)))=0,\\
X_{2q}^\top \begin{bmatrix}x_{2q+1}& x_{2q+2}\end{bmatrix}=0,\\
\begin{bmatrix}x_{2q+1}& x_{2q+2}\end{bmatrix}^\top \begin{bmatrix}x_{2q+1}& x_{2q+2}\end{bmatrix}=I_2.
\end{array}
\right.
\end{align}
Our goal now is to minimize \eqref{dep.} subject to \eqref{eqxq}. By writing
$\delta_{j+1,k}(\lambda_1)=\frac{\delta_2}{\delta_1}$ with $0\neq\delta_1\in\mathbb{R}$ and $\delta_2\in\mathbb{R}$,
it follows from \cite{GCQX} that the restriction \eqref{eqxq} is equivalent to
\begin{align}\label{eqxq2}
\left\{\begin{array}{l}
M_{2q, p}\begin{bmatrix}\tilde{x}_{2q+1}+i\tilde{x}_{2q+2}\\ \tilde{v}_{2q+1}+i\tilde{v}_{2q+2}\end{bmatrix}=0,\\
\begin{bmatrix}\tilde{x}_{2q+1}& \tilde{x}_{2q+2}\end{bmatrix}^\top \begin{bmatrix}\tilde{x}_{2q+1}& \tilde{x}_{2q+2}\end{bmatrix}=
\diag( 1/\delta_1^2, 1/\delta_2^2),\\
x_{2q+1}=\delta_1\tilde{x}_{2q+1}, \  x_{2q+2}=\delta_2\tilde{x}_{2q+2}, \\
v_{2q+1}=\delta_1\tilde{v}_{2q+1}, \  v_{2q+2}=\delta_2\tilde{v}_{2q+2},
\end{array}
\right.
\end{align}
where
\begin{align}\label{M_2q^p}
M_{2q, p}=\begin{bmatrix}Q_2^\top(A-\lambda_1I_n)& -Q_2^\top X_p\\ X_{2q}^\top&0\end{bmatrix}.
\end{align}
Let the columns of
\[
S_{2q,p}=\begin{bmatrix}S_{2q,p}^{(1)}\\S_{2q,p}^{(2)}\end{bmatrix}
\begin{array}{l}n\\p\end{array}
\]
be an orthonormal basis of $\mathcal{N}(M_{2q,p})$. We shall distinguish three cases  upon $\dim(\mathcal{R}(S_{2q,p}^{(1)}))$ to reveal the assigning process, i.e., to  compute $x_{2q+1},x_{2q+2},v_{2q+1}$ and $v_{2q+2}$  such that \eqref{dep.} is optimized.

\begin{itemize}
\item {\bf Case \romannumeral3: $\dim(\mathcal{R}(S_{2q,p}^{(1)}))\geq 2$.}\quad
Let $S_{2q,p}^{(1)}=U_{2q,p}\Sigma_{2q,p}V_{2q,p}^{*}$ be the SVD of $S_{2q,p}^{(1)}$ with $\sigma_1$, $\sigma_2$ being the first two largest singular values of $S_{2q,p}^{(1)}$ and let $\tilde{x}_1=\mbox{Re}(U_{2q,p}e_1)$, $\tilde{y}_1=\mbox{Im}(U_{2q,p}e_1)$.
If $\tilde{x}_1^\top\tilde{y}_1=0$ and
$\|\tilde{x}_1\|_2=\|\tilde{y}_1\|_2=\frac{\sqrt{2}}{2}$,
we take
\begin{align*}
&x_{2q+1}=\sqrt{2}\tilde{x}_1, \quad v_{2q+1}=\sqrt{2}\mbox{Re}(S_{2q,p}^{(2)}V_{2q,p}e_1/\sigma_1),  \\
&x_{2q+2}=\sqrt{2}\tilde{y}_1, \quad
v_{2q+2}=\sqrt{2}\mbox{Im}(S_{2q,p}^{(2)}V_{2q,p}e_1/\sigma_1).
\end{align*}
With such choice, \eqref{eqxq} is satisfied with $\delta_{j+1,k}(\lambda_1)=1$,
which results in the third term in the function defined in \eqref{dep.} vanishing  and  the
first two terms achieving  $2\frac{1-\sigma_1^2}{\sigma_1^2}$, a value that is a comparable
multiple (less that $2$) of its minimum (Please refer to \cite{GCQX} for details.).
Otherwise, that is  $\tilde{x}_1^\top\tilde{y}_1\neq0$ or
$\|\tilde{x}_1\|_2\neq\|\tilde{y}_1\|_2$, the suboptimal technique for assigning complex conjugate poles in \cite{GCQX} is applied. Specifically, denote
$\tilde{x}_2=\mbox{Re}(U_{2q,p}e_2)$, $\tilde{y}_2=\mbox{Im}(U_{2q,p}e_2)$
and define $\tilde{X}_{2q,p}=\begin{bmatrix}\tilde{x}_1&\tilde{x}_2\end{bmatrix}$,
$\tilde{Y}_{2q,p}=\begin{bmatrix}\tilde{y}_1&\tilde{y}_2\end{bmatrix}$,
$w_1={S_{2q,p}^{(2)}V_{2q,p}e_1}/{\sigma_1}$, $w_2={S_{2q,p}^{(2)}V_{2q,p}e_2}/{\sigma_2}$, then we set
\begin{align*}
x_{2q+1}&=\begin{bmatrix}\tilde{X}_{2q,p}&-\tilde{Y}_{2q,p}\end{bmatrix}
\begin{bmatrix}\gamma_1&\gamma_2&\zeta_1& \zeta_2\end{bmatrix}^{\top}, \\
x_{2q+2}&=\begin{bmatrix}\tilde{Y}_{2q,p}&\tilde{X}_{2q,p}\end{bmatrix}
\begin{bmatrix}\gamma_1&\gamma_2&\zeta_1& \zeta_2\end{bmatrix}^{\top},\\
v_{2q+1}&=[\begin{smallmatrix}\mbox{Re}(w_1)& \mbox{Re}(w_2)&-\mbox{Im}(w_1)& -\mbox{Im}(w_2)\end{smallmatrix}]
\begin{bmatrix}\gamma_1&\gamma_2&\zeta_1& \zeta_2\end{bmatrix}^{\top},\\
v_{2q+2}&=[\begin{smallmatrix}\mbox{Im}(w_1)& \mbox{Im}(w_2)&\mbox{Re}(w_1)& \mbox{Re}(w_2)\end{smallmatrix}]
\begin{bmatrix}\gamma_1&\gamma_2&\zeta_1& \zeta_2\end{bmatrix}^{\top},
\end{align*}
where $\begin{bmatrix}\gamma_1&\gamma_2&\zeta_1& \zeta_2\end{bmatrix}^{\top}\in\mathbb{R}^4$
is to be  chosen such that the function defined in \eqref{dep.}
is optimized in some sense. We refer readers to \cite{GCQX} for more details on this suboptimal technique.
Overall, the resulted $T_{2q+2}$ will be in the form of \eqref{updatereal_a_1} in this case.

\item{\bf Case \romannumeral4: $\dim(\mathcal{R}(S_{2q,p}^{(1)}))=1$ and $\mbox{Re}(u),\mbox{Im}(u)$ are linearly independent.} \quad
    Here $u$ is the left singular vector of $S_{2q,p}^{(1)}$ corresponding to its unique nonzero singular value $\sigma_1$. In this case, suppose that $S_{2q,p}^{(1)}\in\mathbb{R}^{n\times r}$, and let $V_{2q,p}\in\mathbb{R}^{r\times r}$ be the right singular vectors matrix of $S_{2q,p}^{(1)}$. Define $\mathcal{N}_1(M_{2q,p})=\{ \begin{bmatrix}u^\top & w^\top \end{bmatrix}^\top: \ w=S_{2q,p}^{(2)}V_{2q,p}\begin{bmatrix}\frac{1}{\sigma_1} & \eta_2& \cdots & \eta_r\end{bmatrix}^\top, \ \eta_2, \ldots, \eta_r\in\mathbb{C}\}$,  then in the sense of nonzero scaling, $\mathcal{N}_1(M_{2q,p})$ is the unique  subset of $\mathcal{N}(M_{2q,p})$ satisfying $z\in \mathbb{C}^n$, $w\in \mathbb{C}^p$,  $z\neq0$ with $\begin{bmatrix}z^\top&w^\top\end{bmatrix}^\top\in \mathcal{N}(M_{2q,p})$.
    Write $u=\mbox{Re}(u)+i \mbox{Im}(u)\in\mathbb{C}^n$, $w=\mbox{Re}(w)+i \mbox{Im}(w)\in\mathbb{C}^{p}$,
    then we have that $\mbox{Re}(u)$, $\mbox{Im}(u)$, $\mbox{Re}(w)$ and $\mbox{Im}(w)$ satisfy
\begin{align*}
\left\{\begin{array}{l}
Q_2^\top (A \begin{bmatrix}\mbox{Re}(u)& \mbox{Im}(u)\end{bmatrix}
-X_{p} \begin{bmatrix}\mbox{Re}(w) &\mbox{Im}(w)\end{bmatrix} \\
\qquad - \begin{bmatrix}\mbox{Re}(u)& \mbox{Im}(u)\end{bmatrix}D_0(\lambda_1))=0,\\
X_{2q}^\top \begin{bmatrix}\mbox{Re}(u)& \mbox{Im}(u)\end{bmatrix}=0,
\end{array}\right.
\end{align*}
and $\|w\|_2^2=\frac{1-\sigma_1^2}{\sigma_1^2}+|\eta_2|^2+\ldots+|\eta_r|^2$.

Since $\mbox{Re}(u)$ and $\mbox{Im}(u)$ are linearly independent, we shall pursue the Jacobi orthogonal transformation in Lemma \ref{Lemma2.2} on them, i.e., $\begin{bmatrix}\tilde{x}_{2q+1}&\tilde{x}_{2q+2}\end{bmatrix}=  \begin{bmatrix}\mbox{Re}(u)&\mbox{Im}(u)\end{bmatrix}\begin{bmatrix}c&s\\-s&c\end{bmatrix}$,
and set $x_{2q+1}, x_{2q+2}$   be  the normalized vectors of $\tilde{x}_{2q+1}, \tilde{x}_{2q+2}$, respectively.
Accordingly, $v_{2q+1}, v_{2q+2}$ are defined as
\begin{equation} \label{v1v2}
\begin{split}
&\begin{bmatrix}v_{2q+1}& v_{2q+2}\end{bmatrix}\\
=&\begin{bmatrix}\mbox{Re}(w)&\mbox{Im}(w)\end{bmatrix}
\begin{bmatrix}c&s\\-s&c\end{bmatrix} \begin{bmatrix}\frac{1}{\|\tilde{x}_{2q+1}\|_2}&\\&\frac{1}{\|\tilde{x}_{2q+2}\|_2}\end{bmatrix}.
\end{split}
\end{equation}
It is worthwhile to stress again that now we have $\breve{v}_{2q+s}=\begin{bmatrix}v_{2q+s}^\top&0\end{bmatrix}^\top$, $v_{2q+s}\in\mathbb{R}^{p}$ for  $s=1,  2$.
Be aware that $w$ is unknown here since those values
$\eta_2, \ldots, \eta_r \in\mathbb{C}$ have not been specified.
Notice that $ D(\delta_{j+1, k}(\lambda_1))$  has already been determined
with $\delta_{j+1,k}(\lambda_1)=\frac{\|\tilde{x}_{2q+1}\|_2}{\|\tilde{x}_{2q+2}\|_2}$, so we are to
choose  appropriate $\eta_2, \ldots, \eta_r$ to minimize
$\|v_{2q+1}\|_2^2+\|v_{2q+2}\|_2^2$, the first two terms of the function defined in  \eqref{dep.}.

  Define $S_{2q,p}^{(2)}V_{2q,p}=$$\begin{bmatrix}w_1& W\end{bmatrix}$ with $w_1\in\mathbb{C}^p$,
  $Y_1=$$\begin{bmatrix}\mbox{Re}(W)& -\mbox{Im}(W)\end{bmatrix}$, $Y_2=$$\begin{bmatrix}\mbox{Im}(W)& \mbox{Re}(W)\end{bmatrix}$,
  and $\mbox{Re}(y)+i\mbox{Im}(y)=y=\begin{bmatrix}\eta_2&\cdots&\eta_r\end{bmatrix}^\top$,
  then with some simple  computations, we have
 \begin{equation}\label{opt}
  \begin{split}
  &\|v_{2q+1}\|_2^2+ \|v_{2q+2}\|_2^2 \\
 =&\begin{bmatrix}\mbox{Re}(y)^\top & \mbox{Im}(y)^\top\end{bmatrix}H
  \begin{bmatrix}\mbox{Re}(y)^\top & \mbox{Im}(y)^\top\end{bmatrix}^\top \\
  & \quad + g^\top \begin{bmatrix}\mbox{Re}(y)^\top & \mbox{Im}(y)^\top\end{bmatrix}^\top + \zeta,
  \end{split}
  \end{equation}
 where
\begin{align*}
H=& \frac{1}{\|\tilde{x}_{2q+1}\|_2^2}(cY_1-sY_2)^\top(cY_1-sY_2)\\
  &  +\frac{1}{\|\tilde{x}_{2q+2}\|_2^2}(sY_1+cY_2)^\top(sY_1+cY_2),\\
g=&\frac{2}{\sigma_1}\left(\frac{c^2}{\|\tilde{x}_{2q+1}\|_2^2}+\frac{s^2}{\|\tilde{x}_{2q+2}\|_2^2}\right)Y_1^\top\mbox{Re}(w_1) \\ &+\frac{2}{\sigma_1}\left(\frac{s^2}{\|\tilde{x}_{2q+1}\|_2^2}+\frac{c^2}{\|\tilde{x}_{2q+2}\|_2^2}\right)Y_2^\top\mbox{Im}(w_1) \\
&+\frac{2cs}{\sigma_1}\left(\frac{1}{\|\tilde{x}_{2q+2}\|_2^2}
-\frac{1}{\|\tilde{x}_{2q+1}\|_2^2}\right)(Y_2^\top\mbox{Re}(w_1)+ Y_1^\top \mbox{Im}(w_1)),\\
\zeta= &\left(\frac{c^2}{\|\tilde{x}_{2q+1}\|_2^2}+\frac{s^2}{\|\tilde{x}_{2q+2}\|_2^2}\right)\frac{\|\mbox{Re}(w_1)\|_2^2}{\sigma_1^2}\\
&+\left(\frac{s^2}{\|\tilde{x}_{2q+1}\|_2^2}+\frac{c^2}{\|\tilde{x}_{2q+2}\|_2^2}\right)\frac{\|\mbox{Im}(w_1)\|_2^2}{\sigma_1^2}\\
&  +\frac{ 2cs}{\sigma_1^2}\left(\frac{1}{\|\tilde{x}_{2q+2}\|_2^2}
  -\frac{1}{\|\tilde{x}_{2q+1}\|_2^2}\right)\mbox{Re}(w_1)^\top \mbox{Im}(w_1).
  \end{align*}
Apparently, $H$ is symmetric semipositive definite. We can further show that $H$ is nonsingular, that is, it is positive definite. Indeed, assume that $f\in\mathbb{R}^{2r-2}$ satisfies $Hf=0$, which is then equivalent to $Y_1f=Y_2f=0$ by the definition of $H$. Using the definitions of $Y_1,Y_2$ and $W$, we have
\begin{align}\label{Y1Y2}
Y_1^\top Y_1+Y_2^\top Y_2=I_{2(r-1)}.
\end{align}
So it must hold that $f=0$, which implies that $H$ is symmetric positive definite.
Consequently, the minimizer of  \eqref{opt} can be given by
\[
\begin{bmatrix}\mbox{Re}(y)^\top & \mbox{Im}(y)^\top\end{bmatrix}^\top=-\frac{1}{2}H^{-1}g.
\]

Accordingly, $v_{2q+1}$ and $v_{2q+2}$ can be computed by \eqref{v1v2}.  In all, in this case, the size of
$D_{kk}(\lambda_1)$ in $T_{2q}$ is increased by $2$, and $T_{2q+2}$ being of   the form of \eqref{updatereal_a_1}
will be obtained.

\item{\bf Case \romannumeral5: $\dim(\mathcal{R}(S_{2q,p}^{(1)}))=1$ and $\mbox{Re}(u),\mbox{Im}(u)$ are linearly dependent, or $\dim(\mathcal{R}(S_{2q,p}^{(1)}))=0$.} \quad
In this case, we cannot find $x_{2q+1}$, $x_{2q+2}$ and $v_{2q+1}, v_{2q+2}\in\mathbb{R}^p$ satisfying \eqref{eqxq},  meaning that $T_{2q+2}$ cannot be chosen in the form of \eqref{updatereal_a_1}. Instead, we set $T_{2q+2}$ in the form of \eqref{updatereal_a_2} to continue the assigning process, which leads to:
\begin{align}\label{eqxq-new}
\left\{\begin{array}{l}
Q_2^\top (A \begin{bmatrix}x_{2q+1}& x_{2q+2}\end{bmatrix}
-X_{2q} \begin{bmatrix}v_{2q+1} &v_{2q+2}\end{bmatrix}\\
 \qquad - \begin{bmatrix}x_{2q+1}& x_{2q+2}\end{bmatrix}D(\delta_{1,k+1}(\lambda_1)))=0,\\
X_{2q}^\top \begin{bmatrix}x_{2q+1}& x_{2q+2}\end{bmatrix}=0,\\
\begin{bmatrix}x_{2q+1}& x_{2q+2}\end{bmatrix}^\top \begin{bmatrix}x_{2q+1}& x_{2q+2}\end{bmatrix}=I_2,
\end{array}
\right.
\end{align}
with $v_{2q+1}, v_{2q+2}\in\mathbb{R}^{2q}$.
Denote $\delta_{1,k+1}(\lambda_1)=\frac{\delta_2}{\delta_1}$ with $0\neq\delta_1\in\mathbb{R}$ and $\delta_2\in\mathbb{R}$,
then \eqref{eqxq-new} is equivalent to some  constraints similar to those in  \eqref{eqxq2},
where the essential difference  here is that the parameter $p$ in \eqref{eqxq2} is replaced by $2q$.
More specifically, the matrix $M_{2q,p}$ in \eqref{M_2q^p} now turns to $M_{2q,2q}$,
where the $(1,2)$ block is $-Q_2^\top X_{2q}$  presently, instead of  $-Q_2^\top X_p$.
Bear in mind that now we have $v_{2q+1}\in\mathbb{R}^{2q}$ and $v_{2q+2}\in\mathbb{R}^{2q}$,
indicating that the $2\times 2$ block $T(2q+1:2q+2, 2q+1:2q+2)$ locates in the $(k+1)$-th
diagonal block $D_{k+1, k+1}(\lambda_1)$ of $T_{2a_1}$.
Now, we are to compute $x_{2q+1}$, $x_{2q+2}$, $v_{2q+1}$ and $v_{2q+2}$
satisfying some nonlinear constraints
such that the corresponding objective function specified as  \eqref{dep.} is optimized.

The forthcoming Theorem~\ref{theoremc} in Subsection \ref{subsection2.3}
demonstrates that $\dim(\mathcal{N}(M_{2q,2q}))=m$ and
there exists $\begin{bmatrix}z^\top&w^\top\end{bmatrix}^\top\in\mathcal{N}(M_{2q,2q})$ with
$z\in\mathbb{R}^n$, $w\in\mathbb{R}^{2q}$
such that $z\neq0$ and $\mbox{Re}(z)$ and $\mbox{Im}(z)$ are linearly independent, meaning  that we can always find
$x_{2q+1},x_{2q+2},v_{2q+1}$ and $v_{2q+2}$ to satisfy \eqref{eqxq-new}.

Suppose that the columns of
$S_{2q,2q}=\begin{bmatrix}S_{2q,2q}^{(1)^\top}&S_{2q,2q}^{(2)^\top}\end{bmatrix}^{\top}$
with ${S_{2q,2q}^{(1)}}\in\mathbb{C}^{n\times m}$,
$S_{2q,2q}^{(2)}\in\mathbb{C}^{2q\times m}$, form an orthonormal basis of $\mathcal{N}(M_{2q,2q})$
and let $S_{2q,2q}^{(1)}=U_{2q,2q}\Sigma_{2q,2q}V_{2q,2q}^{*}$ be the SVD of $S_{2q,2q}^{(1)}$, with the
singular values in decreasing order. Different placing strategies  based on $\rank(S_{2q,2q}^{(1)})$
 will be employed to acquire the $(2q+1)$-th and  $(2q+2)$-th  columns of $X$ and $T$. Notice that Theorem~\ref{theoremc} ensures that $\rank(S_{2q,2q}^{(1)})\ge 1$.

If $\rank(S_{2q,2q}^{(1)})=1$, then $S_{2q,2q}^{(1)}$ has only one nonzero singular value $\sigma_1$ with $u=U_{2q,2q}e_1$ being its corresponding left singular vector. Theorem~\ref{theoremc} assures that $\mbox{Re}(u)$ and $\mbox{Im}(u)$ must be linearly independent.
Then the assigning procedure is similar as that in \textbf{Case \romannumeral4}.
While $\rank(S_{2q,2q}^{(1)})>1$, the assigning procedure is similar as that  in  \textbf{Case \romannumeral3}.

Accordingly, in either situation, we can compute $x_{2q+1},x_{2q+2},v_{2q+1},v_{2q+2}$ with
$T_{2q+2}$ in the form of \eqref{updatereal_a_2}.
Moreover, in this case, $n_k$ is fixed, and $n_{k+1}$ is initially set to be 1.
\end{itemize}

The above placing process can be proceeded with until all $\{\lambda_1, \bar{\lambda}_1\}$ have been assigned.
From the assigning process, we can see that if
$T_{2q}=D_{11}(\lambda_1)$ in \eqref{T_2q},
$M_{2q, p}$ defined in \eqref{M_2q^p} would be
\[
M_{2q,0}=\begin{bmatrix}Q_2^\top(A-\lambda_1I_n)\\ X_{2q}^\top\end{bmatrix},
\]
where $\rank(M_{2q,0})\leq (n-m)+2q$. Thus provided that
$q\leq \lfloor\frac{m}{2}\rfloor-1$, we have $\dim(\mathcal{N}(M_{2q,0}))\geq 2$,
which will lead the resulted $(2q+2)\times (2q+2)$ leading principal submatrix
$T_{2q+2}$ of $T$ in the form of \eqref{updatereal_a_1}, i.e.,
$T_{2q+2}=\diag(T_{2q}, D(\delta_{q+1, 1}(\lambda_1)))$, suggesting
that the size of the first diagonal block in $T_{2a_1}$ is increased by $2$.
Consequently, in the case of $a_1\leq\lfloor\frac{m}{2}\rfloor$,
both $\lambda_1$ and $\bar{\lambda}_1$ can be placed with $g_1=a_1$,
that is, they are assigned as semi-simple eigenvalues of $A_c=A+BF$.

The procedure assigning $\{\lambda_1, \bar{\lambda}_1\}$ is summarized in the following Algorithm \ref{algorithm1c}.

\begin{algorithm}
\caption{ \ Assigning complex conjugate $\{\lambda_1, \bar{\lambda}_1\}$ }
\renewcommand{\algorithmicrequire}{\textbf{Input:}}
\renewcommand\algorithmicensure {\textbf{Output:} }
\begin{algorithmic}[1]\label{algorithm1c}
\REQUIRE ~~\\
$A, Q_2 $, $\lambda_1\in\mathbb{C}$ with $\mbox{Im}(\lambda_1)\neq0$ and $a_1$ (the multiplicity of $\lambda_1$).
\ENSURE ~~\\
Orthogonal $X_{2a_1}$ and upper quasi-triangular $T_{2a_1}$.
\STATE Find $S=S_1+iS_2$, whose columns form an orthonormal basis of $\mathcal{N}(Q_2^\top(A-\lambda_1I_n))$.
\IF {$S_1^\top S_2=0 \ \text{and}\  S_1^\top S_1=\frac{1}{2}I_m$}
    \STATE Set $y_1, y_2 \in \mathbb{R}^{m}$ be any vectors with $\|y_1\|_2=\|y_2\|_2=1$; compute $x_1, x_2$ by \eqref{x1x2} and set
    $T_2=D_0(\lambda_1)$.
\ELSE
   \STATE Compute $x_1, x_2$ by \eqref{x1x2} with  $y_1, y_2 \in \mathbb{R}^{m}$ defined as in \eqref{gammazeta}; normalize $x_1$, $x_2$
   and set $T_2=D_0(\lambda_1)$.
\ENDIF
\STATE Set $j=2, k=0$.
\WHILE {$j<2a_1$}
    \STATE  Find 
        \[
        S_{j,k}=\begin{bmatrix}S_{j,k}^{(1)}\\S_{j,k}^{(2)}\end{bmatrix}
        \begin{array}{l}n\\ k\end{array},
        \]
        whose columns form an orthonormal basis of the null space of
        $M_{j,k}=\begin{bmatrix}Q_2^\top(A-\lambda_1I_n)&-Q_2^\top X_k\\  X_{j}^\top&0\end{bmatrix}$; compute the  SVD of $S_{j,k}^{(1)}=U_{j,k}\Sigma_{j,k}V_{j,k}^*$.
    \IF { $\rank( S_{j,k}^{(1)})\geq2$ }
        \STATE Compute the $(j+1)$-th and $(j+2)$-th columns of $X_{2a_1}$ and $T_{2a_1}$
          as  in {\bf Case \romannumeral3}; set $j=j+2$.
    \ELSIF {$\rank( S_{j,k}^{(1)})=1$ and $\mbox{Re}(U_{j,k}e_1)$ and $\mbox{Im}(U_{j,k}e_1)$ are linearly independent}
         \STATE Compute the $(j+1)$-th and $(j+2)$-th columns of  $X_{2a_1}$ and $T_{2a_1}$
         as in {\bf Case \romannumeral4}; set $j=j+2$.
    \ELSE
         \STATE Find 
            \[
            S_{j,j}=\begin{bmatrix}S_{j,j}^{(1)}\\S_{j,j}^{(2)}\end{bmatrix}
            \begin{array}{l}n\\j\end{array},
            \]
            whose columns form an orthonormal basis of  the null space of              $M_{j,j}=\begin{bmatrix}Q_2^\top(A-\lambda_1I_n)&-Q_2^\top X_{j}\\  X_j^\top&0\end{bmatrix}$; compute the $(j+1)$-th and $(j+2)$-th columns of  $X_{2a_1}$ and $T_{2a_1}$ as in {\bf Case \romannumeral5}; set $k=j$ and $j=j+2$.
    \ENDIF
\ENDWHILE
\end{algorithmic}
\end{algorithm}

\subsection{Assigning repeated poles $\lambda_{j+1}\,(j\ge 1)$}\label{subsection2.2}

Suppose that the poles $\lambda_1,\dots,\lambda_j$ have been assigned.
Here the set $\{\lambda_1,\dots,\lambda_j\}$ is closed under complex conjugate.
That is, we have already obtained
$X_{r_0}=\begin{bmatrix}x_1&x_2&\cdots&x_{r_0}\end{bmatrix} \in\mathbb{R}^{n\times r_0}$
and the ${r_0\times r_0}$ leading principal submatrix $T_{r_0}$ of $T$ satisfying
\begin{align*}
Q_2^\top (AX_{r_0} - X_{r_0}T_{r_0})=0, &\qquad
X_{r_0}^\top X_{r_0}=I_{r_0},
\end{align*}
where $r_0=\sum_{k=1}^{j}a_k$ with $a_1, \ldots, a_j$ being the multiplicities of $\lambda_1, \ldots, \lambda_j$, respectively, and
\[
\lambda(T_{r_0})=\{\underbrace{\lambda_1, \ldots, \lambda_1}_{a_1}, \ldots,
\underbrace{\lambda_j, \ldots, \lambda_j}_{a_j} \}\subset\mathfrak{L}.
\]
Then we are to assign $\lambda_{j+1}$ with  multiplicity $a_{j+1}$.
Here we assume $a_{j+1}>1$. Similarly, we will again distinguish into two different cases when $\lambda_{j+1}$ is real or non-real.

\bigskip

\subsubsection{$\lambda_{j+1}$ is real}\label{subsection3.2.1}
To make the geometric multiplicity of $\lambda_{j+1}$ as large as possible, we take
$T(r_0+1:r_0+a_{j+1}, r_0+1:r_0+a_{j+1})$, the block diagonal part in $T$ corresponding to $\lambda_{j+1}$,
in the  special form of
\begin{equation}\label{T_form2}
\begin{split}
&T(r_0+1:r_0+a_{j+1}, r_0+1:r_0+a_{j+1})\\
=&\begin{array}{lll}
& \begin{array}{llll}
\qquad n_1&\qquad \quad n_2  & \qquad  \cdots&\quad  \ \  n_l
\end{array}&\\
&\begin{bmatrix}D_{11}(\lambda_{j+1})&*&\cdots&*\\ & D_{22}(\lambda_{j+1})&\cdots&*\\&&\ddots&\vdots\\&&&D_{ll}(\lambda_{j+1})\end{bmatrix}&
\begin{array}{l}
n_1\\n_2\\\vdots\\n_l
\end{array}
\end{array},
\end{split}
\end{equation}
where $D_{kk}(\lambda_{j+1})=\lambda_{j+1}I_{n_k}$, $k=1, \ldots, l$,  and $\sum_{k=1}^{l}n_k=a_{j+1}$.
With this form, the geometric multiplicity of $\lambda_{j+1}$ will be no less than
$\max\{n_k:\ k=1, \ldots, l\}$.
Theoretically, if $n_1=a_{j+1}$, $\lambda_{j+1}$ achieves its
maximum geometric multiplicity and serves as a semi-simple eigenvalue of $A_c$,
which is the  most desirable. However, $n_1$ can not be chosen to be equal to $a_{j+1}$ in some cases.

The assigning process of obtaining the columns of $X$ and $T$ corresponding to the
first diagonal block $D_{11}(\lambda_{j+1})$ in \eqref{T_form2} is as below. By noting the form of
$T(r_0+1: r_0+a_{j+1},r_0+1: r_0+a_{j+1})$ in \eqref{T_form2}, then
comparing the $(r_0+1)$-th to  the $(r_0+n_1)$-th columns
of \eqref{eqXT} shows that the corresponding columns of $X$ and $T$ must satisfy
$\begin{bmatrix}x_{r_0+k}^\top&v_{r_0+k}^\top\end{bmatrix}^\top\in\mathcal{N}(M_{r_0, r_0})$ for $1\leq k\leq n_1$, where
\begin{align}\label{M-a^{j+1}}
M_{r_0,r_0}=\begin{bmatrix}Q_2^\top(A-\lambda_{j+1}I_n)& -Q_2^\top X_{r_0}\\X_{r_0}^\top&0\end{bmatrix},
\end{align}
and $\breve{v}_{r_0+k}=\begin{bmatrix}v_{r_0+k}^\top&0\end{bmatrix}^\top$, $v_{r_0+k}\in\mathbb{R}^{r_0}$ for $1\leq k\leq n_1$.
Let the columns of
\begin{align}\label{S1q0}
S_{r_0, r_0}=\begin{bmatrix}\ S_{r_0,r_0}^{(1)}\\S_{r_0,r_0}^{(2)}\end{bmatrix}
\begin{array}{l}n\\r_0\end{array}
\end{align}
be an orthonormal basis of $\mathcal{N}(M_{r_0, r_0})$.
Write $r_{r_0}=\rank(S_{r_0,r_0}^{(1)})$, which indicates that
we can  select at most $r_{r_0}$ linearly independent
vectors from $\mathcal{R}(S_{r_0,r_0}^{(1)})$. That is, $n_1$ cannot exceed $r_{r_0}$.
Similarly as  the previous subsection --- Subsection \ref{subsection3.1.1}, $r_{r_0}$ must be nonzero to assure
that the assigning procedure would not interrupt.
The  related results are summarized in Theorem \ref{theorem} in Subsection \ref{subsection2.3}.
In the following, two different cases will be  disposed separately.

\begin{itemize}
\item {\bf Case \romannumeral1: $a_{j+1}\leq r_{r_0}$.}\quad
In this case, we set $n_1=a_{j+1}$. With this choice,
$\lambda_{j+1}$ will act as a semi-simple eigenvalue of $A_c$.
Then to get a small departure from normality of $A_c$,
it is natural to consider the following optimization problem:
\begin{subequations}\label{opt2}
\begin{align}
&\min_{}\|\begin{bmatrix}v_{r_0+1}& v_{r_0+2}&\ldots&v_{r_0+a_{j+1}}\end{bmatrix}\|_F^2 \label{eqreal-opt-j+1}\\
\mbox{s.t.}&\left\{ \begin{array}{l}
 M_{r_0,r_0}\begin{bmatrix}x_{r_0+1}& \cdots &x_{r_0+a_{j+1}}\\
 v_{r_0+1}& \cdots&v_{r_0+a_{j+1}} \end{bmatrix} =0,\\
 \\
\begin{bmatrix}x_{r_0+1}& \cdots&x_{r_0+a_{j+1}}\end{bmatrix}^{\top}\\
\qquad \begin{bmatrix}x_{r_0+1}& \cdots&x_{r_0+a_{j+1}}\end{bmatrix}=I_{a_{j+1}}.\label{eqreal-opt-constrain-j+1}
\end{array}
\right.
\end{align}
\end{subequations}
Apparently, it can be solved by the same method that solves \eqref{opt1}.
Once the solution is obtained, $X_{r_0}$ and $T_{r_0}$ will be  updated as
\begin{equation}\label{updatereal_a_j_1}
\begin{split}
&X_{r_0+a_{j+1}}\\
=&\begin{bmatrix}X_{r_0}&x_{r_0+1}&\ldots&x_{r_0+a_{j+1}}
\end{bmatrix}\in\mathbb{R}^{n \times ( r_0+a_{j+1} )},\\
\\
&T_{r_0+a_{j+1}}\\
=&\left[\begin{array}{c|c}
T_{r_0} \mathbf & \begin{array}{ccc}
v_{r_0+1}&\cdots&v_{r_0+a_{j+1}}\\
\end{array}\\
& \\[-2mm]
\hline
& \\[-2mm]
\mathbf &\lambda_{j+1}I_{a_{j+1}}\\
\end{array}\right]
\in\mathbb{R}^{{(r_0+a_{j+1})}\times {(r_0+a_{j+1})}},
\end{split}
\end{equation}
where $T_{r_0+a_{j+1}}$ is the $(r_0+a_{j+1})\times (r_0+a_{j+1})$ leading
principal submatrix of $T$.

\smallskip

\item {\bf Case \romannumeral2: $a_{j+1}> r_{r_0}$.} \quad
In this case, the maximum possible value of $n_1$ is $r_{r_0}$, and we then  set $n_1=r_{r_0}$.
Similarly to {\bf Case \romannumeral2} in Subsection \ref{subsection3.1.1},
let $S_{r_0,r_0}^{(1)}=U_{r_0, r_0}\Sigma_{r_0, r_0}V_{r_0, r_0}^{\top}$ be the SVD of
$S_{r_0,r_0}^{(1)}$ with $\sigma_{1,r_0}, \ldots, \sigma_{ r_{r_0}, r_0}$ being its singular values, then we take
\begin{align*}
\begin{bmatrix}x_{r_0+1}&\cdots&x_{r_0+r_{r_0}}\end{bmatrix}=
U_{r_0, r_0}\begin{bmatrix}e_1&\cdots&e_{r_{r_0}}\end{bmatrix},
\end{align*}
\begin{align*}
&\begin{bmatrix}v_{r_0+1}&\ldots&v_{r_0+r_{r_0}}\end{bmatrix}\\
=&S_{r_0,r_0}^{(2)}V_{r_0,r_0}\begin{bmatrix}e_1&\ldots&e_{r_{r_0}}\end{bmatrix}
\diag( \frac{1}{\sigma_{1,r_0}}, \ldots, \frac{1}{\sigma_{ r_{r_0}, r_0} }) ,
\end{align*}
and update $X_{r_0}$ and $T_{r_0}$ as
\begin{equation}\label{updatereal_a_j_1_r}
\begin{split}
&X_{r_0+n_1}=X_{r_0+r_{r_0}}\\
=&\begin{bmatrix}X_{r_0}&x_{r_0+1}&\cdots&x_{r_0+r_{r_0}}
\end{bmatrix}\in\mathbb{R}^{n \times (r_0+r_{r_0} )},\\
\\
&T_{r_0+n_1}=T_{r_0+r_{r_0}}\\
=&\left[\begin{array}{c|c}
T_{r_0} \mathbf & \begin{array}{ccc}
v_{r_0+1}&\cdots&v_{r_0+r_{r_0}}\\
\end{array}\\
& \\[-2mm]
\hline
& \\[-2mm]
\mathbf &\lambda_{j+1}I_{r_{r_0}}\\
\end{array}\right]
\in\mathbb{R}^{{(r_0+r_{r_0})}\times {(r_0+r_{r_0})}}.
\end{split}
\end{equation}
\end{itemize}

\bigskip

Hence, if $a_{j+1}\le r_{r_0}$, all $\lambda_{j+1}$ have been assigned,
and we can continue with $\lambda_{j+2}$;
while in the case of $a_{j+1}>r_{r_0}$, 
we still need to perform a similar procedure as {\bf Case \romannumeral1} and
{\bf Case \romannumeral2}   until all $\lambda_{j+1}$ are assigned.
Ultimately, we would acquire the $(r_0+a_{j+1}) \times (r_0+a_{j+1})$ leading principal submatrix of $T$ being of the form
\begin{align*}
&T_{r_0+a_{j+1}}\\
=&\begin{bmatrix}T_{r_0}&*&\cdots&*\\
&\lambda_{j+1}I_{n_1}&\cdots&*\\ & & \ddots &\vdots\\ &&&
\lambda_{j+1}I_{n_l}\end{bmatrix}\in\mathbb{R}^{(r_0+a_{j+1})\times (r_0+a_{j+1})},
\end{align*}
where $\sum_{k=1}^{l}n_k=a_{j+1}$. Furthermore, the geometric multiplicity
$g_{j+1}$ of $\lambda_{j+1}$ satisfies $\max\{n_k:\ k=1, \ldots, l\} \leq  g_{j+1}\leq m $.
We synthesize the assigning process of $\lambda_{j+1}$ in
Algorithm \ref{algorithm2}.

\begin{algorithm}
\caption{ \ Assigning real $\lambda_{j+1}$  }
\renewcommand{\algorithmicrequire}{\textbf{Input:}}
\renewcommand\algorithmicensure {\textbf{Output:} }
\begin{algorithmic}[1]\label{algorithm2}
\REQUIRE ~~\\
$A, Q_2, X_{r_0}, T_{r_0}$, $\lambda_{j+1}\in\mathbb{R}$ and $a_{j+1}$ (the multiplicity of $\lambda_{j+1}$).
\ENSURE ~~\\
Orthogonal $X_{r_0+a_{j+1}}$ and upper quasi-triangular $T_{r_0+a_{j+1}}$.
\STATE Set $q=0$.
\WHILE {$q<a_{j+1}$}
     \STATE Find          \[
                          \begin{array}{lll}
                          S=&\begin{bmatrix}S_1\\S_2\end{bmatrix} &\begin{array}{l}n \\ r_0+q \end{array}
                          \end{array},
                           \]
                  whose columns form an orthonormal basis of $\mathcal{N}(M_{r_0+q, r_0+q})$,  where
                  \begin{align*}
                   M_{r_0+q, r_0+q}=\begin{bmatrix}Q_2^\top(A-\lambda_{j+1}I_n)& -Q_2^\top X_{r_0+q}\\X_{r_0+q}^\top&0\end{bmatrix};
                  \end{align*}
            \IF{ $(a_{j+1}-q)\leq \rank(S_1)$ }
              \STATE Solve the optimization problem \eqref{opt2}
                     with $r_0$ replaced by $(r_0+q)$ and $a_{j+1}$ by $(a_{j+1}-q)$;
              \STATE Update $X_{r_0+q}$ and $T_{r_0+q}$ similarly as \eqref{updatereal_a_j_1}, set $q=a_{j+1}$.
           \ELSE
              \STATE Update $X_{r_0+q}$ and $T_{r_0+q}$ similarly as  \eqref{updatereal_a_j_1_r}, set $q=q+\rank(S_1)$.
           \ENDIF
\ENDWHILE
\end{algorithmic}
\end{algorithm}

\bigskip

\subsubsection{$\lambda_{j+1}$ is non-real}\label{subsection3.2.2}
Let $\lambda_{j+1}=\alpha_{j+1}+i\beta_{j+1}$ with $\alpha_{j+1},  \beta_{j+1}\in \mathbb{R}$ and $\beta_{j+1}\neq0$.
In this part, we shall sketch the process of assigning all  complex conjugate pairs $\{\lambda_{j+1}, \bar{\lambda}_{j+1}\}$.
Denote the algebraic multiplicity  and geometric multiplicity of $\lambda_{j+1}$ (and $\bar{\lambda}_{j+1}$)
by $a_{j+1}$ and $g_{j+1}$, respectively.
To make the geometric multiplicity $g_{j+1}$ as large as possible,  similarly as $T_{2a_1}$ in Subsection \ref{subsection3.1.2},
we take $T(r_0+1:r_0+2a_{j+1}, r_0+1:r_0+2a_{j+1})$ in the special form of
\begin{equation}\label{T_2j1}
\begin{split}
&T(r_0+1:r_0+2a_{j+1}, r_0+1:r_0+2a_{j+1})\\
=&\begin{array}{ll}
\begin{array}{llll}
\quad 2n_1& \qquad \quad \ 2n_2&\qquad  \cdots&\quad  2n_l
\end{array}&\\
\begin{bmatrix}D_{11}(\lambda_{j+1})&*&\cdots&*\\ & D_{22}(\lambda_{j+1})&\cdots&*\\&&\ddots&\vdots\\&&&D_{ll}(\lambda_{j+1})\end{bmatrix}&
\begin{array}{l}
2n_1\\2n_2\\\ \vdots\\2n_l
\end{array}
\end{array},
\end{split}
\end{equation}
where $D_{kk}(\lambda_{j+1})=\diag(D(\delta_{1,k}(\lambda_{j+1})), \ldots, D(\delta_{n_k,k}(\lambda_{j+1})))$ with
\begin{equation}\label{D_delta}
\begin{split}
&D(\delta_{p,k}(\lambda_{j+1}))\\
=&\begin{bmatrix}\mbox{Re}(\lambda_{j+1})&\delta_{p,k}(\lambda_{j+1})\mbox{Im}(\lambda_{j+1})\\
-\frac{1}{\delta_{p,k}(\lambda_{j+1})}\mbox{Im}(\lambda_{j+1})&\mbox{Re}(\lambda_{j+1})\end{bmatrix},
\end{split}
\end{equation}
$0 \neq \delta_{p,k}(\lambda_{j+1})\in\mathbb{R}$, $p=1, \ldots, n_k$, $k=1, \ldots, l$,  and $\sum_{k=1}^{l}n_k=a_{j+1}$.
Apparently, as eigenvalues of $A_c$,
the geometric multiplicity $g_{j+1}$ of $\lambda_{j+1}$ (and $\bar{\lambda}_{j+1}$)
is no less than $\max\{n_k:\ k=1, \ldots, l\}$.

Similarly as that in  Subsection \ref{subsection3.1.2}, we shall place one complex conjugate pair
$\{\lambda_{j+1}, \bar{\lambda}_{j+1}\}$ at a time,
obtaining two columns of $T$ and $X$ corresponding to the
$2\times 2$ matrix $D(\delta_{p,k}(\lambda_{j+1}))$ concurrently.

Firstly, we dispose the issue that  how to  obtain the $(r_0+1)$-th and $(r_0+2)$-th columns of $X$ and $T$.
Notice that $T(r_0+1: r_0+2, r_0+1: r_0+2)=D(\delta_{1,1}(\lambda_{j+1}))$.
Define $\delta_{1,1}(\lambda_{j+1})=\frac{\delta_2}{\delta_1}$ with $0\neq \delta_1\in\mathbb{R}$ and $\delta_2\in\mathbb{R}$,
then it follows from \cite{GCQX} that
\begin{align}\label{M_r0}
M_{r_0,r_0}
\begin{bmatrix}\frac{1}{\delta_1}x_{r_0+1}+i\frac{1}{\delta_2}x_{r_0+2}\\
\frac{1}{\delta_1}v_{r_0+1}+i\frac{1}{\delta_2}v_{r_0+2}\end{bmatrix}=0,
\end{align}
where the definition of $M_{r_0,r_0}$ is analogous to that specified in \eqref{M-a^{j+1}}
and  $\breve{v}_{r_0+k}=\begin{bmatrix}v_{r_0+k}^\top&0\end{bmatrix}^\top$,
$v_{r_0+k}\in\mathbb{R}^{r_0}$ for $k=1, \ 2$.
And  the intrinsical changing  on  $M_{r_0,r_0}$  is  that now $\lambda_{j+1}\in\mathbb{C}$ with
$\mbox{Im}(\lambda_{j+1})\neq 0$.
Accordingly, to get proper $x_{r_0+1}$, $x_{r_0+2}$, $v_{r_0+1}$, $v_{r_0+2}$,
$\delta_1$ and $\delta_2$, we need to minimize the  function defined in \eqref{dep.}
subject to the two  constraints \eqref{M_r0} and
$\begin{bmatrix}x_{r_0+1}&x_{r_0+2}\end{bmatrix}^\top \begin{bmatrix}x_{r_0+1}&x_{r_0+2}\end{bmatrix}=I_2$.

Theorem \ref{theoremc} in the forthcoming Subsection \ref{subsection2.3} shows that
$\dim(\mathcal{N}(M_{r_0,r_0}))=m$ and there exists
$\begin{bmatrix}z^\top & w^\top \end{bmatrix}^\top\in\mathcal{N}(M_{r_0,r_0})$
with $0\neq z\in\mathbb{C}^n$, $w\in\mathbb{C}^{r_0}$ and $\mbox{Re}(z)$, $\mbox{Im}(z)$ being linearly independent.
Define  $S_{r_0,r_0}=\begin{bmatrix}S_{r_0,r_0}^{(1)^\top} & S_{r_0,r_0}^{(2)^\top}\end{bmatrix}^\top$ with $S_{r_0,r_0}^{(1)}\in\mathbb{C}^{n\times m}$, $S_{r_0,r_0}^{(2)}\in\mathbb{C}^{r_0\times m}$,
whose columns form an  orthonormal basis of $\mathcal{N}(M_{r_0,r_0})$, the
placing process will be realized through addressing
two distinct cases upon $\rank(S_{r_0,r_0}^{(1)})$.
For convenience, we denote the left  and right singular vectors of
$S_{r_0,r_0}^{(1)}$, corresponding to its largest singular value $\sigma_1$, by $u$ and $v$, respectively.

If $\rank(S_{r_0,r_0}^{(1)})\geq 2$, a similar placing  process as that in {\bf Case \romannumeral3 }
in  Subsection  \ref{subsection3.1.2} will be implemented. That is, if
$\mbox{Re}(u)^\top\mbox{Im}(u)=0$ and $\|\mbox{Re}(u)\|_2=\|\mbox{Im}(u)\|_2=\frac{\sqrt{2}}{2}$,
we set $x_{r_0+1}=\sqrt{2}\mbox{Re}(u)$, $x_{r_0+2}=\sqrt{2}\mbox{Im}(u)$, and
$v_{r_0+1}=\sqrt{2}\mbox{Re}(S_{r_0,r_0}^{(2)}v/\sigma_1)$, $v_{r_0+2}=\sqrt{2}\mbox{Im}(S_{r_0,r_0}^{(2)}v/\sigma_1)$.
Otherwise,  the complex conjugate pair placing strategy in \cite{GCQX} would be applied.
When  $\rank(S_{r_0,r_0}^{(1)})=1$, Theorem \ref{theoremc} in the following subsection
would guarantee that $\mbox{Re}(u)$ and $\mbox{Im}(u)$
are linearly independent. We then apply the Jacobi  orthogonal  transformation in Lemma \ref{Lemma2.2} to orthogonalize
$\mbox{Re}(u)$ and $\mbox{Im}(u)$, and then normalize the resulted vectors as $x_{r_0+1}$ and  $x_{r_0+2}$.
Furthermore, $v_{r_0+1}$ and $v_{r_0+2}$ will be obtained  by minimizing
some function defined similarly as that  in \eqref{opt}.
The process resembles that in {\bf Case \romannumeral4 } in Subsection \ref{subsection3.1.2}, and we omit details here.
%

Now assume that we have obtained $2q$ ( $1\leq q<a_{j+1}$) columns of $X$ and $T$
corresponding to $\{\lambda_{j+1}, \bar{\lambda}_{j+1}\}$,
we then proceed to compute the $(r_0+2q+1)$-th and   $(r_0+2q+2)$-th columns of $X$ and $T$,  which virtually are
associated with the diagonal block $T(r_0+2q+1:r_0+2q+2, r_0+2q+1:r_0+2q+2)$ in $T$.
The whole procedure is similar to what we do
to get the  $(2q+1)$-th and   $(2q+2)$-th columns of $X$ and $T$ in
Subsection \ref{subsection3.1.2}, and we just give a concise presentation.

Assume that
\begin{align*} 
T_{r_0+2q}=\left[\begin{array}{c|c}
T_{r_0} \mathbf &
\begin{array}{ccc}
*& \qquad \cdots& \qquad *\\
\end{array}  \\
& \\[-2mm]
\hline
& \\[-2mm]
\mathbf &
\begin{array}{ccc}
D_{11}(\lambda_{j+1})&\cdots&*\\
& \ddots&\vdots\\
&&D_{tt}(\lambda_{j+1})
\end{array} \\
\end{array}\right],
\end{align*}
where $D_{kk}(\lambda_{j+1})\in\mathbb{R}^{2n_k\times 2n_k}$, $k=1, \ldots, t$, and
$T(r_0+2q-1:r_0+2q,  r_0+2q-1:r_0+2q)=D(\delta_{s,t}(\lambda_{j+1}))$,
indicating that  $T(r_0+2q-1:r_0+2q, r_0+2q-1:r_0+2q)$ is the
$s$-th $2\times 2$ diagonal block in $D_{tt}(\lambda_{j+1})$.
Denote $p=r_0+2n_1+\cdots+2n_{t-1}$.
Then like $T_{2q+2}$ in Subsection \ref{subsection3.1.2}, $T_{r_0+2q+2}$ could be in the form of
\begin{equation}\label{updatereal_j+1}
\begin{split}
&T_{r_0+2q+2}=\left[\begin{array}{c|c|c}
T_p \mathbf & *   \mathbf &
\begin{array}{cc}
v_{r_0+2q+1}&v_{r_0+2q+2}\\
\end{array}  \\
& \\[-2mm]
\hline
& \\[-2mm]
\mathbf &  D_{tt}(\lambda_{j+1}) \mathbf & 0\\
& \\[-2mm]
\hline
& \\[-2mm]
\mathbf &  \mathbf & D(\delta_{s+1, t}(\lambda_{j+1}))\\
\end{array}\right],\\
& v_{r_0+2q+1},\ v_{r_0+2q+2}\in\mathbb{R}^{p},
\end{split}
\end{equation}
or
\begin{align*}
&T_{r_0+2q+2}=\left[\begin{array}{c|c}
T_{r_0+2q} \mathbf & \begin{array}{cc}
v_{r_0+2q+1}&v_{r_0+2q+2}\\
\end{array}\\
& \\[-2mm]
\hline
& \\[-2mm]
\mathbf & D(\delta_{1,t+1}(\lambda_{j+1}))\\
\end{array}\right],\\
&v_{r_0+2q+1},\ v_{r_0+2q+2}\in\mathbb{R}^{r_0+2q}.
\end{align*}
And to get a large $g_{j+1}$, we incline to $T_{r_0+2q+2}$ being  of the form in \eqref{updatereal_j+1},
which suggests that we need to regard the  null space of $M_{r_0+2q,p}$, where
\begin{align}\label{M_a_j+1_1}
M_{r_0+2q,p}=
\begin{bmatrix}Q_2^{\top}(A-\lambda_{j+1}I_n)&-Q_2^\top X_{p}\\
X_{r_0+2q}^\top&0\end{bmatrix}.
\end{align}
Suppose  that the columns of
\[
S_{r_0+2q,p}=\begin{bmatrix}S_{r_0+2q,p}^{(1)}\\S_{r_0+2q,p}^{(2)}\end{bmatrix}
\begin{array}{l}n\\p\end{array}
\]
form an orthonormal basis of $\mathcal{N}(M_{r_0+2q,p})$.
Then the assigning procedure is similar as that in Subsection \ref{subsection3.1.2}, which is accomplished
by distinguishing  three different cases:
$\rank(S_{r_0+2q,p}^{(1)})\geq 2$, $\rank(S_{r_0+2q,p}^{(1)})=1$ and
$\mbox{Re}(u)$ and $\mbox{Im}(u)$ are linearly independent with $u$ being the left singular vector
of $S_{r_0+2q,p}^{(1)}$ corresponding to its only nonzero singular value, and otherwise. 

%

Guaranteed by Theorem \ref{theoremc} below,
we can proceed with the above assigning procedure till all columns of $X$ and $T$ corresponding to
$\{\lambda_{j+1}, \bar{\lambda}_{j+1}\}$ are acquired,
which eventually yields  $T(r_0+1:r_0+2a_{j+1}, r_0+1:r_0+2a_{j+1})$
being of  the special form specified in \eqref{T_2j1}.
And we recapitulate the assigning process
of the repeated complex poles $\{\lambda_{j+1}, \bar{\lambda}_{j+1}\}$ in Algorithm \ref{algorithm2c}.

\begin{algorithm}
\caption{\ Assigning complex conjugate $\{\lambda_{j+1}, \bar{\lambda}_{j+1}\}$}
\renewcommand{\algorithmicrequire}{\textbf{Input:}}
\renewcommand\algorithmicensure {\textbf{Output:} }
\begin{algorithmic}[1]\label{algorithm2c}
\REQUIRE ~~\\
$A$, $Q_2$, $X_{r_0}$, $T_{r_0}$, $\lambda_{j+1}\in\mathbb{C}$ with $\mbox{Im}(\lambda_{j+1})\neq0$ and $a_{j+1}$ (the  multiplicity of $\lambda_{j+1}$).
\ENSURE ~~\\
Orthogonal $X_{r_0+2a_{j+1}}$ and upper quasi-triangular $T_{r_0+2a_{j+1}}$.
\STATE Set $l=k=r_0$.
\WHILE {$l<r_0+2a_{j+1}$}
    \STATE  Find 
                  \[
                  S_{l,k}=\begin{bmatrix}S_{l,k}^{(1)}\\S_{l,k}^{(2)}\end{bmatrix}
                  \begin{array}{l}n\\k\end{array},
                  \]
                  whose columns form an orthonormal basis of the null space of
                  $M_{l, k}=\begin{bmatrix}Q^\top_2(A-\lambda_{j+1} I_n)&-Q_2^\top X_{k}\\
                     X_{l}^\top&0\end{bmatrix}$;
                     compute the SVD of $S_{l,k}^{(1)}=U_{l, k}\Sigma_{l, k}V_{l,k}^*$.
    \IF { $\rank(S_{l, k}^{(1)})\geq2$}
        \STATE Compute the $(l+1)$-th and $(l+2)$-th columns of $X_{r_0+2a_{j+1}}$ and $T_{r_0+2a_{j+1}}$
            as in {\bf Case \romannumeral3} in  Subsection \ref{subsection3.1.2}; set $l=l+2$;
    \ELSIF {$\rank(S_{l, k}^{(1)})=1$ and $\mbox{Re}(U_{l, k}e_1)$, $\mbox{Im}(U_{l,k}e_1)$ are linearly independent}
         \STATE Compute the $(l+1)$-th and $(l+2)$-th columns of $X_{r_0+2a_{j+1}}$ and $T_{r_0+2a_{j+1}}$
                as in {\bf Case \romannumeral4} in Subsection \ref{subsection3.1.2}; set $l=l+2$;
    \ELSE
         \STATE Find 
               \[
               S_{l,l}=\begin{bmatrix}S_{l,l}^{(1)}\\S_{l,l}^{(2)}\end{bmatrix}
               \begin{array}{l}n\\l\end{array},
               \]
               whose columns form an orthonormal basis of the null space of
               $M_{l,l}=\begin{bmatrix}Q^\top_2(A-\lambda_{j+1} I_n)&-Q_2^\top X_{l}\\
                X_{l}^\top&0\end{bmatrix}$;
                compute the $(l+1)$-th and $(l+2)$-th columns of $X_{r_0+2a_{j+1}}$ and $T_{r_0+2a_{j+1}}$
                as in {\bf Case \romannumeral5} in Subsection \ref{subsection3.1.2}; set $k=l$ and $l=l+2$.
     \ENDIF
\ENDWHILE
\end{algorithmic}
\end{algorithm}

\subsection{Theoretical support}\label{subsection2.3}
While assigning repeated real poles, the assigning procedure described in Subsections~\ref{subsection3.1.1} and \ref{subsection3.2.1} can be carried on only if the ranks of $S_{q,q}^{(1)}$ in \eqref{S1q} and $S_{r_0,r_0}^{(1)}$ in \eqref{S1q0} are nonzero, which is guaranteed by the following theorem.

\begin{theorem}\label{theorem}
Assume that $(A,B)$ is controllable.
Suppose that the poles $\lambda_1, \ldots, \lambda_j\in\mathfrak{L}$,
with  multiplicities $a_1, \ldots, a_j$, respectively, have been assigned.
Let $x_1, \ldots, x_{r}$ be the corresponding columns of $X$
obtained from the assigning process in former subsections,
where $r=\sum_{k=1}^{j}a_k$.
Assume that $\lambda\in\mathbb{R}$ is distinct from $\lambda_1, \ldots, \lambda_j$,
and has been assigned $q$ times with the corresponding columns $x_{r+1},\dots,x_{r+q}\,(r+q<n)$ in $X$
being obtained. Denote $X_{r+q}=\begin{bmatrix}x_1&\cdots&x_{r+q}\end{bmatrix}$ and
\begin{align*}
M_{r+q, r+q}=\begin{bmatrix}Q_2^\top(A-\lambda I_n)& -Q_2^\top X_{r+q}\\X_{r+q}^\top&0\end{bmatrix}.
\end{align*}
Let the columns of
\begin{align*}
S=\begin{bmatrix}\ S_{1}\\S_{2}\end{bmatrix}
\begin{array}{l}n\\r+q\end{array}
\end{align*}
be an orthonormal basis of $\mathcal{N}(M_{r+q,r+q})$.
Then $\dim(\mathcal{R}(S))=m$ and $S_1\neq 0$.
\end{theorem}

\begin{IEEEproof}
The conclusion $\dim(\mathcal{N}(M_{r+q,r+q}))=m$ is just that $M_{r+q, r+q}$ is of full row rank.
Assume that $u\in\mathbb{R}^{n-m}$ and $v\in\mathbb{R}^{r+q}$
satisfy  $\begin{bmatrix}u^\top&v^\top\end{bmatrix}M_{r+q,r+q}=0$, that is,
\begin{subequations}
\begin{align}
u^\top Q_2^\top(A-\lambda I_n)+v^\top X_{r+q}^\top=0,\label{eqx_k1}\\
u^\top Q_2^\top X_{r+q}=0.\label{eqx_k2}
\end{align}
\end{subequations}
Post-multiplying $X_{r+q}$ on both sides of \eqref{eqx_k1} gives
\begin{align}\label{equv}
u^\top Q_2^\top(A-\lambda I_n)X_{r+q}+v^\top =0.
\end{align}
Substituting $Q_2^\top AX_{r+q}=Q_2^\top X_{r+q}T_{r+q}$ into \eqref{equv} leads to
$v=0$ and $u^\top Q_2^\top(A-\lambda I_n)=0$ by \eqref{eqx_k2}.
Thus  $u=0$ since $(A, B)$ is controllable.
So $M_{r+q,r+q}$ is of full row rank, and hence $\dim(\mathcal{N}(M_{r+q,r+q}))=m$.

Now we are to prove $S_1\neq 0$.
It holds obviously if $(r+q)<m$. We now consider the case when $(r+q)\geq m$.
Assume that $S_1=0$,
then $\rank(S_2)=m$ and $Q_2^\top X_{r+q}S_2=0$. Hence
there must exist a nonsingular matrix $W\in\mathbb{R}^{m\times m}$ such that
\begin{align}\label{X-m-rB}
X_{r+q}S_2=BW.
\end{align}
Since $Q_2^{\top}AX_{r+q}=Q_2^{\top}X_{r+q}T_{r+q}$
with $T_{r+q}$ being the $(r+q) \times (r+q)$ leading principal submatrix of $T$,
so there must exist a matrix $K\in\mathbb{R}^{m\times (r+q)}$ such that
\begin{align}\label{XT-m-r}
AX_{r+q}=X_{r+q}T_{r+q}+BK.
\end{align}
Post-multiplying $S_2$ on both sides of \eqref{XT-m-r} and
substituting \eqref{X-m-rB} into it give
$ABW=X_{r+q}T_{r+q}S_2+BKS_2$. Noticing that $W$ is nonsingular, so
\[
AB=X_{r+q}T_{r+q}S_2W^{-1}+X_{r+q}S_2W^{-1}KS_2W^{-1}.
\]
Denote $G_1=T_{r+q}S_2W^{-1}+S_2W^{-1}KS_2W^{-1}$,
then it can be simply verified by induction that
$A^kB=X_{r+q}G_k$ with $G_k=T_{r+q}G_{k-1}+S_2W^{-1}KG_{k-1}$.
And this eventually leads to
\[
\begin{bmatrix}B&AB&\cdots&A^{n-1}B\end{bmatrix}=X_{r+q}L
\]
for some $L\in\mathbb{R}^{(r+q)\times mn}$, which implies that
$\rank(\begin{bmatrix}B&AB&\cdots&A^{n-1}B\end{bmatrix})\leq (r+q)<n$,
contradicting with the controllability of $(A, B)$.
Hence $S_1\neq 0$.
\end{IEEEproof}

\bigskip

While assigning non-real repeated poles, 
continuing the assigning process is based on the facts that
the matrix $M_{j,j}$, appearing in Step 15 in Algorithm~\ref{algorithm1c},
satisfies that
$\dim(\mathcal{N}(M_{j,j}))=m$  
and there exists
$\begin{bmatrix}z^\top&w^\top\end{bmatrix}^\top\in\mathcal{N}(M_{j,j})$ with
$z\in\mathbb{R}^n$, $w\in\mathbb{R}^{j}$,
such that $z\neq0$ and $\mbox{Re}(z)$ and $\mbox{Im}(z)$ are linearly independent.
This also applies to Step 9 in Algorithm~\ref{algorithm2c}. The following Theorem then ensures that these processes
can be continued.

\begin{theorem}\label{theoremc}
Assume that $(A,B)$ is controllable. Let $\{\lambda_1, \ldots, \lambda_j\}\subset\mathfrak{L}$ be a self-conjugate subset with
$a_1, \ldots, a_j$ being the  multiplicities of $\lambda_1, \ldots, \lambda_j$, respectively,
and let $x_1, \ldots, x_{r}$ be the associate columns of $X$
obtained from the assigning process in previous subsections, where $r=\sum_{k=1}^{j}a_k$.
Assume that $\lambda=\alpha+i\beta \in \mathbb{C}$ $(\beta\neq 0)$
is some pole distinct from $\lambda_1, \ldots, \lambda_j$,
and $x_{r+1}$, $x_{r+2}$, $\ldots$, $x_{r+2q-1}, x_{r+2q}$ ($r+2q<n$) are the columns of $X$
corresponding to complex conjugate paris $\{\lambda, \bar{\lambda}\}$.
Define
\begin{align*}
M_{r+2q, r+2q}=\begin{bmatrix}Q_2^\top(A-\lambda I_n)& -Q_2^\top X_{r+2q}\\X_{r+2q}^\top&0\end{bmatrix},
\end{align*}
and let the columns of
\[
S=\begin{bmatrix}S_1\\S_2\end{bmatrix}
\begin{array}{l}n\\r+2q\end{array}
\]
be  an orthonormal basis of $\mathcal{N}(M_{r+2q,r+2q})$,
then we have
\begin{enumerate}
 \item[{\em(1)}] $\dim(\mathcal{R}(S))= m$;
 \item[{\em (2)}] $S_1\neq 0$;
 \item[{\em (3)}] there exist
$0\neq z=\mbox{Re}(z)+i\mbox{Im}(z)\in\mathbb{C}^n$ and
$w\in\mathbb{C}^{r+2q}$
with $\mbox{Re}(z)$ and $\mbox{Im}(z)$ being linearly independent,
such that $\begin{bmatrix}z^\top &w^\top\end{bmatrix}^\top\in\mathcal{R}(S)$.
\end{enumerate}
\end{theorem}

\begin{IEEEproof}
We can prove the {\em (1)}, {\em (2)} results by the  method proving  Theorem
\ref{theorem}, and we skip the  proof process here.

Regarding {\em (3)}, if $\dim(\mathcal{N}(Q_2^\top X_{r+2q}))<(m-1)$,
then there exist two vectors
$\begin{bmatrix}z_1^\top & w_1^\top\end{bmatrix}^\top, \begin{bmatrix}z_2^\top & w_2^\top\end{bmatrix}^\top\in\mathcal{R}(S)$
with $0 \neq z_1\in\mathbb{C}^n, \  0 \neq z_2\in\mathbb{C}^n$,  and $z_1$, $z_2$ being linearly independent.
Let $\begin{bmatrix}z^\top & w^\top\end{bmatrix}^\top=
(\xi_1+i\eta_1)\begin{bmatrix}z_1^\top & w_1^\top\end{bmatrix}^\top+(\xi_2+i\eta_2)\begin{bmatrix}z_2^\top & w_2^\top \end{bmatrix}^\top$, $\xi_1, \xi_2, \eta_1, \eta_2\in\mathbb{R}$,
then we can always find suitable $\xi_1, \xi_2, \eta_1, \eta_2$ such that the real part and the imaginary part of
the resulted $z$ are linearly independent. If $\dim(\mathcal{N}(Q_2^\top X_{r+2q}))=(m-1)$,
assume that $w_1, \ldots, w_{m-1}\in\mathbb{C}^{r+2q}$
form an orthonormal basis of $\mathcal{N}(Q_2^\top X_{r+2q})$ and
$0\neq z=(1+i\zeta)y$, $y\in\mathbb{R}^n$, $w\in\mathbb{C}^{r+2q}$ satisfy
$\begin{bmatrix}z^\top & w^\top\end{bmatrix}^\top\in\mathcal{R}(S)$ with $\|z\|_2^2+\|w\|_2^2=1$.
Obviously, it holds that  $Q_2^\top(A-\alpha I_n)y+\beta \zeta Q_2^\top y=Q_2^\top X_{r+2q}\mbox{Re}(w)$ and
$\zeta Q_2^\top (A-\alpha I_n)y-\beta Q_2^\top y=Q_2^\top X_{r+2q}\mbox{Im}(w)$. Thus there exist
$u, v\in\mathbb{R}^m$ such that
\begin{align}\label{y}
\left\{\begin{array}{ll}
(A-\alpha I_n)y+\beta \zeta y- X_{r+2q}\mbox{Re}(w)&=Bu, \\
\zeta(A-\alpha I_n)y-\beta y- X_{r+2q}\mbox{Im}(w)&=Bv.\end{array}
\right.
\end{align}
It follows from \eqref{y} that
\begin{subequations}
\begin{align}
&\beta(1+\zeta^2)y+X_{r+2q}(\mbox{Im}(w)-\zeta\mbox{Re}(w))\notag\\
=&\zeta Bu-Bv, \label{y1}\\
&(1+\zeta^2)(A-\alpha I_n)y-X_{r+2q}(\zeta\mbox{Im}(w)+\mbox{Re}(w))\notag\\
=&Bu+\zeta Bv. \label{y2}
\end{align}
\end{subequations}
Since $Q_2^\top X_{r+2q}\begin{bmatrix}w_1&\cdots&w_{m-1}\end{bmatrix}=0$, hence
\begin{align}\label{X-2k-BZ}
X_{r+2q}\begin{bmatrix}w_1&\cdots&w_{m-1}\end{bmatrix}=BG
\end{align}
for some $G\in\mathbb{R}^{m\times (m-1)}$ with $\rank(G)=m-1$. And it
follows from $Q_2^\top AX_{r+2q}=Q_2^\top X_{r+2q}T_{r+2q}$ that
\begin{equation}\label{AX-2k-B}
AX_{r+2q}= X_{r+2q}T_{r+2q}+BZ
\end{equation}
for some $Z\in\mathbb{R}^{m\times (r+2q)}$. Now define
\[
\begin{array}{ll}
Y&=\left[\begin{array}{c|c}
\begin{array}{ccc}
w_{1}&\cdots&w_{m-1}
\end{array}\mathbf & \mbox{Im}(w)-\zeta\mbox{Re}(w)\\
& \\[-2mm]
\hline
& \\[-2mm]
0 \mathbf &\beta(1+\zeta^2)\\
\end{array}\right], \\
\\
L&=\begin{bmatrix}G&\zeta u-v\end{bmatrix},\\
M&=\begin{bmatrix}T_{r+2q}&\frac{1}{1+\zeta^2}(\zeta \mbox{Im}(w)+\mbox{Re}(w))\\0&\alpha\end{bmatrix}, \\
E&=\begin{bmatrix}Z&\frac{1}{1+\zeta^2}(u+\zeta v)\end{bmatrix}.
\end{array}
\]
Noting  \eqref{y1}, \eqref{y2}, \eqref{X-2k-BZ} and \eqref{AX-2k-B},
then the following equations
\begin{equation}\label{A-KLMN}
\begin{split}
&\begin{bmatrix}X_{r+2q}&y\end{bmatrix}Y=BL, \\
&A\begin{bmatrix}X_{r+2q}&y\end{bmatrix}=\begin{bmatrix}X_{r+2q}&y\end{bmatrix}M+BE
\end{split}
\end{equation}
hold, where  $L$ is nonsingular since $\begin{bmatrix}X_{r+2q}&y\end{bmatrix}$ is of full column rank.
Then \eqref{A-KLMN} shows that
$AB=\begin{bmatrix}X_{r+2q}&y\end{bmatrix}H_1$ with $H_1=MYL^{-1}+YL^{-1}EYL^{-1}$.
Hence by induction, we will get that $A^{l+1}B=\begin{bmatrix}X_{r+2q}&y\end{bmatrix}H_{l+1}$,
where $H_{l+1}=MH_l+YL^{-1}EH_l$  with $l\geq 1$.
Eventually, $\begin{bmatrix}B&AB&\cdots&A^{n-1}B\end{bmatrix}=
\begin{bmatrix}X_{r+2q}&y\end{bmatrix}\begin{bmatrix}YL^{-1}&H_1&\cdots&H_{n-1}\end{bmatrix}$,
suggesting  that
\[
\text{rank}(\begin{bmatrix}B&AB&\cdots&A^{n-1}B\end{bmatrix})<n.
\]
This contradicts with the
assumption that $(A,B)$ is controllable. Thus we have proved {\em (3)}.
\end{IEEEproof}

\subsection{Algorithm}\label{alg:Framwork_new}\label{subsection2.4}

The framework of our algorithm referred to as ``\verb|Schur-multi|" is given in this subsection.
We assume that repeated real  poles appear together in $\mathfrak{L}$, while repeated complex conjugate poles appear in pairs, that is, they appear as
$\underbrace{\{\lambda, \bar{\lambda}\}, \ldots, \{\lambda, \bar{\lambda}\}}_a$ in  $\mathfrak{L}$ adjacently, where $a$ is the counting time
(the algebraic multiplicity) of $\lambda$ (and $\bar{\lambda}$) in  $\mathfrak{L}$. The  \verb|Schur-multi| algorithm below  combines techniques designed for simple poles in \cite{GCQX} and techniques for repeated poles in this paper. Again, we denote the multiplicity of $\lambda_j\in\mathfrak{L}$ by $a_j$.

\begin{algorithm}
\caption{Framework of our \texttt{ Schur-multi} algorithm.}
\renewcommand{\algorithmicrequire}{\textbf{Input:}}
\renewcommand\algorithmicensure {\textbf{Output:} }
\begin{algorithmic}[1]
\REQUIRE ~~\\
$A, B$ and $\mathfrak{L}=\{\lambda_1,\dots,\lambda_n\}$.
\ENSURE ~~\\
The feedback matrix $F$.
\STATE Compute the QR decomposition of
$B=Q\begin{bmatrix}R^{\top}& 0\end{bmatrix}^{\top}=\begin{bmatrix}Q_1&Q_2\end{bmatrix}\begin{bmatrix}R^{\top}& 0\end{bmatrix}^{\top}=Q_1R$.

\IF {$a_1=1$}
\STATE Compute the initial columns of $X$ and $T$ by \verb|Schur-rob| \cite{GCQX};
       set $j=1$ for $\lambda_1\in\mathbb{R}$ and $j=2$ for $\lambda_1\in\mathbb{C}$.
\ELSIF {$\lambda_1\in\mathbb{R}$}
\STATE Compute $X_{a_1}$ and $T_{a_1}$ by  Algorithm \ref{algorithm1}; set $j=a_1$.
\ELSE
\STATE Compute $X_{2a_1}$ and $T_{2a_1}$ by  Algorithm \ref{algorithm1c}; set $j=2a_1$.
\ENDIF
\WHILE{$j<n$}
\IF {$a_{j+1}=1$ }
    \STATE  Compute the corresponding columns of $X$ and $T$ by \verb|Schur-rob| \cite{GCQX};
            set $j=j+1$ for $\lambda_{j+1}\in\mathbb{R}$ and $j=j+2$ for $\lambda_{j+1}\in\mathbb{C}$.
\ELSIF {$\lambda_{j+1}\in\mathbb{R}$}
\STATE Compute $X_{j+a_{j+1}}$ and $T_{j+a_{j+1}}$ by  Algorithm \ref{algorithm2}; set $j=j+a_{j+1}$.
\ELSE
    \STATE Compute $X_{j+2a_{j+1}}$ and $T_{j+2a_{j+1}}$ by  Algorithm \ref{algorithm2c};
              set $j=j+2a_{j+1}$.

\ENDIF
\ENDWHILE
\STATE Compute $F$ by $F=R^{-1}Q_1^{\top}(X_nT_nX_n^{\top}-A)$.
\end{algorithmic}
\end{algorithm}

\section{Numerical examples}\label{section4}
In this section, we illustrate the performance of our \verb|Schur-multi| method by comparing
with the MATLAB functions \verb|place| \cite{KNV},
\verb|robpole| \cite{Tits} and the \verb|Schur-rob| method \cite{GCQX} on some examples.

Similarly to \cite{GCQX}, we define
\[
precs=\left\lceil\max_{1\le j\le n}(\log(|\frac{\lambda_j-\hat{\lambda}_j}{\lambda_j}|))
\right\rceil
\]
to characterize the precision of the assigned poles, where $\hat{\lambda}_j$, $j=1, \ldots, n$, are  the
computed eigenvalues of the obtained closed-loop system matrix $A_c=A+BF$.
Actually, $precs$ is the ceiling value of the exponent of the maximum relative error of
$\hat{\lambda}_j$ $(j=1, \ldots, n)$, relative to the entries in $\mathfrak{L}$.
Obviously, smaller $precs$ would imply more accurately computed poles.
Regarding the robustness of the closed-loop system, different measures are used in these methods for solving the {\bf SFRPA}. We will compare three measures  for all methods. Specifically, assume that the spectral decomposition and the real Schur decomposition of $A_c=A+BF$ respectively are
\[
A+BF=X\Lambda X^{-1},\qquad A+BF=UTU^\top,
\]
where $\Lambda$ is diagonal, $T$ is upper quasi-triangular and $U$ is orthogonal.
Then the measures adopted in \verb|place| and \verb|robpole| are closely related to the condition number
of the eigenvectors matrix $X$, i.e. $\kappa_F(X)=\|X\|_F\|X^{-1}\|_F$,
while \verb|Schur-rob| and our \verb|Schur-multi|  aim to
minimize the departure from normality of $A_c$ (denoted by $``dep."$).
We also display the Frobenius norm of the feedback matrix $F$ (denoted by $`` \|F\|_F"$),
which is also regarded as a measure of robustness in some literature.
In addition, the CPU time for all methods is also presented.
When \verb|robpole| is applied, the maximum number of sweep is set to be the default value $5$ for all examples.
All calculations are carried out by running MATLAB R2012a,
with machine epsilon $\epsilon\approx 2.2\times 10^{-16}$,
on an Intel\textregistered    Core\texttrademark i3, dual core, 2.27 GHz machine, with $2.00$ GB RAM.

The first illustrative set includes CARE examples $1.6$, $2.9$ \#$1$\cite{AB}  and DARE example $1.12$ \cite{AB2},
in which some poles are repeated and real. Additional, in the following  TABLE \ref{table_1} and  TABLE \ref{table_2},
we will use $\alpha(k)$ to represent $\alpha \times 10^{k}$ for space saving.

\begin{Example}\label{example3}
The three examples in this test set come from the SLICOT CARE/DARE benchmark collections \cite{AB}, \cite{AB2}.
The numerical results on precision and robustness for these four algorithms are
exhibited in TABLE \ref{table_1}. Concerning the CARE example $2.9$ \#$1$,
compared with \verb|Schur-rob|, our \verb|Schur-multi| does not make improvement on $``precs"$.
The reason might be that some poles are rather close to the imaginary axis.
This is a weakness of the Schur-type methods.
Note that we do not list the  $``precs"$ values for the DARE example $1.12$
since some algorithms could not achieve any relative accuracy for certain assigned poles.
And in TABLE \ref{table_2}, we display the differences between the placed poles and the eigenvalues of the
computed $A_c$ obtained from distinct methods.
The ``\verb|exact| \verb|poles|" column gives the exact values of the poles to be assigned.
TABLE \ref{table_2} shows that our \verb|Schur-multi| produces the best result on this example.

\tabcolsep 0.08in
\begin{table*}
\normalsize
\centering
\caption{ Numerical results for four algorithms on CARE/DARE examples}
\begin{tabular}{c|cccc|cccc|ccc}
\Xhline{1.0pt}
\multirow{2}*{} &  \multicolumn{4}{c|}{CARE example 1.6}&
\multicolumn{4}{c|}{CARE example 2.9 \#1}&\multicolumn{3}{c}{DARE example 1.12} \\
\Xcline{2-12}{0.5pt}
& $precs$&$dep.$ &$\kappa_F(X)$ & $\|F\|_F$&$precs$&$dep.$ &$\kappa_F(X)$ & $\|F\|_F$& $dep.$ &$\kappa_F(X)$ & $\|F\|_F$\\
\Xhline{0.7pt}
\verb|place|&-11&1.5(6)& 1.7(15)&2.2(3)&-11&2.9(6)& 8.5(4)&2.8(1)&4.3(7)& 9.2(292)&4.3(7)\\
\verb|robpole|&-13&7.5(5)& 2.2(7)&2.2(2)&-12&2.9(6)& 8.9(4)&2.8(1)&3.9(12)&1.3(308)&3.9(12)\\
\verb|Schur-rob|&-8&1.1(5)&9.0(7)& 1.2(2)&-9&7.3(6)&2.0(6)& 2.9(1)&9.8(0)&5.6(292)&6.5(0)\\
\verb|Schur-multi|&-11& 2.6(5)&1.3(7)&4.5(2)&-9&2.6(6)&1.2(6)&2.8(1)&9.1(0)&3.2(295)&5.5(0)\\
\Xhline{1.0pt}
\end{tabular}
\label{table_1}
\end{table*}

\tabcolsep 0.1in
\begin{table*}
\normalsize
\centering
\caption{ Accuracy of the assigned poles for DARE example $1.12$}
\begin{tabular}{c|c|c|c|c|c}
\Xhline{1.0pt}
\multirow{2}*{} &  \multicolumn{1}{c|}{ }&
\multicolumn{4}{c}{ \textbf{ $\lambda_j-\hat{\lambda}_j$}} \\
\Xcline{3-6}{0.5pt}
$num.$& \verb|exact poles|& \verb|place|& \verb|robpole|& \verb|Schur-rob|&\verb|Schur-multi|\\
\Xhline{0.7pt}
1&8.1(-1)&	-3.3(-16)& -3.3(-16) &-3.3(-16)&-3.3(-16)\\
2&5.8(-1)& -2.5(-7) &3.6(-5)&-1.4(-12)&2.3(-13)\\
3&1.1(-3)&8.4(-4)&2.9(-4)&-1.5(-4)&-6.4(-5)\\
4&0&-3.4(-17)&-3.4(-17)&-3.4(-17)&-3.4(-17)\\
5&0&-5.2(-17)&-5.2(-17)&-5.2(-17)&-5.2(-17)\\
6&7.6(-1)+i$\times$1.4(-1)&1.9(-7)-i$\times$1.2(-7)&-4.6(-5)-i$\times$3.7(-6)&-7.1(-13)-i$\times$1.3(-13)&6.2(-13)+i$\times$4.8(-13)\\
7&7.6(-1)-i$\times$1.4(-1)&1.9(-7)+i$\times$1.2(-7)&-4.6(-5)+i$\times$3.7(-6)&-7.1(-13)+i$\times$1.3(-13)&6.2(-13)-i$\times$4.8(-13)\\
8&6.4(-1)+i$\times$2.3(-1)&-2.5(-8)-i$\times$3.1(-8)&-4.0(-5)-i$\times$1.6(-5)&9.3(-13)-i$\times$1.1(-12)&-6.4(-13)+i$\times$3.5(-13)\\
9&6.4(-1)-i$\times$2.3(-1)&-2.5(-8)+i$\times$3.1(-8)&-4.0(-5)+i$\times$1.6(-5)&9.3(-13)+i$\times$1.1(-12)&-6.4(-13)-i$\times$3.5(-13)\\
10&-9.0(-4)+i$\times$6.6(-4)&-8.3(-4)+i$\times$6.6(-4)&-9.0(-4)+i$\times$6.6(-4)&1.2(-4)-i$\times$8.8(-5)&5.2(-5)-i$\times$3.9(-5)\\
11&-9.0(-4)-i$\times$6.6(-4)&2.0(-3)-i$\times$6.6(-4)&4.0(-4)-i$\times$6.6(-4)&1.2(-4)+i$\times$8.8(-5)&5.2(-5)+i$\times$3.9(-5)\\
12&3.5(-4)+i$\times$1.1(-3)&7.1(-4)-i$\times$1.2(-4)&-8.1(-1)+i$\times$1.1(-3)&-4.7(-5)-i$\times$1.4(-4)&-2.1(-5)-i$\times$6.1(-5)\\
13&3.5(-4)-i$\times$1.1(-3)&7.1(-4)+i$\times$1.2(-4)&3.5(-4)-i$\times$1.1(-3)&-4.7(-5)+i$\times$1.4(-4)&-2.1(-5)+i$\times$6.1(-5)\\
\Xhline{1.0pt}
\end{tabular}
\label{table_2}
\end{table*}

\end{Example}

All test sets in the following two examples are randomly
generated by the ``\verb|randn|" command in MATLAB,
where  $\mathfrak{L}$ contains some repeated poles (real or non-real).

\begin{Example}\label{example1}
This example consists of two test sets. The first test set, which is to illustrate
the performance of all methods when repeated poles are all real,
contains $70$ random examples, where $n$ varies from 3 to 13 increased by 2,
and $m$ is set to be $2, \lfloor\frac{n}{2}\rfloor, n-1$ for each $n$.
For each fixed $(n, m)$, the greatest multiplicity $a_{max}$ of all
real poles increases from $1$ to $m$ in increment of $1$.
All examples are generated as follows. We first randomly generate a
nonsingular matrix $Y\in\mathbb{R}^{n\times n}$ and
$B\in\mathbb{R}^{n\times m}, F\in\mathbb{R}^{m\times n}$ by the MATLAB function \verb|randn|
and the assigned poles
$\mathfrak{L}=\{ \verb|randn|\times \verb|ones|(1, a_{max}), \verb|randn|(1, n-a_{max})\}$,
then set $A=Y\Lambda Y^{-1}-BF$, where the diagonal elements of the diagonal matrix $\Lambda$ are those in $\mathfrak{L}$.
Taking $A,B$ and $\mathfrak{L}$ as the input, we apply the methods \verb|place|, \verb|robpole|,
\verb|Schur-rob| and \verb|Schur-multi| to these examples,
where the poles are assigned in ascendant order.

For concision, we only list results for $n=13$.
Results for other examples are quite similar. Specifically, Fig.~\ref{fig1} to  Fig.~\ref{fig4}
show the three measures of robustness and the precision of the computed poles by all four methods,
and Fig.~\ref{fig5} plots the ratios of the CPU time costs of  \verb|place|, \verb|robpole|
and \verb|Schur-rob| with respect to that of \verb|Schur-multi|.
In each figure, the three subfigures correspond to $m=2, 6$ and $12$, respectively.
The $x$-axis represents $a_{max}$, and the values in the $y$-axis are
mean values over $50$ trials for a certain triple $(13, m, a_{max})$.

On these examples, our method is comparable with \verb|place| and \verb|robpole|,
but with much less time cost. Comparing with \verb|Schur-rob|, \verb|Schur-multi| does
improve the relative accuracy of the assigned poles when some poles to be assigned are repeated and real.

\begin{figure*}[!t]
\centering
\begin{tabular}{ccc}
    \begin{minipage}[t]{0.3\textwidth}
    \includegraphics[width=2.2in]{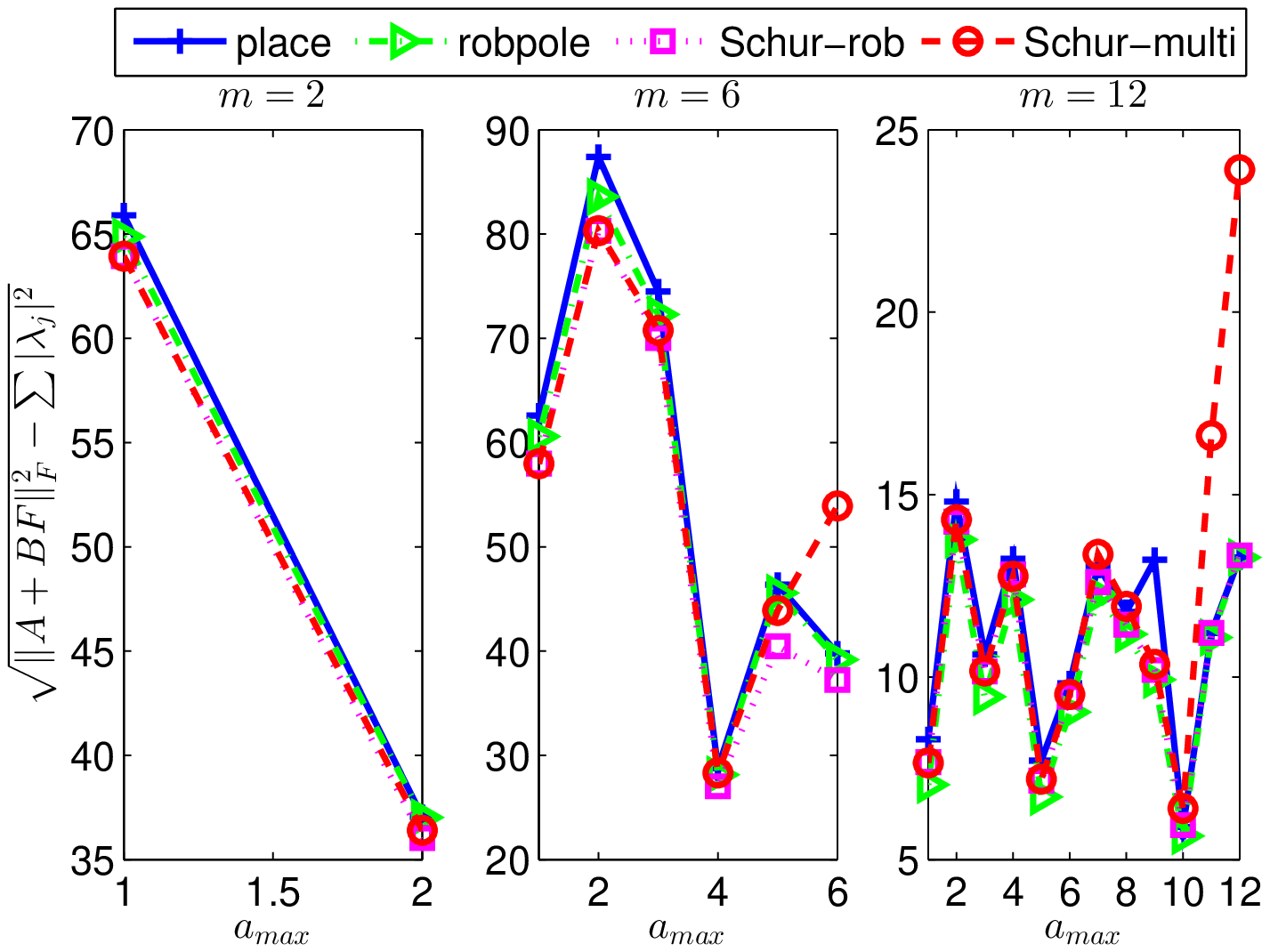}
    \caption{$dep.$ (Example \ref{example1} with real repeated poles)}\label{fig1}
    \end{minipage}&
    \begin{minipage}[t]{0.3\textwidth}
    \includegraphics[width=2.2in]{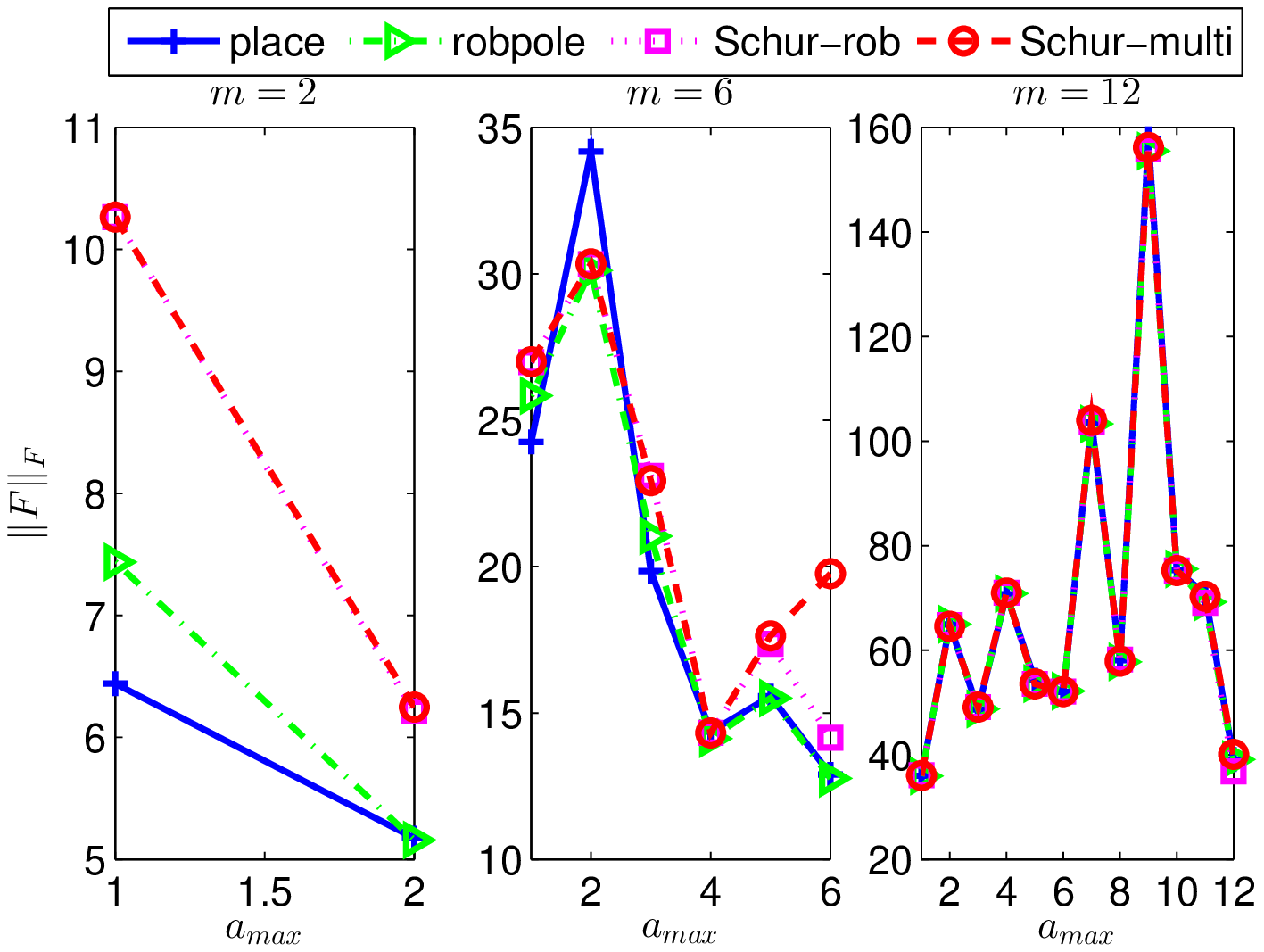}
    \caption{$\|F\|_F$ (Example \ref{example1} with real repeated poles)}\label{fig2}
    \end{minipage}&
    \begin{minipage}[t]{0.3\textwidth}
    \includegraphics[width=2.2in]{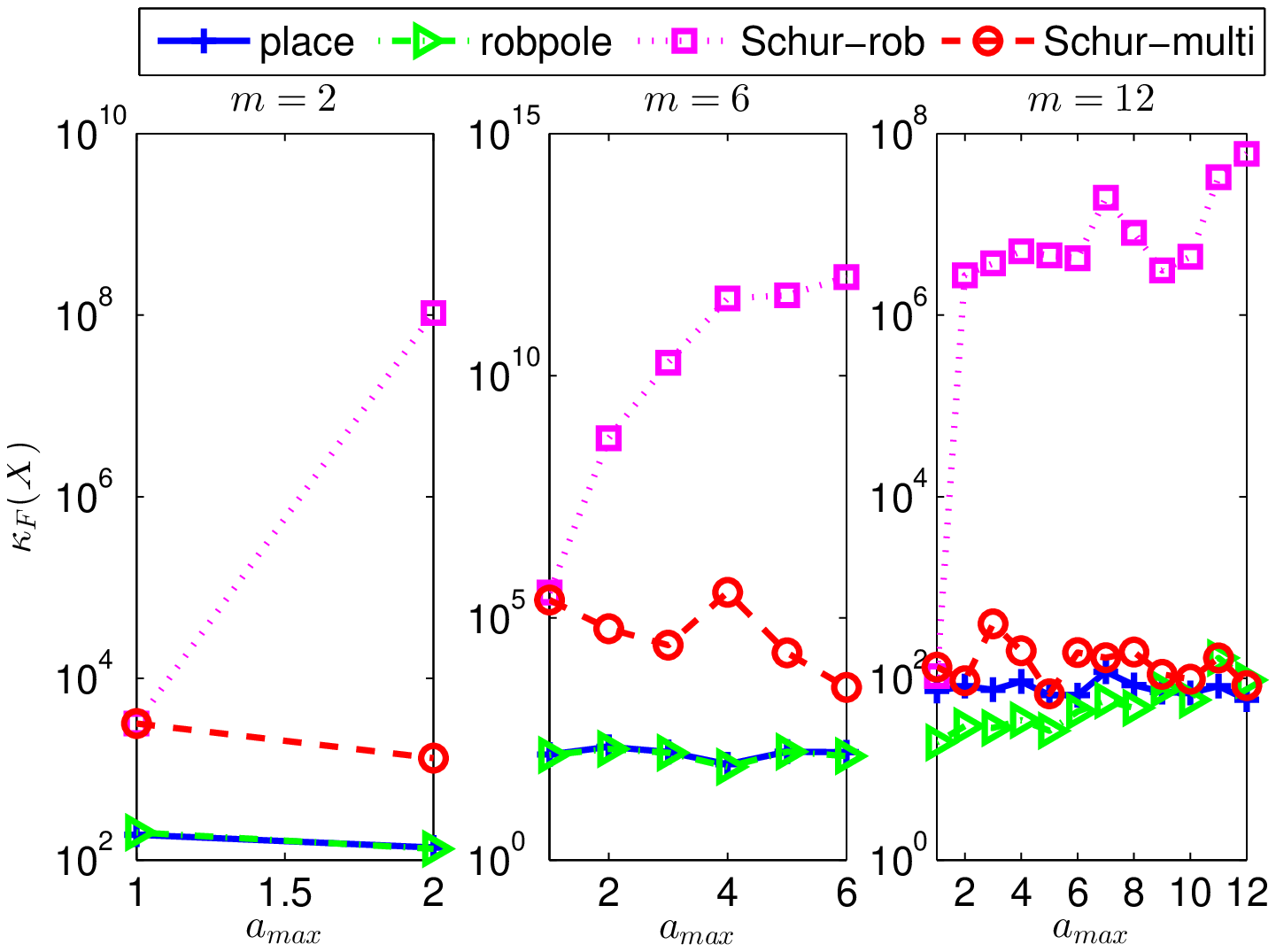}
    \caption{$\kappa_F(X)$ (Example \ref{example1} with real repeated poles)}\label{fig3}
    \end{minipage}
\end{tabular}
\end{figure*}

\begin{figure*}[!t]
\centering
\begin{tabular}{cc}
    \begin{minipage}[t]{0.4\textwidth}
    \includegraphics[width=2.2in]{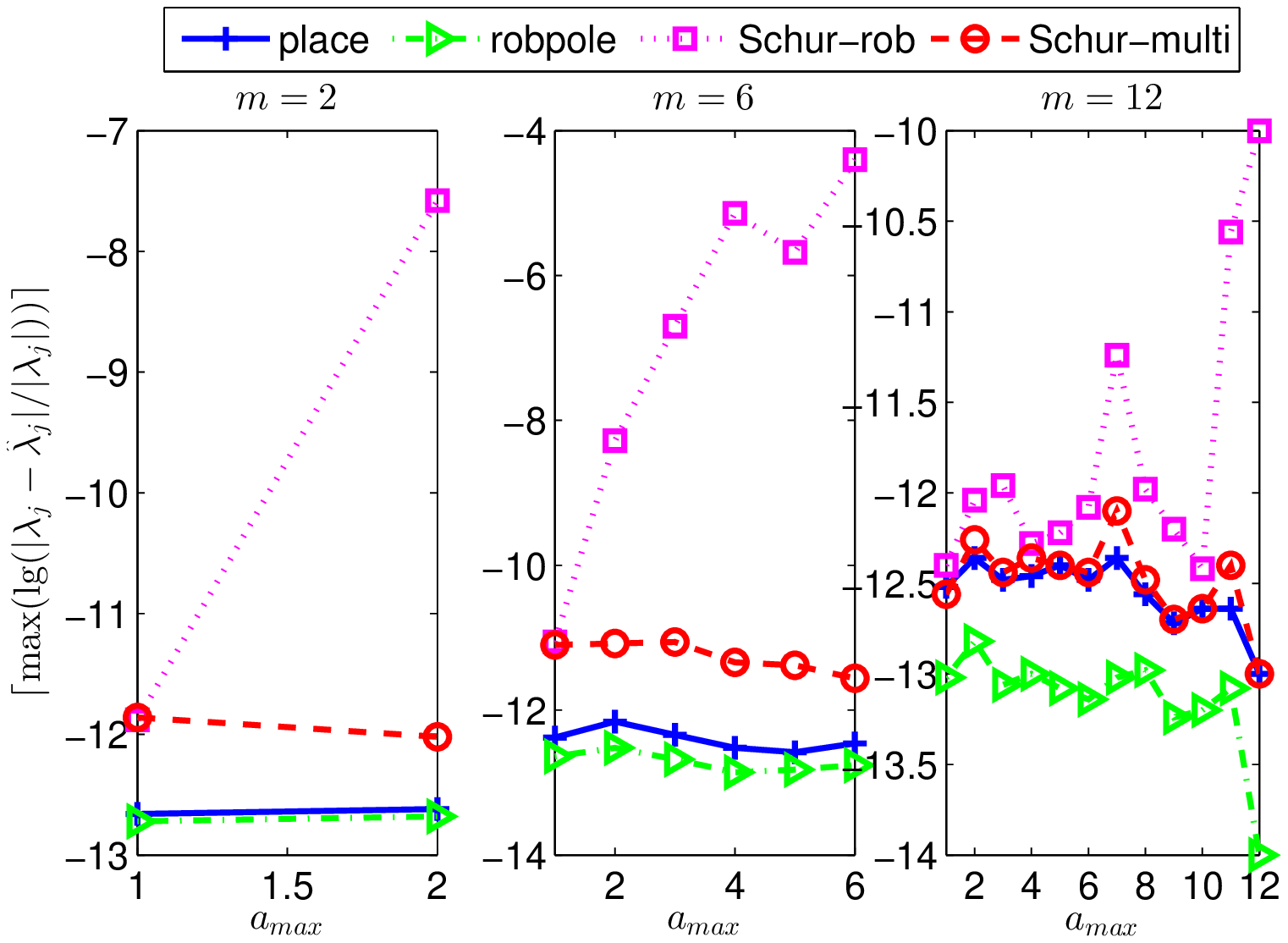}
    \caption{$precs$ (Example \ref{example1} with real repeated poles)}\label{fig4}
    \end{minipage}&
    \begin{minipage}[t]{0.4\textwidth}
    \includegraphics[width=2.2in]{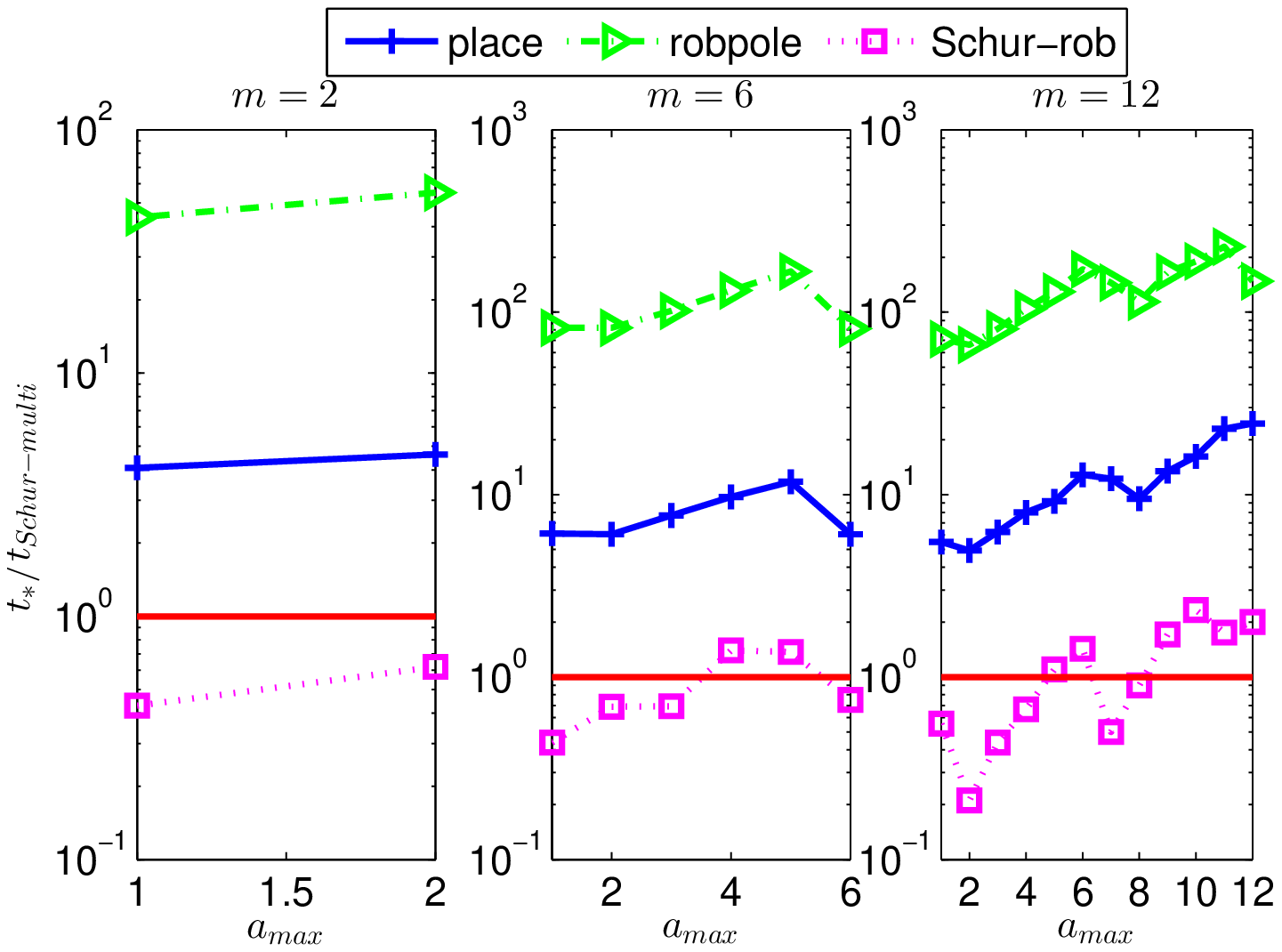}
    \caption{CPU time ratio (Example \ref{example1} with real repeated poles)}\label{fig5}
    \end{minipage}
\end{tabular}
\end{figure*}




\smallskip

The second test set consists of $82$ random examples,
which is to demonstrate the performance of all methods when non-real repeated poles are contained in $\mathfrak{L}$.
Here, we take $n$ varying from $7$ to $19$ with an increment of $2$,
and $m$ is set to be $3, \lfloor\frac{n}{2}\rfloor, n-1$ for each $n$.
For fixed $(n, m)$, the largest multiplicity $a_{max}$ of
all complex poles increases from $2$ to $\min\{\lfloor\frac{n}{2}\rfloor, m\}$.
All examples are generated as follows.
First, we randomly generate the placed poles
$\mathfrak{L}=\{ \verb|randn|(1, n-2a_{max}), \lambda\times \verb|ones|(1, a_{max}),
\bar{\lambda}\times \verb|ones|(1, a_{max})\}$
with $\lambda=\verb|randn|+i\times \verb|randn|$, and three matrices $Y\in\mathbb{R}^{n\times n}$,
$B\in\mathbb{R}^{n\times m}, F\in\mathbb{R}^{m\times n}$ using the MATLAB function \verb|randn|.
Compute the QR decomposition of $Y$ as $Y=Q_YR_Y$, and we reset the diagonal and subdiagonal entries of
$R_Y$ such that it is upper quasi-triangular with its eigenvalues being those in $\mathfrak{L}$.
Then set $A=Q_YR_YQ^{\top}_Y-BF$.
Thereafter,  the algorithms \verb|place|, \verb|robpole|, \verb|Schur-rob|
and \verb|Schur-multi| are applied on all examples with $A,B$ and $\mathfrak{L}$ taken as the input.

Fig.~\ref{fig1c} to Fig.~\ref{fig5c} exhibit the numerical results on $dep.$, $\|F\|_F$, and $\kappa_F(X)$,
$precs$  and the CPU time ratio for $n=19$, respectively,
where the $x$-axis and the $y$-axis own the some  meanings as those in the first test set.
Each figure includes three subfigures,  where the first one displays the results for $m=3$,
the second for  $m=9$ and  the third for $m=18$.
Note that for the CPU time, we still adopt the time cost of \verb|Schur-multi| as
the standard of comparison, and present the ratios of \verb|place|, \verb|robpole|
and \verb|Schur-rob| to it.

All figures show that when $a_{max}$ is no more than $\lfloor\frac{m+1}{2}\rfloor$,
then compared with \verb|robpole|, our approach produces comparable results on
the robustness and the  precision of the assigned poles,
but with much less time consumption.
However, if there exists at least one complex pole with
its  multiplicity being larger than $\lfloor\frac{m+1}{2}\rfloor$,
the closed-loop system matrix obtained by \verb|Schur-multi| can not be diagonalized and it would not be
as robust as that computed by \verb|robpole|.
Notice that for our  \verb|Schur-multi| method,
there are sharp jumps in Fig.~\ref{fig1c} and Fig.~\ref{fig2c} for $m = 9, \ 18$ cases,
where $a_{max}=\lfloor\frac{m+1}{2}\rfloor$. And the explanation for those jumps is:
$\lfloor\frac{m+1}{2}\rfloor$  actually is  a threshold  that distinguishes
if the repeated non-real pole acts as a semi-simple eigenvalue or not,
hence  those  repeated complex poles, whose multiplicities
equal to $\lfloor\frac{m+1}{2}\rfloor$,  would be more sensitive to perturbations;
and such behavior eventually reflects in $dep.$ and $\|F\|_F$.
In addition, compared with \verb|Schur-rob|, \verb|Schur-multi| does make some improvements
on the precision of the assigned repeated complex conjugate poles.
The undisplayed results for other  different $n$  show similar behavior.

\begin{figure*}[!t]
\centering
\begin{tabular}{ccc}
    \begin{minipage}[t]{0.3\textwidth}
    \includegraphics[width=2.2in]{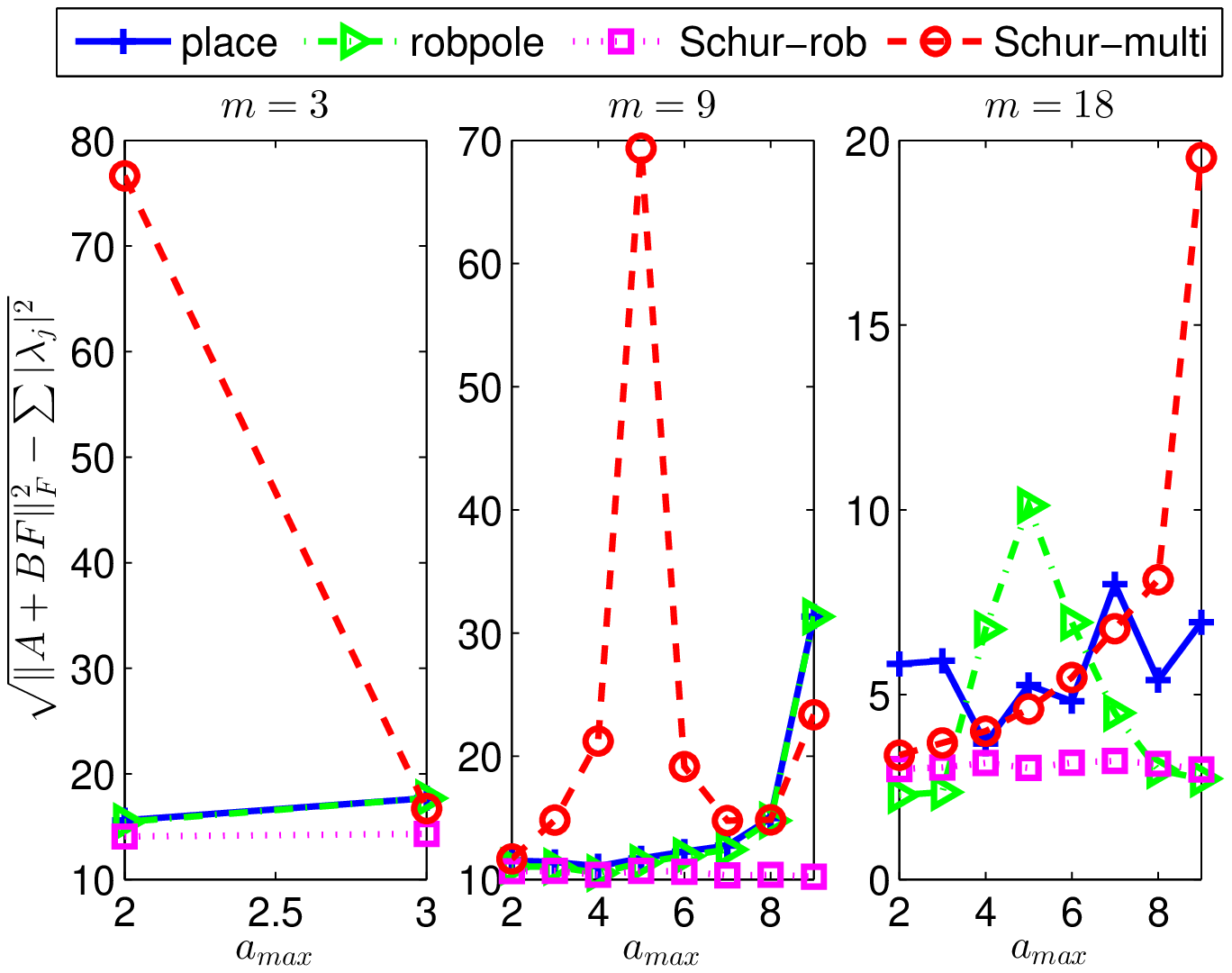}
    \caption{$dep.$ (Example \ref{example1} with non-real  repeated poles)}\label{fig1c}
    \end{minipage}&
    \begin{minipage}[t]{0.3\textwidth}
    \includegraphics[width=2.2in]{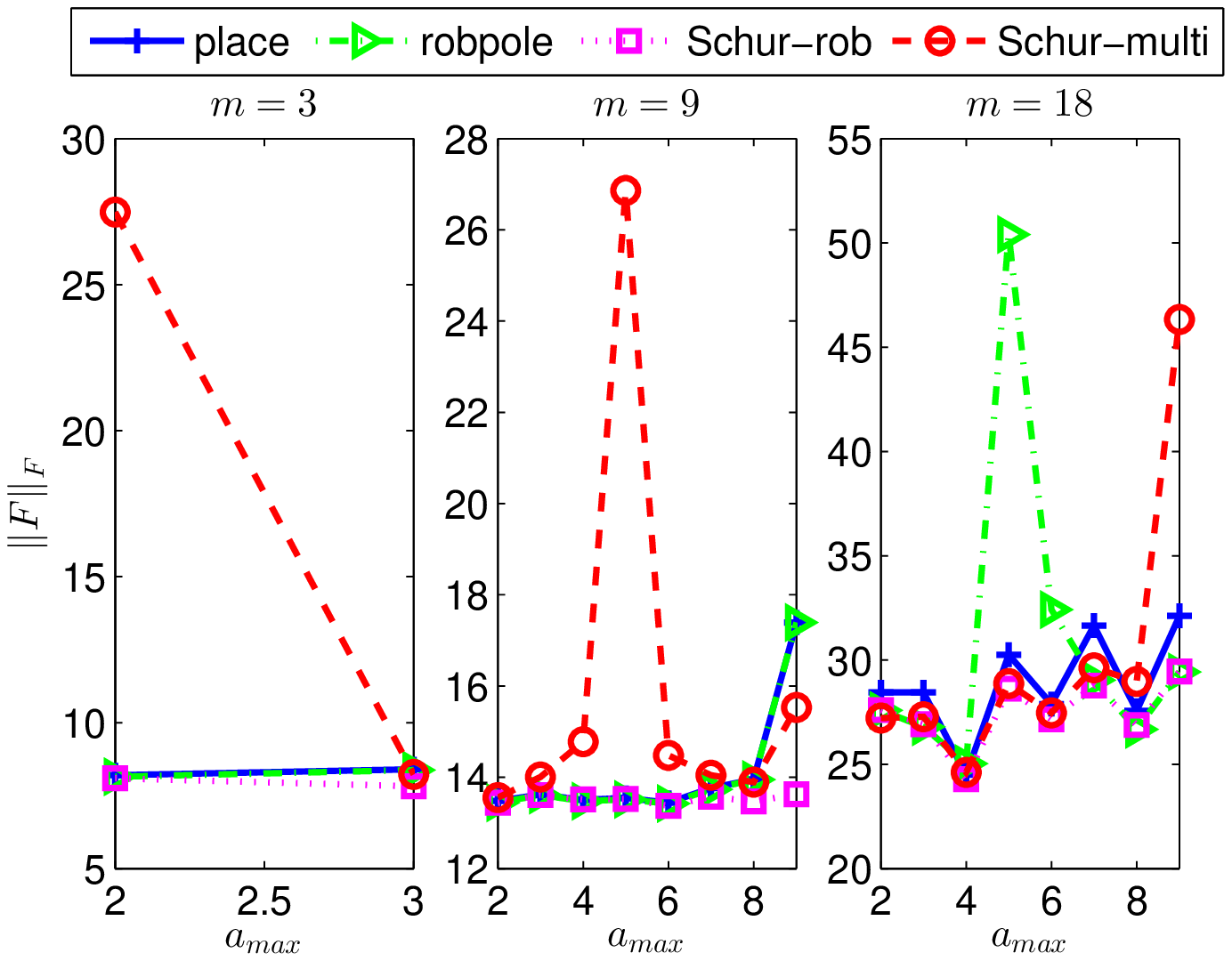}
    \caption{ $\|F\|_F$ (Example \ref{example1} with non-real repeated poles)}\label{fig2c}
    \end{minipage}&
    \begin{minipage}[t]{0.3\textwidth}
    \includegraphics[width=2.2in]{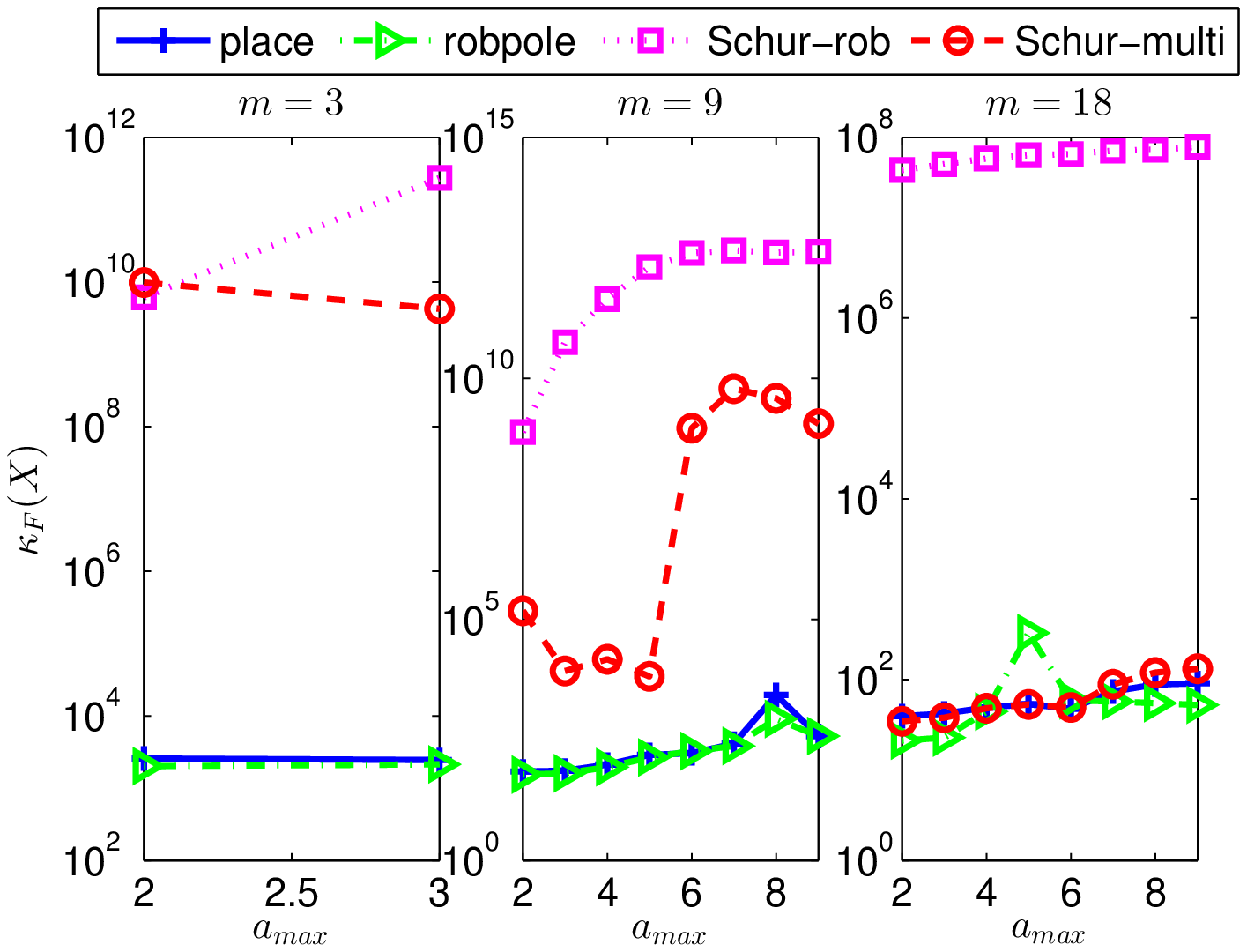}
    \caption{ $\kappa_F(X)$ (Example \ref{example1} with non-real repeated poles)}\label{fig3c}
    \end{minipage}
\end{tabular}
\end{figure*}

\begin{figure*}[!t]
\centering
\begin{tabular}{cc}
    \begin{minipage}[t]{0.4\textwidth}
    \includegraphics[width=2.2in]{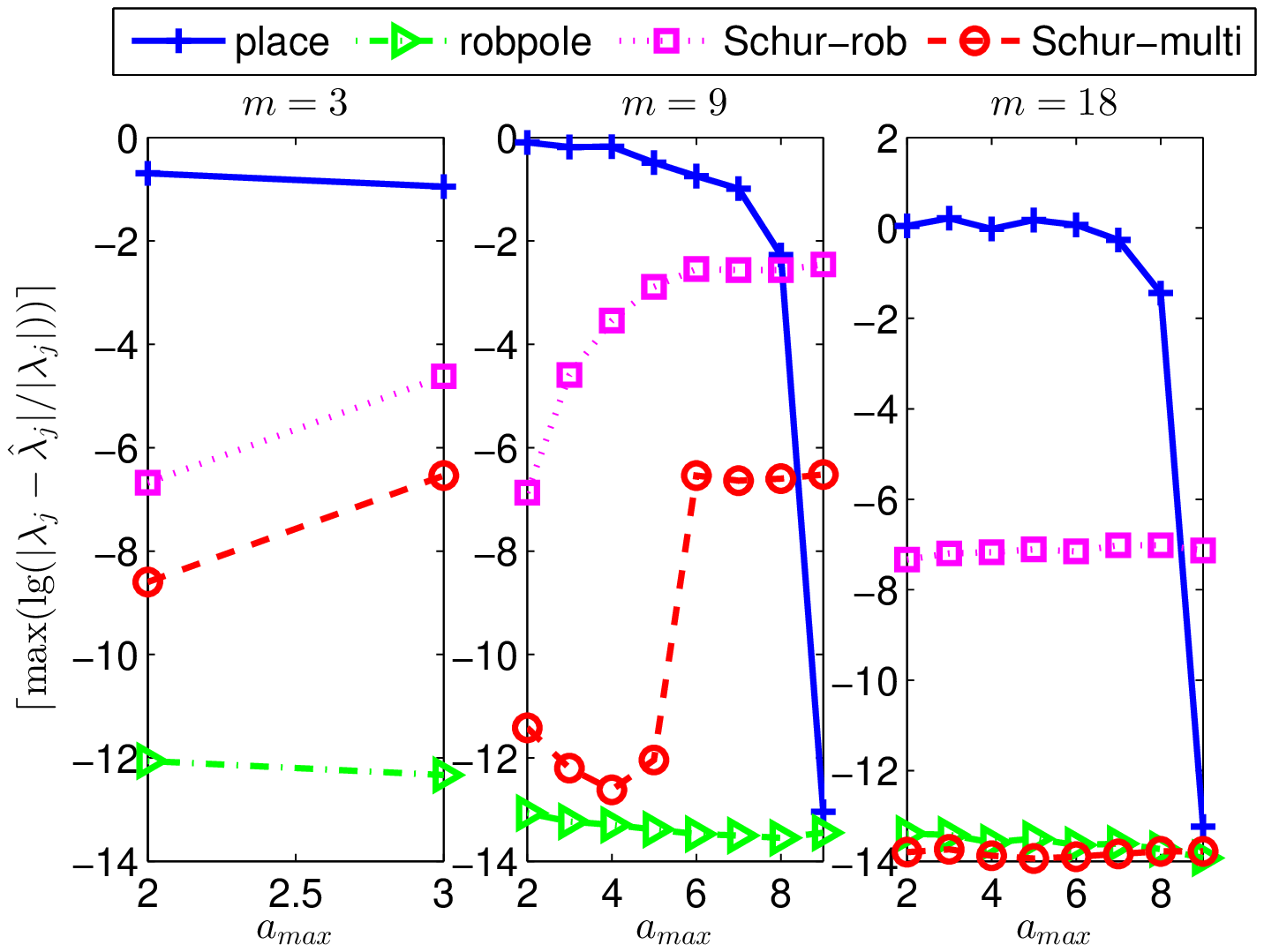}
    \caption{$precs$ (Example \ref{example1} with non-real repeated poles)}\label{fig4c}
    \end{minipage}&
    \begin{minipage}[t]{0.4\textwidth}
    \includegraphics[width=2.2in]{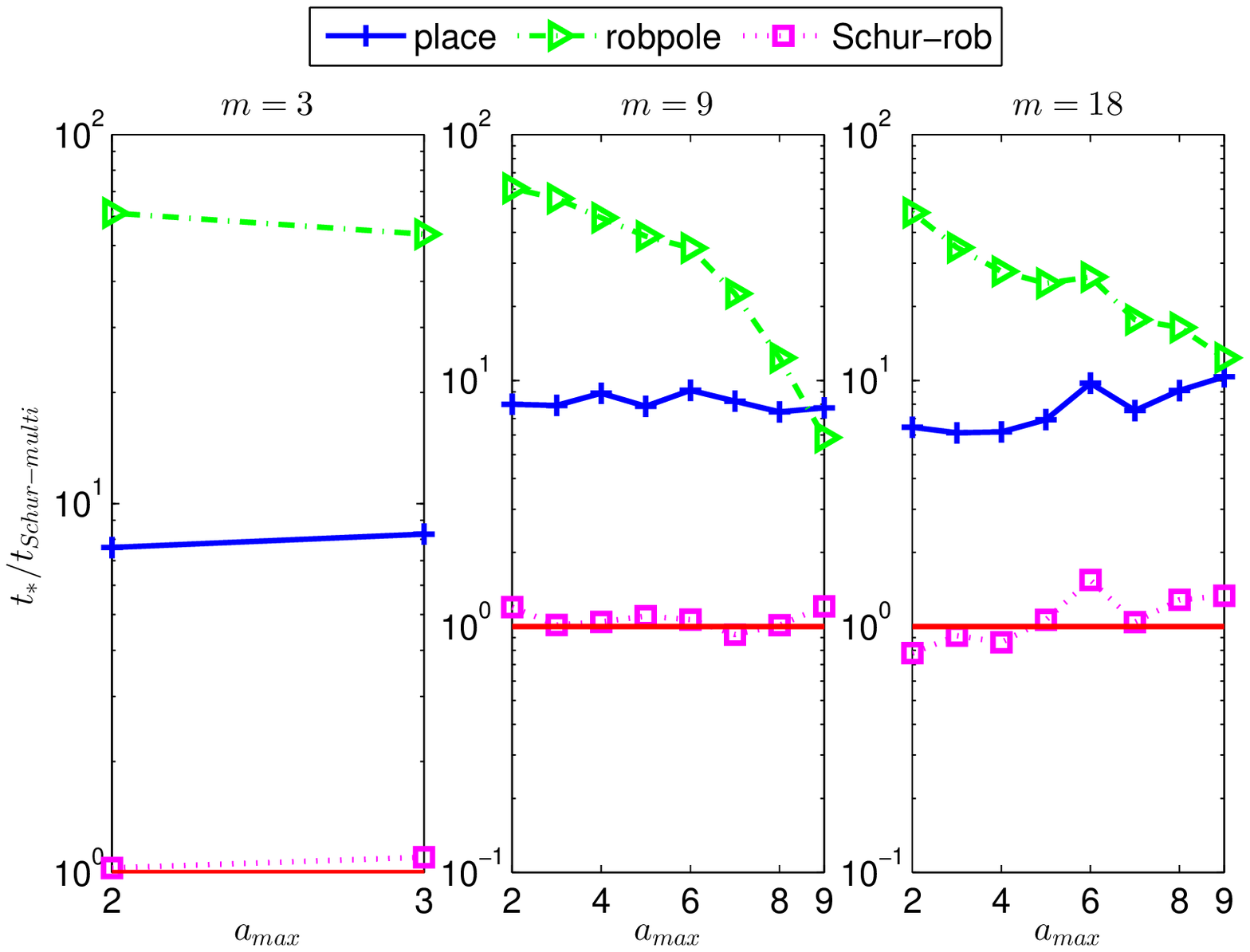}
    \caption{CPU time ratio (Example \ref{example1} with non-real repeated poles)}\label{fig5c}
    \end{minipage}
\end{tabular}
\end{figure*}

%
%

\end{Example}

\bigskip

It is well known that \verb|place| and \verb|robpole| can not solve the
{\bf SFRPA} if the  multiplicity of some pole is greater than $m$,
while \verb|Schur-rob| and our \verb|Schur-multi| can still work.
The following randomly generated examples are to reveal the behavior of
\verb|Schur-rob| and \verb|Schur-multi| on examples in which the multiplicity of some repeated pole might be greater than $m$.

\begin{Example}\label{example2}
This example also consists of two test sets. The first test set, where the repeated poles are all real,
is comprised of $270$ random examples with $n$ increasing from $7$ to $27$ in increment of $4$,
and $m$ being $2, \lfloor\frac{n}{2}\rfloor, n-1$ for each $n$.
For fixed $(n, m)$, the greatest  multiplicity of the assigned  repeated real poles $a_{max}$ varies from $2$ to $n-1$.
All examples are generated as below. We first randomly generate the assigned poles
$\mathfrak{L}=\{ \verb|randn|\times \verb|ones|(1, a_{max}), \verb|randn|(1, n-a_{max})\}$ and $Y\in\mathbb{R}^{n\times n}$,
$B\in\mathbb{R}^{n\times m}, F\in\mathbb{R}^{m\times n}$ by the MATLAB function \verb|randn|.
Then we compute the  QR decomposition of $Y$ as $Y=Q_YR_Y$,
reset the diagonal elements of the upper triangular matrix $R_Y$
be those in $\mathfrak{L}$, and set $A=Q_YR_YQ^{\top}_Y-BF$.
Taking $A, B$ and $\mathfrak{L}$ as the input, we then apply \verb|Schur-rob|
and \verb|Schur-multi| to all generated examples.
The poles in $\mathfrak{L}$ are also assigned in  ascendant order.
Note that when applying \verb|place| and \verb|robpole| on these examples, they fail to give results for some examples.
For instance, when $m=2$ and $a_{max}>2=m$, they fail to output solutions.

Both algorithms produce fairly similar $dep.$ and  $\|F\|_F$, and we omit the interrelated results here.
The  numerical results on $\kappa_F(X)$ and $precs$ with respect to $a_{max}$ for $n=19$ are displayed in Fig.~\ref{fig8} and  Fig.~\ref{fig9}, respectively, where the $x$-axis and $y$-axis own the same meanings as those in Example \ref{example1}.
In each figure, the three subfigures correspond to $m=2, 9$ and $18$, respectively.

From Fig.~\ref{fig8} and Fig.~\ref{fig9},  we know that the condition numbers of the
eigenvectors matrices obtained by \verb|Schur-multi| are smaller than those by \verb|Schur-rob|,
and the eigenvalues of $A_c$ computed by \verb|Schur-multi| are more accurate than those by \verb|Schur-rob|.
The differences become more significant when $a_{max}$ is no greater than $m$. If $a_{max}$ is greater than $m$, that is, some eigenvalues of $A_c$ are defective, the precision of the  poles diminishes.
For other $(n, m, a_{max})$, $\kappa_F(X)$ and $precs$ show quite similar variation tendency.

\smallskip

It is shown in Subsection \ref{subsection3.1.1} that if the repeated  real pole with multiplicity $a_{max}$ is assigned as the initial $\lambda_1$, then its geometric multiplicity is theoretically $\min\{m, a_{max}\}$.
However, if it is not assigned foremost, we cannot prove such result in theory.
We then compute the geometric multiplicity (denoted as $``g_{multi}"$)
of the repeated  real  pole by using the SVD of $(A_c-\lambda I_n)$,
where $A_c$ is the computed closed-loop system matrix and $\lambda\in\mathfrak{L}$.
Note that in our experiments, the poles are assigned in ascendant order.
That is, the  repeated  real pole may not be the first one to be placed. However, the numerical
results for $n=19$ listed in TABLE \ref{table_0} show that
$g_{multi}$ obtained by \verb|Schur-multi| always equals to $\min\{m,a_{max}\}$.
The unshown results for other different $(n, m, a_{max})$ behave similarly.


\begin{figure*}[!t]
\centering
\begin{tabular}{ccc}
    \begin{minipage}[t]{0.3\textwidth}
    \includegraphics[width=2.2in]{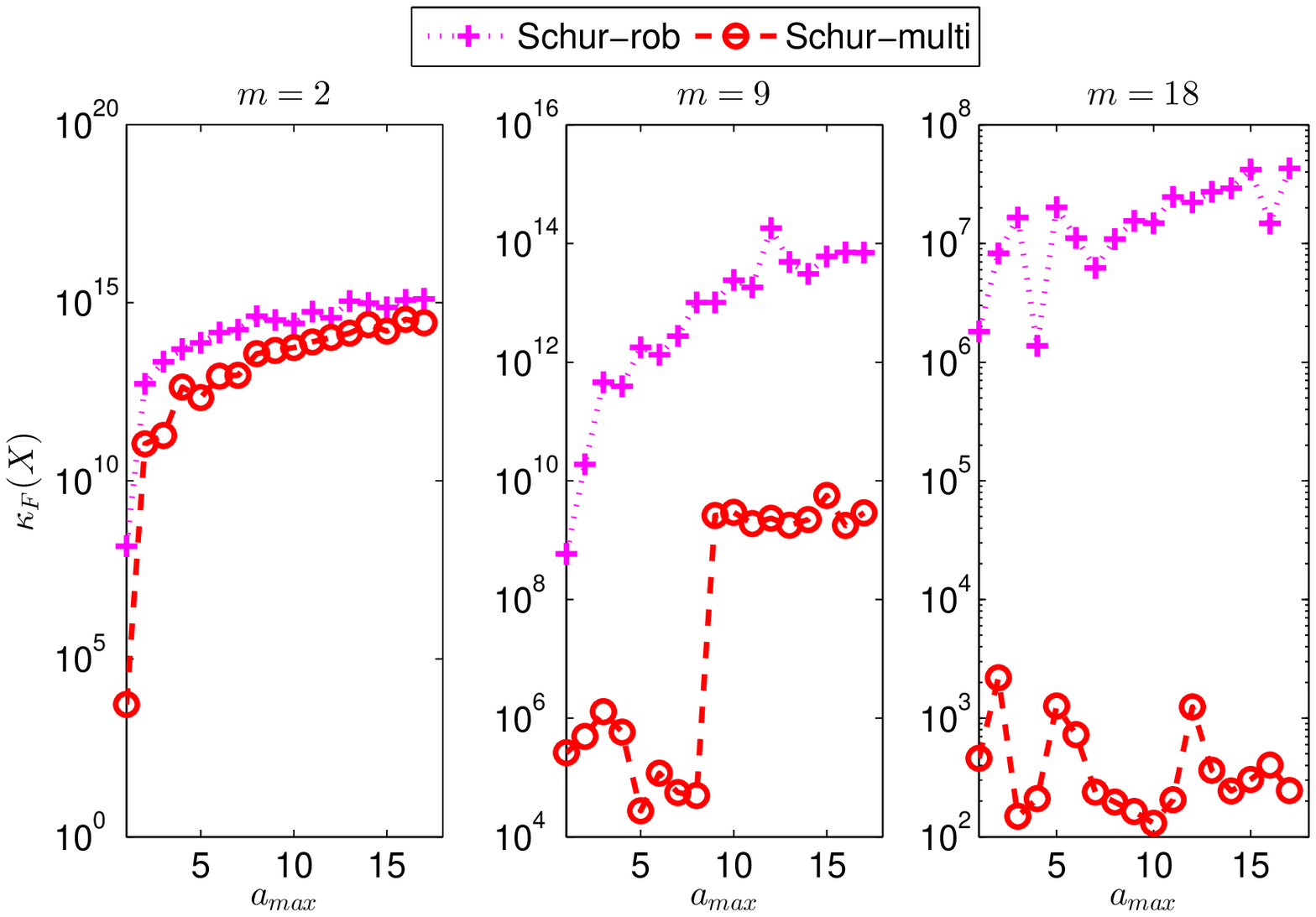}
    \caption{ $\kappa_F(X)$ (Example \ref{example2} with real repeated poles)}\label{fig8}
    \end{minipage}&
    \begin{minipage}[t]{0.3\textwidth}
    \includegraphics[width=2.2in]{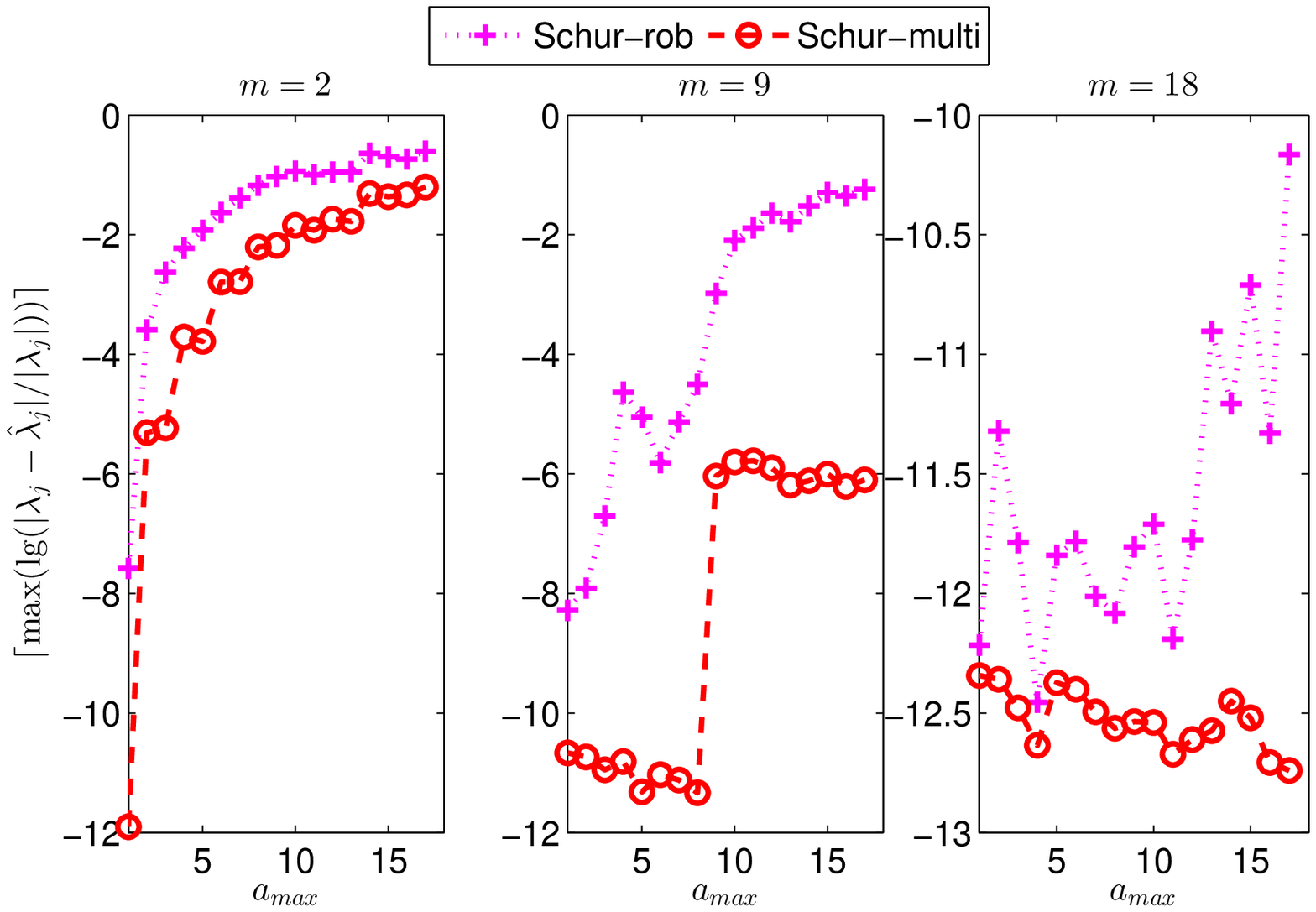}
    \caption{$precs$ (Example \ref{example2} with real repeated poles)}\label{fig9}
    \end{minipage}&
    \begin{minipage}[t]{0.3\textwidth}
    \includegraphics[width=2.2in]{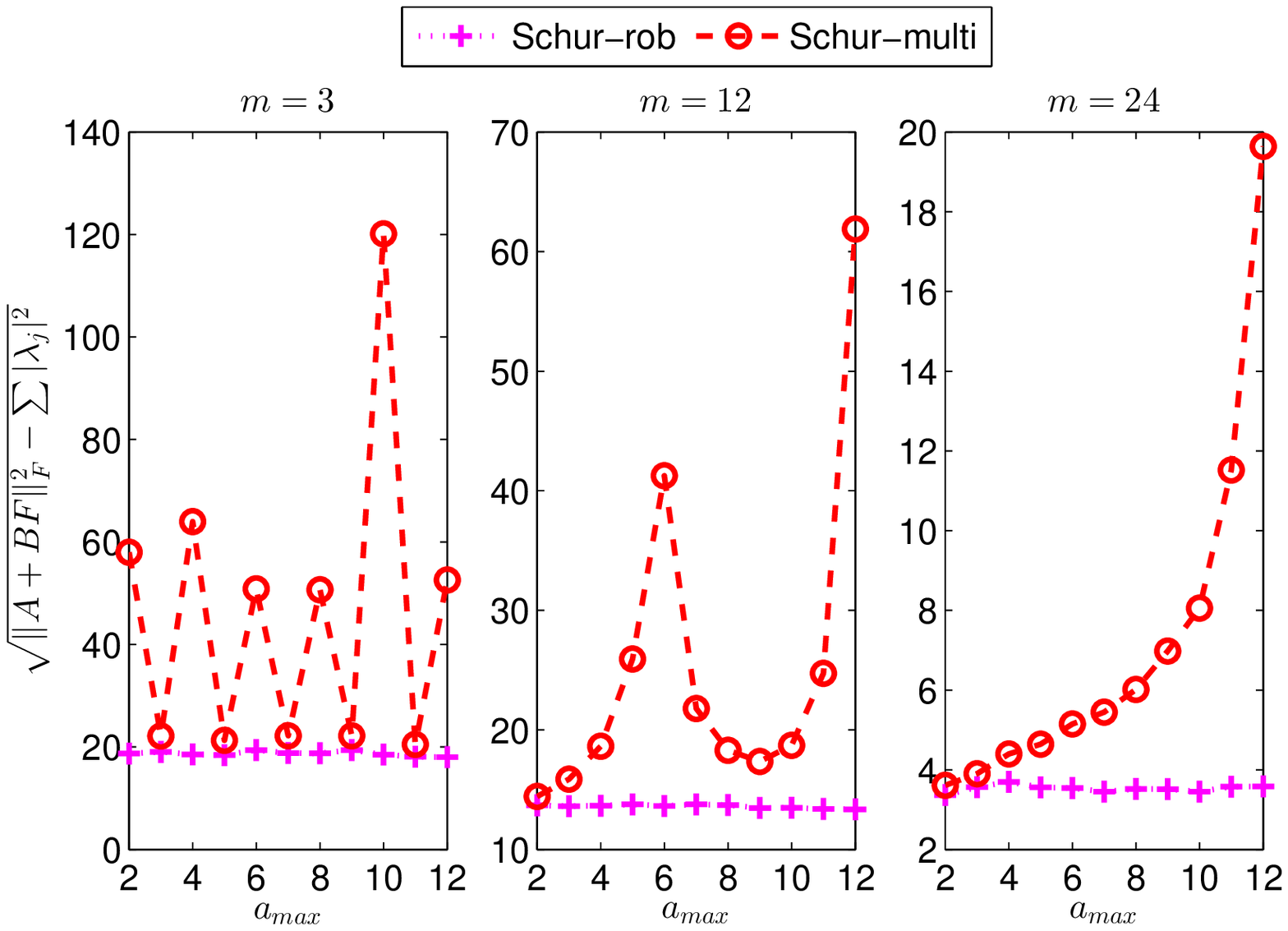}
    \caption{$dep.$ (Example \ref{example2} with non-real repeated poles)}\label{fig6c}
    \end{minipage}
\end{tabular}
\end{figure*}



\tabcolsep 0.1in
\begin{table*}
\small
\centering
\caption{ Geometric multiplicity over $50$ trials (real repeated poles)  }
\begin{tabular}{c|cc|cc|cc}
\Xhline{1.0pt}
\multirow{3}*{} & \multicolumn{6}{c}{$g_{multi}$ \ \text{for} \  $n=19$}\\
\Xcline{2-7}{0.7pt}
 &  \multicolumn{2}{c|}{$m=2$}&
\multicolumn{2}{c|}{{$m=\lfloor\frac{n}{2}\rfloor$}}&\multicolumn{2}{c}{$m=n-1$} \\
\Xcline{2-7}{0.5pt}
$a_{max}$ & \verb|Schur-rob|& \verb|Schur-multi|& \verb|Schur-rob|& \verb|Schur-multi|& \verb|Schur-rob|& \verb|Schur-multi|\\
\Xhline{0.7pt}
2 &     1.04 &   2.00 &  1.44 &   2.00 &   1.96  &  2.00\\
3&     1.06&   2.00  &  2.10 &   3.00 &   2.80  &  3.00\\
4&     1.04 &   2.00  &  2.44 &   4.00 &   3.86 &  4.00\\
5&     1.06 &   2.00  &  2.22 &   5.00&   4.98 &  5.00\\
6&     1.06 &   2.00  &  2.90 &   6.00 &   5.90  &  6.00\\
7&     1.08 &   2.00  &  4.24 &   7.00 &   6.88  &  7.00\\
8&     1.06 &   2.00  &  4.28 &   8.00&   7.92  &  8.00\\
9&     1.08 &   2.00 &  4.42 &   9.00 &   8.92  &  9.00\\
10&    1.02 &   2.00 &  5.06 &   9.00 &   9.86 & 10.00\\
11&    1.16 &   2.00  &  4.98 &   9.00 &  10.84  & 11.00\\
12&    1.10 &   2.00 &  5.54 &   9.00 &  11.90  & 12.00\\
13&    1.14 &   2.00 &  5.60&   9.00 &  12.84 & 13.00\\
14&    1.14 &   2.00  &  6.66 &   9.00 &  13.70  & 14.00\\
15&    1.20 &   2.00 &  6.62&   9.00 &  14.76  & 15.00\\
16&    1.28 &   2.00 &  7.78 &   9.00 &  15.66  & 16.00\\
17&    1.30 &   2.00  &  8.20 &   9.00 &  16.66  & 17.00\\
18&    1.44 &   2.00  &  8.46&   9.00 &  17.32  & 18.00\\
\Xhline{1.0pt}
\end{tabular}
\label{table_0}
\end{table*}

\smallskip

All numerical examples in the second  test set are designed to illustrate the behavior of
both Schur-type approaches when $\mathfrak{L}$ contains some repeated complex conjugate poles with their multiplicities exceeding $m$.
There are $193$ random illustrative examples in this test set,
with $n$ increasing from $7$ to $25$ in an increment of $2$,
and $m$ taking $3, \lfloor\frac{n}{2}\rfloor, n-1$ for each $n$.
With $(n, m)$ fixed, the largest  multiplicity of the assigned complex
poles varies from $2$ to $\lfloor\frac{n}{2}\rfloor$. All these examples
are generated in the same way as those in the second test set in Example~\ref{example1}.
Regarding $A, B$ and $\mathfrak{L}$ as the input, \verb|Schur-rob|
and \verb|Schur-multi| are then applied to each example.

Here, we just exhibit the numerical results for $n=25$. Numerical results on $dep., \|F\|_F$ and $\kappa_F(X)$
for both algorithms are shown in Fig.~\ref{fig6c} to Fig.~\ref{fig8c}, and
Fig.~\ref{fig9c} displays the relative accuracy $precs$ of the assigned poles.
Each figure includes three subfigures, corresponding to $m=3, 12$ and $24$, respectively.
The $x$-axis and $y$-axis own the same meanings as those in Example \ref{example1}.
From these figures we can see that \verb|Schur-multi| produces slightly worse, but comparable $dep.$ and $\|F\|_F$ as \verb|Schur-rob|,
while $\kappa_F(X)$ and $precs$  produced by \verb|Schur-multi| are much better than those by  \verb|Schur-rob|.
Numerical results for other $n$ behave similarly.

\smallskip

When  the largest  multiplicity of the repeated non-real poles is
larger than $\lfloor\frac{m+1}{2}\rfloor$, for the computed $A_c$ by \verb|Schur-multi|,
there exist defective complex conjugate eigenvalues.
Consequently, the relative accuracy of the placed  repeated complex conjugate poles would be not that high.
To show the  geometric multiplicity (denoted as $``g_{multi}"$) of non-real   repeated eigenvalues of $A_c$ visually,
just as what we do in the first test set, we shall compute it by using  the SVD of $(A_c-\lambda I_n)$,
where $A_c$ is the computed closed-loop system matrix and $\lambda\in\mathfrak{L}$ with $\mbox{Im}(\lambda)\neq 0$.
Typically,  relevant results for $n=25$ are displayed in  TABLE \ref{table_0c},
which shows that $g_{multi}$ obtained from \verb|Schur-multi| equals to the smaller value
between its corresponding algebraic multiplicity and  $\lfloor\frac{m+1}{2}\rfloor$.
The unshown results for other different $(n, m, a_{max})$ are quite similar.

\begin{figure*}[!t]
\centering
\begin{tabular}{ccc}
    \begin{minipage}[t]{0.3\textwidth}
    \includegraphics[width=2.2in]{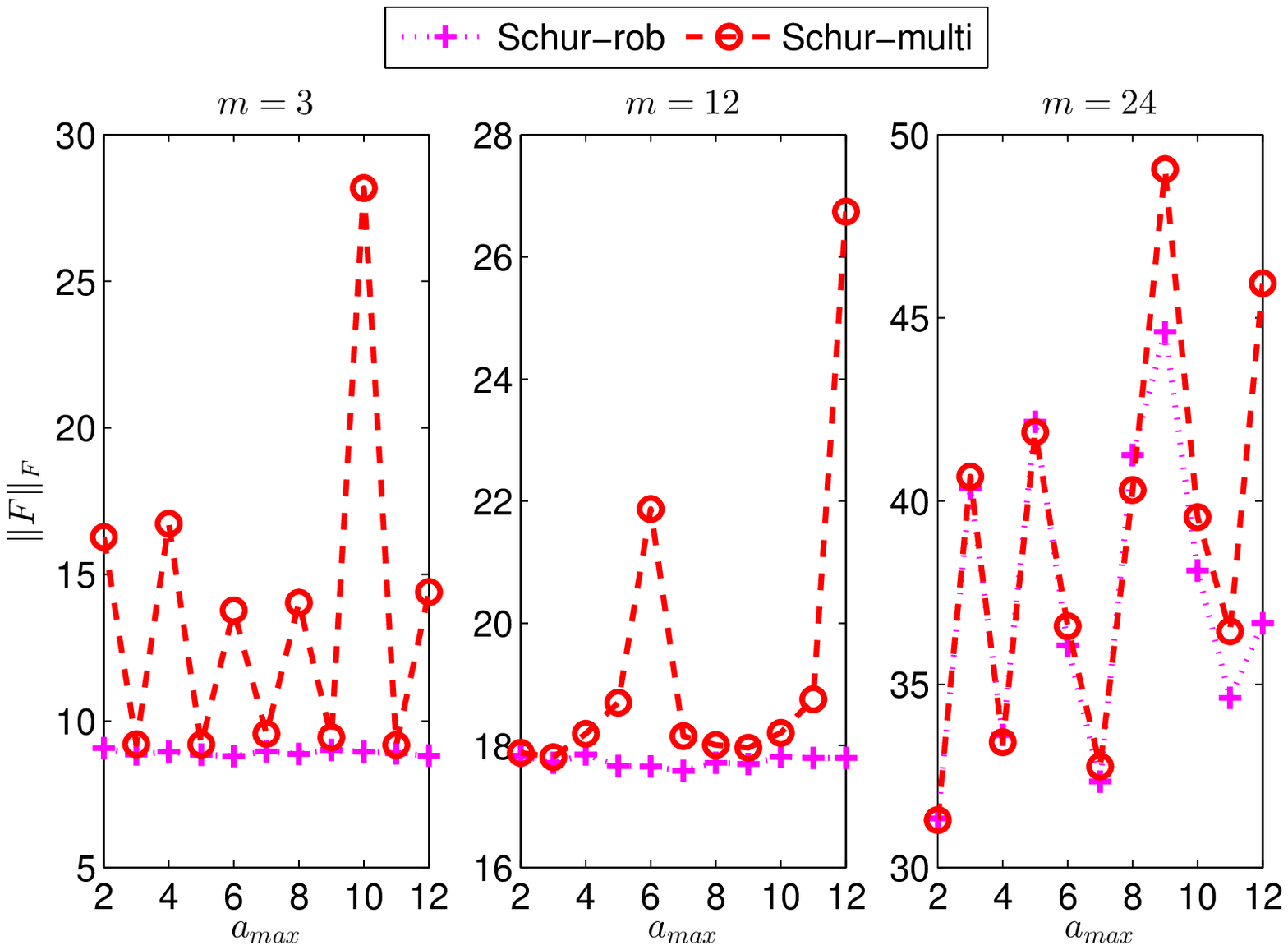}
    \caption{ $\|F\|_F$ (Example \ref{example2} with non-real repeated poles)}\label{fig7c}
    \end{minipage}&
    \begin{minipage}[t]{0.3\textwidth}
    \includegraphics[width=2.2in]{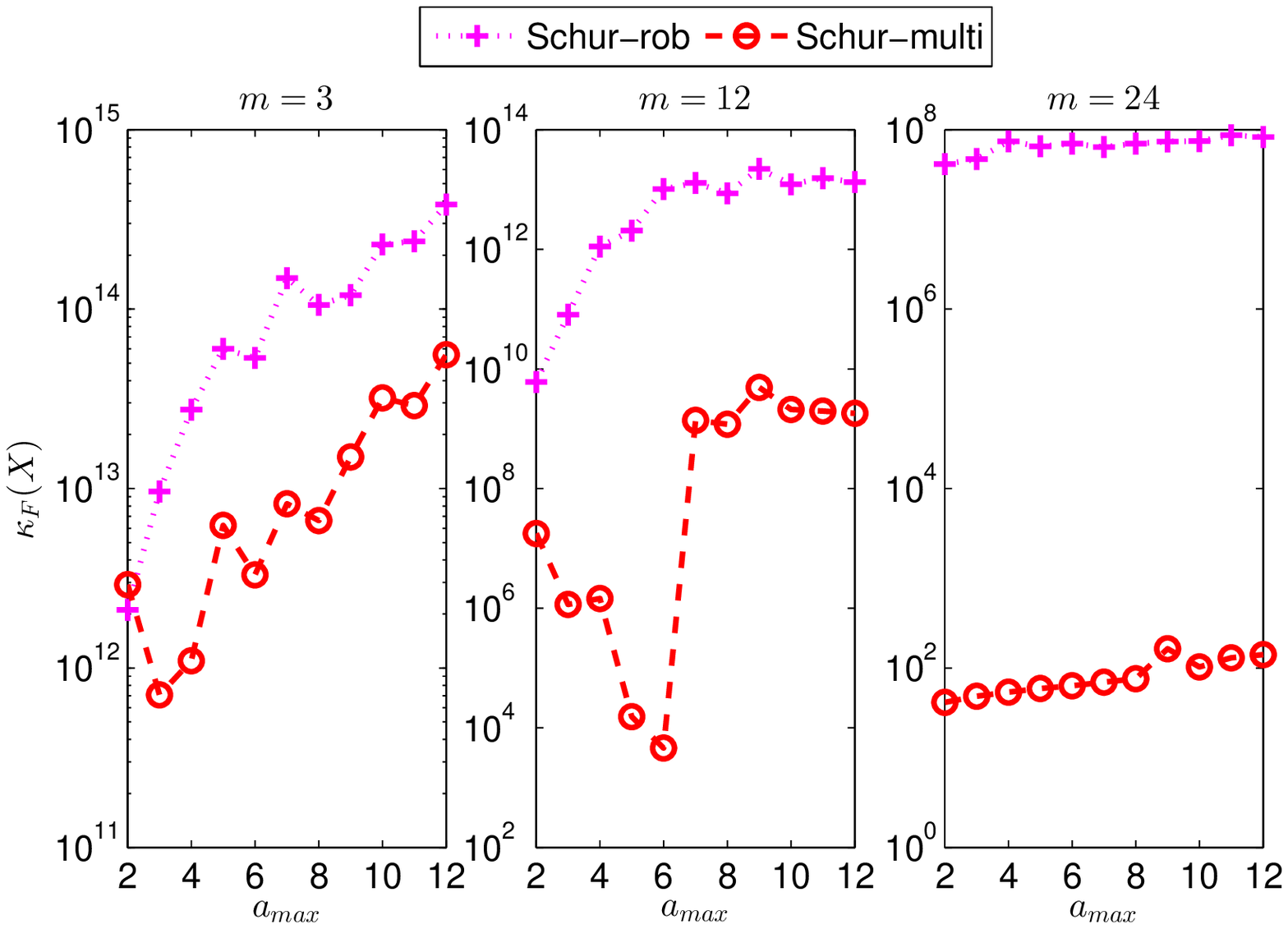}
    \caption{ $\kappa_F(X)$ (Example \ref{example2} with non-real  repeated poles)}\label{fig8c}
    \end{minipage}&
    \begin{minipage}[t]{0.3\textwidth}
    \includegraphics[width=2.2in]{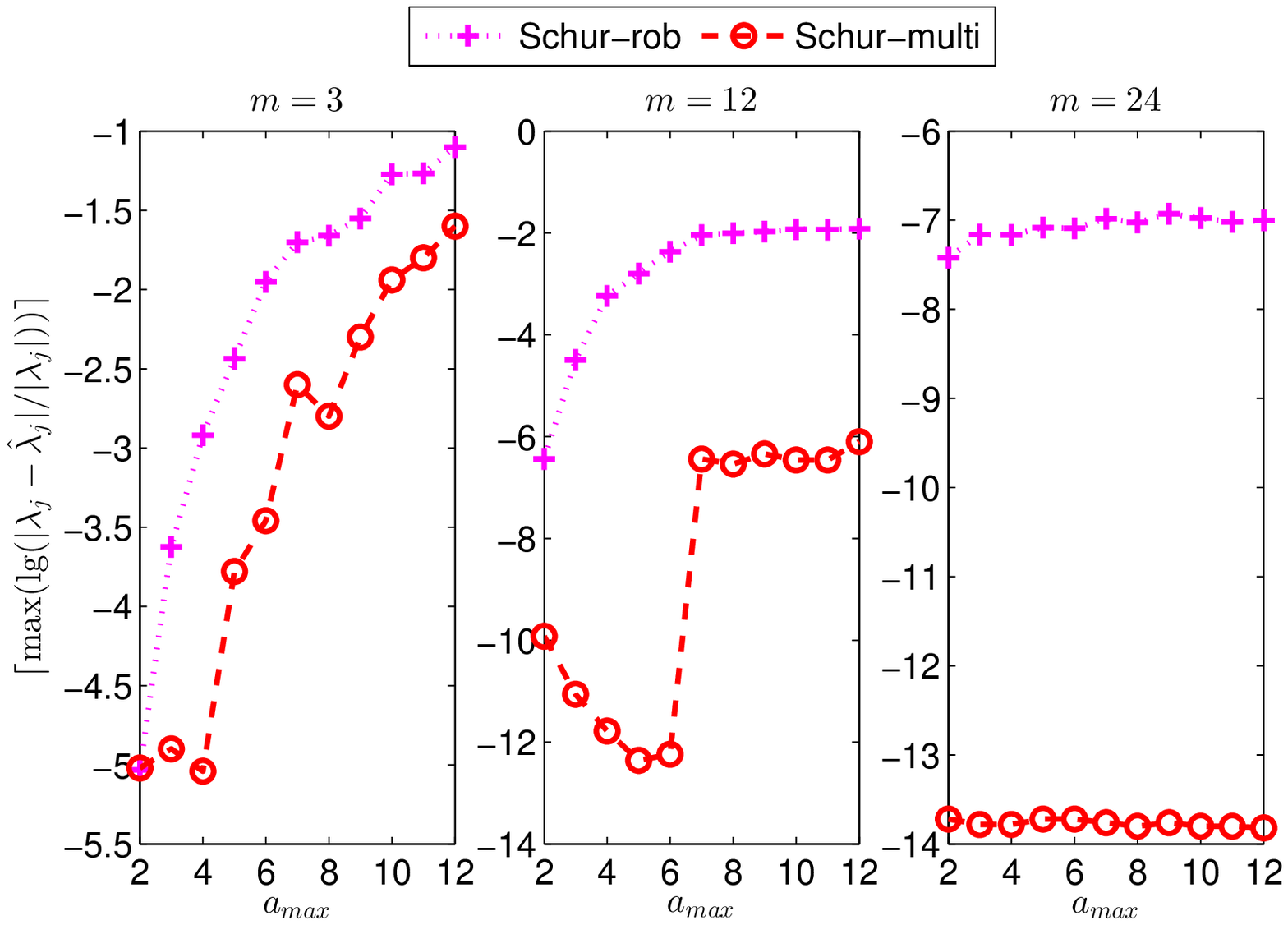}
    \caption{$precs$ (Example \ref{example2} with non-real repeated poles)}\label{fig9c}
    \end{minipage}
\end{tabular}
\end{figure*}

\tabcolsep 0.1in
\begin{table*}
\small
\centering
\caption{ Geometric multiplicity  over $50$ trials (non-real repeated poles)  }
\begin{tabular}{c|cc|cc|cc}
\Xhline{1.0pt}
\multirow{3}*{} & \multicolumn{6}{c}{$g_{multi} \  \text{for} \ n=25$}\\
\Xcline{2-7}{0.7pt}
 &  \multicolumn{2}{c|}{$m=3$}&
\multicolumn{2}{c|}{{$m=\lfloor\frac{n}{2}\rfloor$}}&\multicolumn{2}{c}{$m=n-1$} \\
\Xcline{2-7}{0.5pt}
$a_{max}$ & \verb|Schur-rob|& \verb|Schur-multi|& \verb|Schur-rob|& \verb|Schur-multi|& \verb|Schur-rob|& \verb|Schur-multi|\\
\Xhline{1pt}
 2&       1.00  & 2.00  &  1.00  &  2.00  &  1.00  &  2.00\\
3&       1.00  & 2.00  &  1.00  &  3.00  &  2.00  &  3.00\\
4&       1.00  & 2.00  &  1.00  &  4.00  &  3.00  &  4.00\\
5&       1.00  & 2.00  &  1.00  &  5.00  &  4.00  &  5.00\\
6&       1.00  & 2.00  &  1.00  &  6.00  &  5.00  &  6.00\\
7&       1.00  & 2.00  &  1.00  &  6.00  &  6.00  &  7.00\\
8&       1.00  & 2.00  &  1.06  &  6.00  &  7.00  &  8.00\\
9&       1.00  & 2.00  &  2.04  &  6.00  &  8.00  &  9.00\\
10&      1.00  & 2.00  &  3.08  &  6.00  &  9.00  & 10.00\\
11&      1.00  & 2.00  &  4.04  &  6.00  & 10.00  & 11.00\\
12&      1.00  & 2.00  &  5.16  &  6.00  & 11.00  & 12.00\\
\Xhline{1.0pt}
\end{tabular}
\label{table_0c}
\end{table*}

\end{Example}

\bigskip

\section{Conclusion}\label{section5}
Based on the \verb|Schur-rob| method \cite{GCQX}, a refined approach is proposed  to solve the {\bf SFRPA}, specifically when some poles to be assigned are repeated. In the proposed \verb|Schur-multi| method, we treat the geometric
multiplicities of the repeated poles as the precedential consideration, and
then try to minimize the departure from normality of the closed-loop system matrix $A_c$. Numerical results show that the \verb|Schur-multi| method does outperform the  \verb|Schur-rob| method for examples with repeated poles. Moreover, our \verb|Schur-multi| method can still produce fairly good results when \verb|place| and \verb|robpole| fail for examples where the  multiplicity of the repeated pole is greater than $m$.


\bibliographystyle{IEEEtran}
\bibliography{IEEEabrv,mybibfile}

\begin{thebibliography}{19}

\bibitem {AB}
{J.~Abels and P.~Benner}, {``CAREX - A collection of benchmark examples
for continuous-time algebraic Riccati equations (Version 2.0)"},
\hskip 1em plus  0.5em minus 0.4em\relax
{Katholieke Universiteit Leuven, ESAT/SISTA, Leuven, Belgium,
SLICOT Working Note 1999-14, Nov. 1999. [Online]. Available:}
\url {http://www.slicot.de/REPORTS/SLWN1999-14.ps.gz.}

\bibitem {AB2}
{J.~Abels and P.~Benner}, {``DAREX - A collection of benchmark examples
for discrete-time algebraic Riccati equations (Version 2.0)"},
\hskip 1em plus 0.5em minus 0.4em\relax
{Katholieke Universiteit Leuven, ESAT/SISTA, Leuven, Belgium,
SLICOT Working Note 1999-16, Dec. 1999. [Online]. Available:}
\url {http://www.slicot.de/REPORTS/SLWN1999-16.ps.gz.}


\bibitem{BD}
{S.P.~Bhattacharyya and E.~De Souza}, {``Pole assignment via Sylvester's equation"},
{\it Systems \verb'&'  Control Letters}, 1(1982), 261--263.


\bibitem{BN}
{R.~Byers and S.G.~Nash}, {``Approaches to robust pole assignment"},
{\it International Journal of Control}, 49(1989), 97--117.


\bibitem{CB}
{R.K.~Cavin and S.P.~Bhattacharyya},
{``Robust and well-conditioned eigenstructure assignment via Sylvester's equation"},
{\it Optimal Control Applications and Methods}, 4(1983), 205--212.


\bibitem{Chu2}
{E.K.W.~Chu}, {``Pole assignment via the Schur form"},
{\it Systems  \verb'&' Control Letters}, 56(2007), 303--314.

\bibitem{Dic}
{A.~Dickman}, {``On the robustness of multivariable linear feedback systems in state-space representation"},
{\it IEEE Transactions on Automatic Control}, 32(1987), 407--410.


\bibitem{FO}
{M.~Fahmy and J.~O'Reilly}, {``On eigenstructure assignment in linear multivariable systems"},
{\it IEEE Transactions on   Automatic Control},  27(1982), 690--693.

\bibitem{FR}
{J.L.~Figueroa and J.A.~Romagnoli}, {``An algorithm for robust pole assignment via polynomial approach"},
{\it IEEE Transactions on   Automatic Control}, 39(1994),  831--835.

\bibitem{GR}
{V.~Gourishankar and  K.~Ramar}, {``Pole assignment with minimum eigenvalue sensitivity to plant parameter variations"},
{\it International Journal of Control}, 23(1976), 493--504.

\bibitem{GCQX}
{Z.C.~Guo, Y.F.~Cai, J.~Qian and S.F.~Xu},
{``A modified Schur method for robust pole assignment in state feedback control"},
{\it Automatica}, 52(2015), 334--339.


\bibitem{Ka}
{S.K.~Katti}, {``Pole placement in multi-input systems via elementary transformations"},
{\it International Journal of Control}, 37(1983), 315--347.


\bibitem{KN}
{J.~Kautsky and N.K.~Nichols}, {``Robust pole assignment in systems subject to structured perturbations"},
{\it Systems \verb'&' Control Letters}, 15(1990), 373--380.

\bibitem{KNV}
{J.~Kautsky, N.K.~Nichols and P.~Van Dooren}, {``Robust pole assignment in linear state feedback"},
{\it International Journal of Control}, 41(1985), 1129--1155.

\bibitem{LV}
{J.~Lam and W.Y.~Van}, {``A gradient flow approach to the robust pole-placement problem"},
{\it International Journal of Robust and Nonlinear Control}, 5(1995), 175--185.

\bibitem{LW}
{X.~Le and J.~Wang}, {``Robust pole assignment for synthesizing feedback control systems using recurrent neural networks"},
{\it IEEE Transactions on Neural Networks and Learning Systems}, 25(2014), 383--393.

%
%
%
%

\bibitem{MP2}
{G.S.~Miminis and  C.C.~Paige}, {``A direct algorithm for pole assignment of
time-invariant multi-input linear systems using state feedback"},
{\it Automatica}, 24(1988), 343--356.


\bibitem{MP3}
{G.S.~Miminis and  C.C.~Paige}, {``A QR-like approach for the eigenvalue assignment problem"},
{in} {\it Proceedings of the 2nd Hellenic Conference on Mathematics and Informatics},
Athens, Greece, Sep.,  1994.


\bibitem{PM}
{R.V.~Patel and P.~Misra}, {``Numerical algorithms for eigenvalue assignment by state feedback"},
{\it Proceedings of the IEEE}, 72(1984), 1755--1764.


\bibitem{PCK1}
{P.Hr.~Petkov, N.D.~Christov and M.M.~Konstantinov},
{``A computational algorithm for pole assignment of linear multiinput systems"},
{\it IEEE Transactions on Automatic Control}, 31(1986), 1044--1047.

\bibitem{RG}
{K.~Ramar and V.~Gourishankar},
{``Utilization of the design freedom of pole assignment feedback controllers of unrestricted rank"},
{\it International Journal of Control}, 24(1976), 423--430.

\bibitem{RFBT}
{M.A.~Rami, S.E.~Faiz, A.~Benzaouia and  F.~Tadeo}, {``Robust exact pole placement via an LMI-based algorithm"},
{\it IEEE Transactions on Automatic Control}, 54(2009), 394--398.

\bibitem{RM}
{D.G.~Retallack and A.G.J.~Macfarlane}, {``Pole-shifting techniques for multivariable feedback systems"},
{\it Proceedings of the Institution of Electrical Engineers}, 117(1970), 1037--1038.

\bibitem{SNNP}
{R.~Schmid, L.~Ntogramatzidis, T.~Nguyen and  A.~Pandey},
{``A unified method for optimal arbitrary pole placement"},
{\it Automatica}, 50(2014),  2150--2154.


\bibitem{SAY}
{V.~Sima, A.L.~Tits and Y.~Yang}, {``Computational experience with robust pole assignment algorithms"},
{in} {\it Computer Aided Control System Design, 2006 IEEE International Conference on Control Applications,
2006 IEEE International Symposium on Intelligent Control, 2006 IEEE},
Munich, Germany,  4-6 Oct., 2006.

\bibitem{SEPB}
{Y.C.~Soh, R.J.~Evans, I.R.~Petersen and R.E.~Betz}, {``Robust pole assignment"},
{\it Automatica}, 23(1987), 601-610.


\bibitem{SSun}
{G.W.~Stewart and J.G.~Sun}, {\it Matrix Perturbation Theory},
\hskip 1em plus 0.5em minus 0.4em\relax
{Academic Press}, New York, 1990.



\bibitem{Tits}
{A.L.~Tits and Y.~Yang}, {``Globally convergent algorithms for robust pole assignment by state feedback"},
{\it IEEE Transactions on Automatic Control}, 41(1996), 1432--1452.

\bibitem{Var}
{A.~Varga}, {``A Schur method for pole assignment"},
{\it IEEE Transactions on Automatic Control}, 26(1981), 517--519.

\bibitem{Var1}
{A.~Varga}, {``Robust pole assignment via Sylvester equation based state feedback parametrization"}
{in} {\it Computer-Aided Control System Design, 2000. CACSD 2000. IEEE International Symposium on}
Anchorage, AK, Sep., 2000.


\bibitem{Wonh}
{W.M.~Wonham}, {\it Linear Multivariable Control: A Geometric Approach}, {3rd ed.},
\hskip 1em plus 0.5em minus 0.4em\relax
{Springer-Verlag}, New York, 1985.

\bibitem{XU}
{S.F.~Xu}, {\it An Introduction to Inverse Algebraic Eigenvalue Problems},
\hskip 1em plus 0.5em minus 0.4em\relax
{Peking University Press, Beijing, and Vieweg, Braunschweig}, 1998.

\end{thebibliography}

\end{document}